\documentclass{amsart}%
\usepackage{amssymb}
\usepackage{amsmath}
\usepackage{graphicx}
\usepackage{amsfonts}%
\newtheorem{theorem}{Theorem}
\theoremstyle{plain}

\newtheorem{corollary}{Corollary}

\newtheorem{definition}{Definition}
\newtheorem{example}{Example}

\newtheorem{lemma}{Lemma}

\newtheorem{proposition}{Proposition}
\newtheorem{remark}{Remark}

\numberwithin{equation}{section}

\newcommand{\M}{\mathcal{M}}

\begin{document}
\title{Diagonals on the Permutahedra, Multiplihedra and Associahedra }
\author{Samson Saneblidze$^{1}$}
\address{A. Razmadze Mathematical Institute\\
Georgian Academy of Sciences\\
M. Aleksidze st., 1\\
0193 Tbilisi, Georgia}
\email{sane@rmi.acnet.ge <mailto:sane@rmi.acnet.ge> }
\author{Ronald Umble$^{2}$}
\address{Department of Mathematics\\
Millersville University of Pennsylvania\\
Millersville, PA. 17551}
\email{Ron.Umble@Millersville.edu {<}mailto:Ron.Umble@Millersville.edu {>}}
\thanks{$^{1}$ This research was funded in part by Award No. GM1-2083 of the U.S.
Civilian Research and Development Foundation for the Independent States of the
Former Soviet Union (CRDF) and by Award No. 99-00817 of INTAS}
\thanks{$^{2}$ This research was funded in part by a Millersville University faculty
research grant.}
\date{June 4, 2004}

\begin{abstract}
We construct an explicit diagonal $\Delta_{P}$ on the permutahedra $P.$
Related diagonals on the multiplihedra $J$ and the associahedra $K$ are
induced by Tonks' projection $P\rightarrow K$ \cite{tonks} and its
factorization through $J.$ We introduce the notion of a permutahedral set
$\mathcal{Z}$ and lift $\Delta_{P}$ to a diagonal on $\mathcal{Z}$. We show
that the double cobar construction $\Omega^{2}C_{\ast}(X)$ is a permutahedral
set; consequently $\Delta_{P}$ lifts to a diagonal on $\Omega^{2}C_{\ast}(X)$.
Finally, we apply the diagonal on $K$ to define the tensor product of
$A_{\infty}$-(co)algebras in maximal generality.

\end{abstract}
\maketitle

\section{Introduction}

A permutahedral set is a combinatorial object generated by permutahedra $P$
and equipped with appropriate face and degeneracy operators. Permutahedral
sets are distinguished from cubical or simplicial set by higher order
(non-quadratic) relations among face and degeneracy operators. In this paper
we define the notion of a permutahedral set and observe that the double cobar
construction $\Omega^{2}C_{\ast}\left(  X\right)  $ is a naturally occurring
example.\ We construct an explicit diagonal $\Delta_{P}:C_{\ast}(P)\rightarrow
C_{\ast}(P)\otimes C_{\ast}(P)$ on the cellular chains of permutahedra and
show how to lift $\Delta_{P}$ to a diagonal on any permutahedral set. We
factor Tonks' projection $\theta:P\rightarrow K$ through the multiplihedron
$J$ and obtain diagonals $\Delta_{J}$ on $C_{\ast}(J)$ and $\Delta_{K}$ on
$C_{\ast}(K).$ We apply $\Delta_{K}$ to define the tensor product of
$A_{\infty}$-(co)algebras in maximal generality; this resolves a long-standing
problem in the theory of operads. Gaberdiel and Zwiebach's open string field
theory \cite{Gaberdiel-Zwiebach} provides a setting in which this tensor
product can be applied.

The paper is organized as follows: Sections 2 and 5 review the families of
polytopes we consider. The diagonal $\Delta_{P}$ is defined in Section 3 and
lifted to general permutahedral sets in Section 4. The related diagonals
$\Delta_{J}$ and $\Delta_{K}$ are obtained in Section 6 and applied in Section
7 to define the tensor product of $A_{\infty}$-(co)algebras in maximal
generality. Sections 5 through 7 do not depend on Section 4.

The first author wishes to acknowledge conversations with Jean-Louis Loday
from which our representation of the permutahedron as a subdivision of the
cube emerged. The second author wishes to thank Millersville University for
its generous financial support and the University of North Carolina at Chapel
Hill for its kind hospitality during work on parts of this project.

\section{The Permutahedra}

Let $S_{n}$ be the symmetric group on $\underline{n}=\left\{  1,2,\ldots
,n\right\}  .$ Recall that the permutahedron $P_{n}$ is the convex hull of
$n!$ vertices $\left(  \sigma(1),\ldots,\sigma(n)\right)  \in\mathbb{R}^{n},$
$\sigma\in S_{n}$ \cite{Coxeter}, \cite{Milgram}, \cite{Zigler}.$\ $As a
cellular complex, $P_{n}$ is an $\left(  n-1\right)  $-dimensional convex
polytope whose $\left(  n-p\right)  $-faces are indexed by (ordered)
partitions $U_{1}|\cdots|U_{p}$ of $\underline{n}$. We shall define the
permutahedra inductively as subdivisions of the standard $n$-cube $I^{n}.$
With this representation the combinatorial connection between faces and
partitions is immediately clear.

Assign the label $\underline{1}$ to the single point $P_{1}.$ If $P_{n-1}$ has
been constructed and $u=U_{1}|\cdots|U_{p}$ is one of its faces, form the
sequence $u_{\ast}=\left\{  u_{0}=0,u_{1},\ldots,u_{p-1},u_{p}=\infty\right\}
$ where $u_{j}=\#\left(  U_{p-j+1}\cup\cdots\cup U_{p}\right)  ,$ $1\leq j\leq
p-1$ and $\#$ denotes cardinality. Define the \emph{subdivision of }%
$I$\emph{\ relative to }$u$ to be
\[
I/u_{\ast}=I_{1}\cup I_{2}\cup\cdots\cup I_{p},
\]
where $I_{j}=\left[  1-\frac{1}{2^{u_{j-1}}},1-\frac{1}{2^{u_{j}}}\right]  $
and $\frac{1}{2^{\infty}}=0.$ Then
\[
P_{n}=\bigcup\limits_{u\in P_{n-1}}u\times I/u_{\ast}%
\]
with faces labeled as follows (see Figures 1 and 2)\vspace{0.1in}:%
\[%
\begin{tabular}
[c]{c|cc}%
\textbf{Face of }$\underset{\ }{u\times I/u_{\ast}}$ & \textbf{Partition of
}$\underline{n}$ & \\\hline
&  & \\
$u\times0$ & $U_{1}|\cdots|U_{p}|n$ & \\
&  & \\
$u\times(I_{j}\cap I_{j+1})$ & $U_{1}|\cdots|U_{p-j}|n|U_{p-j+1}|\cdots
|U_{p},$ & $1\leq j\leq p-1$\\
&  & \\
$u\times1$ & $n|U_{1}|\cdots|U_{p}$ & \\
&  & \\
$u\times I_{j}$ & $U_{1}|\cdots|U_{p-j+1}\cup n|\cdots|U_{p}.$ &
\end{tabular}
\
\]

A \emph{cubical }vertex of $P_{n}$ is a vertex common to both $P_{n}$ and
$I^{n-1}.$ Note that $u$ is a cubical\ vertex of $P_{n-1}$ if and only if
$u|n$ and $n|u$ are cubical\ vertices of $P_{n}.$ Thus the cubical vertices of
$P_{3}$ are $1|2|3,$ $2|1|3,$ $3|1|2$ and $3|2|1$ since $1|2$ and $2|1$ are
cubical vertices of $P_{2}.$\vspace{0.1in}

\begin{center}
\setlength{\unitlength}{0.0004in}\begin{picture}
(2975,2685)(3126,-2038) \thicklines \put(3601,239){\line(
1,0){1800}} \put(5401,239){\line( 0,-1){1800}}
\put(5401,-1561){\line(-1, 0){1800}} \put(3601,-1561){\line(
0,1){1800}} \put(3601,239){\makebox(0,0){$\bullet$}}
\put(3601,-661){\makebox(0,0){$\bullet$}}
\put(3601,-1561){\makebox(0,0){$\bullet$}}
\put(5401,239){\makebox(0,0){$\bullet$}}
\put(5401,-661){\makebox(0,0){$\bullet$}}
\put(5401,-1561){\makebox(0,0){$\bullet$}}
\put(4500,-680){\makebox(0,0){$123$}}
\put(2980,-1861){\makebox(0,0){$1|2|3$}}
\put(2980,-699){\makebox(0,0){$1|3|2$}}
\put(2980,464){\makebox(0,0){$3|1|2$}}
\put(6000,-1861){\makebox(0,0){$2|1|3$}}
\put(6000,-699){\makebox(0,0){$2|3|1$}}
\put(6000,464){\makebox(0,0){$3|2|1$}}
\put(3040,-1260){\makebox(0,0){$1|23$}}
\put(4550,530){\makebox(0,0){$3|12$}}
\put(3040,-111){\makebox(0,0){$13|2$}}
\put(5960,-111){\makebox(0,0){$23|1$}}
\put(5960,-1260){\makebox(0,0){$2|13$}}
\put(4550,-1890){\makebox(0,0){$12|3$}}
\end{picture}\vspace{0.1in}

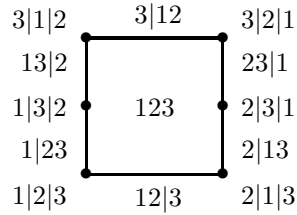
Figure 1: $P_{3}$ as a subdivision of $P_{2}\times I$.\ \vspace{0.4in}

\setlength{\unitlength}{0.00023in}\begin{picture}
(7500,7500) \thicklines
\put(3000,4800){\line( 0,-1){4800}}
\put(3000,4800){\makebox(0,0){$\bullet$}}
\put(3000,2400){\makebox(0,0){$\bullet$}} \put(3000,0){\makebox
(0,0){$\bullet$}} \put(3000,3600){\makebox(0,0){$\bullet$}}
\put(3000,0){\line( 1, 0){4800}}
\put(7800,0){\line( 0, 1){4800}}
\put(7800,4800){\line(-1, 0){4800}}
\put(7800,4800){\makebox(0,0){$\bullet$}}
\put(7800,3600){\makebox(0,0){$\bullet$}} \put(7800,2400){\makebox
(0,0){$\bullet$}} \put(7800,0){\makebox(0,0){$\bullet$}}
\put(3000,2400){\line( 1, 0){4800}}
\put(0,6800){\line( 0,-1){4800}}
\put(0,6800){\makebox(0,0){$\bullet$}}
\put(0,5600){\makebox(0,0){$\bullet$}}
\put(0,4400){\makebox(0,0){$\bullet$}}
\put(0,2000){\makebox(0,0){$\bullet$}}
\put(0,6800){\line( 1, 0){4800}}
\put(0,5600){\line( 1, 0){1200}} \put(2000,5600){\line( 1,
0){2800}}
\put(0,2000){\line( 1, 0){1200}} \put(1800,2000){\line( 1,
0){900}} \put(3300,2000){\line( 1, 0){1500}}
\put(4800,5000){\line( 0,1){1700}} \put(4800,2600){\line(
0,1){2050}} \put(4800,2000){\line( 0,1){200}}
\put(4800,6800){\makebox(0,0){$\bullet$}}
\put(4800,5600){\makebox(0,0){$\bullet$}} \put(4800,4400){\makebox
(0,0){$\bullet$}} \put(4800,2000){\makebox(0,0){$\bullet$}}
\put(0,2000){\line( 3,-2){3000}}
\put(3000,4800){\line(-3,2){3000}}
\put(1500,5800){\line( 0,-1){4800}}
\put(1500,5800){\makebox(0,0){$\bullet$}}
\put(1500,4600){\makebox(0,0){$\bullet$}} \put(1500,3400){\makebox
(0,0){$\bullet$}} \put(1500,1000){\makebox(0,0){$\bullet$}}
\put(3000,3600){\line(-3, 2){1500}}
\put(1500,3400){\line(-3, 2){1500}}
\put(6300,3400){\line(-3, 2){1500}}
\put(6300,5800){\line( 0,-1){800}} \put(6300,4600){\line(
0,-1){2000}} \put(6300,2180){\line( 0,-1){1200}}
\put(6300,5800){\makebox(0,0){$\bullet$ }}
\put(6300,3400){\makebox(0,0){$\bullet$}} \put(6300,1000){\makebox
(0,0){$\bullet$}} \put(6300,4600){\makebox(0,0){$\bullet$}}
\put(7800,3600){\line(-3, 2){1500}}
\put(4800,2000){\line(3, -2){3000}}
\put(4800,6800){\line( 3,-2){3000}}
\put(-1000,1000){\makebox(0,0){$(0,1,0)$}}
\put(5100,7600){\makebox(0,0){$(1,1,1)$}}
\put(8750,-1000){\makebox(0,0){$(1,0,0)$}} \put
(2700,-1000){\makebox(0,0){$(0,0,0)$}}
\end{picture}\vspace{0.4in}

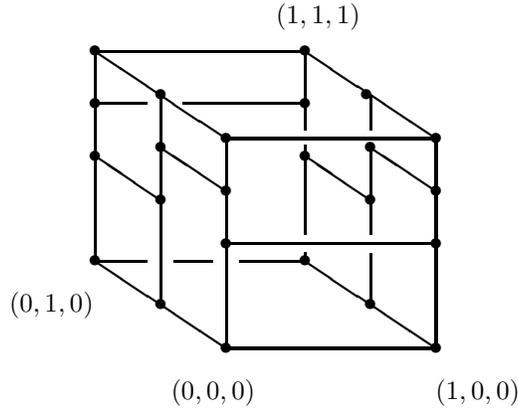
Figure 2a: $P_{4}$ as a subdivision of $P_{3}\times I.\vspace{0.2in}$

\setlength{\unitlength}{0.009in}\begin{picture}
(500,-500) \thicklines
\put(0,-120){\line( 0,-1){120}} \put(120,0){\line( 0,-1){360}}
\put(180,-120){\line( 0,-1){120}} \put(240,0){\line( 0,-1){360}}
\put(360,-120){\line( 0,-1){120}} \put(420,-120){\line(
0,-1){120}} \put(480,-120){\line( 0,-1){120}}
\put(120,0){\line( 1,0){120}} \put(0,-120){\line( 1,0){480}} \put
(120,-150){\line( 1,0){60}} \put(180,-180){\line( 1,0){60}} \put
(240,-150){\line( 1,0){120}} \put(420,-150){\line( 1,0){60}} \put
(0,-180){\line( 1,0){120}} \put(360,-180){\line( 1,0){60}}
\put(0,-240){\line( 1,0){480}} \put(120,-360){\line( 1,0){120}}
\put(120,0){\makebox(0,0){$\bullet$}}
\put(180,0){\makebox(0,0){$\bullet$}}
\put(240,0){\makebox(0,0){$\bullet$}}
\put(0,-120){\makebox(0,0){$\bullet$}}
\put(120,-120){\makebox(0,0){$\bullet$}} \put(180,-120){\makebox
(0,0){$\bullet$}} \put(240,-120){\makebox(0,0){$\bullet$}} \put
(360,-120){\makebox(0,0){$\bullet$}}
\put(420,-120){\makebox(0,0){$\bullet$}}
\put(480,-120){\makebox(0,0){$\bullet$}} \put(120,-150){\makebox
(0,0){$\bullet$}} \put(180,-150){\makebox(0,0){$\bullet$}} \put
(240,-150){\makebox(0,0){$\bullet$}}
\put(360,-150){\makebox(0,0){$\bullet$}}
\put(0,-180){\makebox(0,0){$\bullet$}}
\put(182.25,-180){\makebox(0,0){$\bullet$ }}
\put(240,-180){\makebox(0,0){$\bullet$}} \put(420,-180){\makebox
(0,0){$\bullet$}} \put(480,-180){\makebox(0,0){$\bullet$}} \put
(0,-150){\makebox(0,0){$\bullet$}}
\put(120,-240){\makebox(0,0){$\bullet$}}
\put(360,-240){\makebox(0,0){$\bullet$}} \put(420,-240){\makebox
(0,0){$\bullet$}} \put(480,-150){\makebox(0,0){$\bullet$}} \put
(0,-240){\makebox(0,0){$\bullet$}}
\put(120,-180){\makebox(0,0){$\bullet$}}
\put(180,-240){\makebox(0,0){$\bullet$}} \put(240,-240){\makebox
(0,0){$\bullet$}} \put(360,-180){\makebox(0,0){$\bullet$}} \put
(420,-150){\makebox(0,0){$\bullet$}}
\put(480,-240){\makebox(0,0){$\bullet$}}
\put(120,-360){\makebox(0,0){$\bullet$}} \put(180,-360){\makebox
(0,0){$\bullet$}} \put(240,-360){\makebox(0,0){$\bullet$}}
\put(275,-100){\makebox(0,0){$(1,1,1)$}}
\put(20,-260){\makebox(0,0){$(0,0,0)$ }}
\put(55,-152){\makebox(0,0){$124|3$}}
\put(55,-210){\makebox(0,0){$12|34$}}
\put(150,-134){\makebox(0,0){$24|13$}}
\put(150,-195){\makebox(0,0){$2|134$}}
\put(177,-55){\makebox(0,0){$4|123$}}
\put(177,-295){\makebox(0,0){$123|4$}}
\put(210,-152){\makebox(0,0){$234|1$}}
\put(210,-210){\makebox(0,0){$23|14$}}
\put(300,-134){\makebox(0,0){$34|12$}}
\put(300,-195){\makebox(0,0){$3|124$}}
\put(390,-152){\makebox(0,0){$134|2$}}
\put(390,-200){\makebox(0,0){$13|24$}}
\put(450,-134){\makebox(0,0){$14|23$}}
\put(450,-194){\makebox(0,0){$1|234$}}
\end{picture}\vspace*{3.4in}

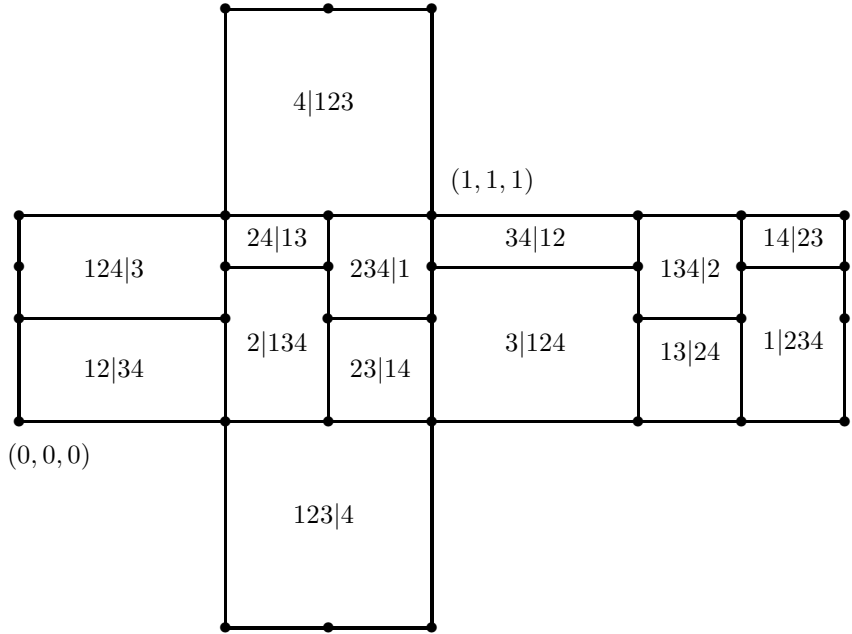
Figure 2b: The $2$-faces of $P_{4}.$
\end{center}

\section{A Diagonal on the Permutahedra}

In this section we construct a combinatorial diagonal on the cellular chains
of the permutahedron $P_{n+1}.$ Given a polytope $X,$ let $\left(  C_{\ast
}\left(  X\right)  ,\partial\right)  $ denote the cellular chains on $X$ with
boundary $\partial.$

\begin{definition}
\label{defn1}A map $\Delta_{X}:C_{\ast}(X)\rightarrow C_{\ast}(X)\otimes
C_{\ast}(X)$ is a \underline{diagonal} on $C_{\ast}\left(  X\right)  $ if
\end{definition}

\begin{enumerate}
\item $\Delta_{X}\left(  C_{\ast}(e)\right)  \subseteq C_{\ast}(e)\otimes
C_{\ast}(e)$\textit{\ for each cell }$e\subseteq X$\textit{\ and}

\item $\left(  C_{\ast}\left(  X\right)  ,\Delta_{X},\partial\right)
$\textit{\ is a DG coalgebra.}
\end{enumerate}

\noindent In general, the DG coalgebra $\left(  C_{\ast}\left(  X\right)
,\Delta_{X},\partial\right)  $ is non-coassociative, non-cocommuta- tive and
non-counital; thus the statement (2) in Definition \ref{defn1} is equivalent
to stating that $\Delta_{X}$ is a chain map. We remark that a diagonal
$\Delta_{P}$ on $C_{\ast}\left(  P_{n+1}\right)  $ is unique if the following
two additional properties hold:\textit{\vspace{0.1in}\ }

\begin{enumerate}
\item \textit{The canonical cellular projection }$\rho_{n+1}:P_{n+1}%
\rightarrow I^{n}$ \textit{\ induces a DG coalgebra map }$C_{\ast}\left(
P_{n+1}\right)  \rightarrow C_{\ast}\left(  I^{n}\right)  $ (see Section 4,
Figures 3 and 4) \textit{and\vspace{0.1in}}

\item \textit{There is a minimal number of} \textit{components} $a\otimes b$
\textit{in} $\Delta_{P}\left(  C_{k}\left(  P_{n+1}\right)  \right)  $
\textit{for }$0\leq k\leq n$.\textit{\vspace*{0.1in}}
\end{enumerate}

\noindent Since the uniqueness of $\Delta_{P}$ is not used in our work,
verification of these facts is left to the interested reader.

\begin{definition}
A partition $A_{1}|\cdots|A_{p}$ is \underline{step increasing} iff
$A_{p}|\cdots|A_{1}$ is \underline{step decreas-} \underline{ing} iff $\min
A_{j}<\max A_{j+1}$ for all $j\leq p-1.$ A \underline{step} partition is
either step increasing or step decreasing.
\end{definition}

Think of $\sigma\in S_{p+q-1}$ as an ordered sequence of positive integers;
let $\overleftarrow{\sigma}_{j}$ and $\overrightarrow{\sigma}_{q-i+1}$ denote
its $j^{th}$ decreasing and $i^{th}$ increasing subsequence of maximal length.
Then $\overleftarrow{\sigma}_{1}|\cdots|\overleftarrow{\sigma}_{p}$ and
$\overrightarrow{\sigma}_{q}|\cdots|\overrightarrow{\sigma}_{1}$ are step
increasing and step decreasing partitions of $\underline{p+q-1},$ respectively
(see Example \ref{example2} below).

\begin{definition}
A pairing of partitions $A_{1}|\cdots|A_{p}\otimes B_{q}|\cdots|B_{1}$ is a
\underline{strong complemen-} \underline{tary pair}(SCP) if there exists
$\sigma\in S_{p+q-1}$ such that $A_{j}=\overleftarrow{\sigma}_{j}$ and
$B_{i}=\overrightarrow{\sigma}_{i}$ as unordered sets for all $i,j.$
\end{definition}

SCP's have a natural matrix representation.

\begin{definition}
A $q\times p$ matrix $O=\left(  o_{ij}\right)  $ is \underline{ordered} if:

\begin{enumerate}
\item $\left\{  o_{i,j}\right\}  =\left\{  0,1,\ldots,p+q-1\right\}  ;$

\item Each row and column of $O$ is non-zero;

\item Non-zero entries in $O$ are distinct and increase in each row and column.
\end{enumerate}
\end{definition}

Let $\mathcal{O}$ denote the set of ordered matrices. Note that the rows and
columns of an ordered matrix $O^{q\times p}$ form a partition of
$\underline{p+q-1}.$

\begin{definition}
Given $O\in\mathcal{O}^{q\times p},$ let $V_{i}=\,$row$_{i}\left(  O\right)
\cap\mathbb{Z}^{+}$ for $i\leq q$ and $U_{j}=\,$col$_{j}\left(  O\right)
\cap\mathbb{Z}^{+}$ for $j\leq p.$ The \underline{row face of $O$} is the face
$r\left(  O\right)  =V_{q}|\cdots|V_{1}\subset P_{p+q-1};$ the \underline
{column face of $O$} is the face $c\left(  O\right)  =U_{1}|\cdots
|U_{p}\subset P_{p+q-1}.$
\end{definition}

\begin{definition}
\label{step}An ordered matrix $E$ is a \underline{step matrix} if:

\begin{enumerate}
\item Non-zero entries in each row of $E$ appear in consecutive columns;

\item Non-zero entries in each column of $E$ appear in consecutive rows;

\item The sub, main and super diagonals of $E$ contain a single non-zero entry.
\end{enumerate}
\end{definition}

Let $\mathcal{E}$ denote the set of step matrices. If $E=\left(
e_{i,j}\right)  \in\mathcal{E}^{q\times p}$, condition (1) in Definition
\ref{step} groups the non-zero entries in each row together in a horizontal
block, condition (2) groups the non-zero entries in each column together in a
vertical block and condition (3) links horizontal and vertical blocks to
produce a \textquotedblleft staircase path\textquotedblright\ of non-zero
entries connecting the lower-left and upper-right entries $e_{q,1}$ and
$e_{1,p}$\ (see Example \ref{example2} below). Clearly, $c\left(  E\right)
\otimes r\left(  E\right)  =\overleftarrow{\sigma}_{1}|\cdots|\overleftarrow
{\sigma}_{p}\otimes\overrightarrow{\sigma}_{q}|\cdots|\overrightarrow{\sigma
}_{1}$ for some $\sigma\in S_{p+q-1},$ so $c\left(  E\right)  \otimes r\left(
E\right)  $ is an SCP. Furthermore, one can recover $E$ from $\sigma=\left(
x_{1}\text{ }x_{2}\text{ }\cdots\text{ }x_{n+1}\right)  $ in the following
way: Set $e_{q,1}=x_{1}.$ Inductively, assume $e_{i,j}=x_{k-1};$ if
$x_{k-1}<x_{k},$ set $e_{i,j+1}=x_{k};$ otherwise, set $e_{i-1,j}=x_{k}.$ Let
$E_{\sigma}$ denote the step matrix given by $\sigma\in S=\lim
\limits_{\rightarrow}S_{n+1}$. We have proved:

\begin{proposition}
\label{ss}There exist one-to-one correspondences%
\[%
\begin{tabular}
[c]{c}%
$\mathcal{E}\leftrightarrow\text{$S$}\leftrightarrow\left\{  \text{Step
increasing partitions}\right\}  \leftrightarrow\left\{  \text{Step decreasing
partitions}\right\}  \leftrightarrow\left\{  \text{SCP's}\right\}  $\\
\\
$E_{\sigma}$ $\ \ \leftrightarrow$ $\ \ \sigma$ $\ \leftrightarrow$
$\ \ \overleftarrow{\sigma}_{1}|\cdots|\overleftarrow{\sigma}_{p}$
$\ \leftrightarrow$ $\ \overrightarrow{\sigma}_{q}|\cdots|\overrightarrow
{\sigma}_{1}$ $\ \leftrightarrow$\ $\ \ \overleftarrow{\sigma}_{1}%
|\cdots|\overleftarrow{\sigma}_{p}\otimes\overrightarrow{\sigma}_{q}%
|\cdots|\overrightarrow{\sigma}_{1}.$%
\end{tabular}
\]

\end{proposition}

\begin{example}
\label{example2}The permutation
\[
\sigma=\left(  9\text{ }7\text{ }1\text{ }3\text{ }8\text{ }4\text{ }6\text{
}5\text{ }2\right)
\]
corresponds to step matrix%
\[
E_{\sigma}=%
\begin{tabular}
[c]{|l|l|l|l|}\hline
&  &  & 2\\\hline
&  &  & 5\\\hline
&  & 4 & 6\\\hline
1 & 3 & 8 & \\\hline
7 &  &  & \\\hline
9 &  &  & \\\hline
\end{tabular}
\text{ .}%
\]
and the SCP%
\[
c\left(  E_{\sigma}\right)  \otimes r\left(  E_{\sigma}\right)
=971|3|84|652\otimes9|7|138|46|5|2.
\]

\end{example}

We now introduce matrix transformations that operate like the vertical and
horizontal shifts one performs in a tableau puzzle. For $\left(  i,j\right)
\in\mathbb{Z}^{+}\times\mathbb{Z}^{+},$ define the \emph{down-shift} and
\emph{right-shift} operators $D_{i,j},R_{i,j}:\mathcal{O}\rightarrow
\mathcal{O}$ on $O^{q\times p}=\left(  o_{i,j}\right)  $ by$\smallskip$

\begin{enumerate}
\item $D_{i,j}O=O$ unless $i\leq q-1,$ $o_{i+1,j}=0,$ $o_{i,j}o_{i,k}>0$ for
some $k\neq j,$ $o_{i,j}>o_{i+1,\ell}$ for $1\leq\ell<j,$ and $o_{i+1,\ell
}>o_{i,j}$ whenever $o_{i+1,\ell}>0$ and $j<\ell\leq p,$ in which case
$D_{i,j}O$ is obtained from $O$ by transposing $o_{i,j}$ and $o_{i+1,j}%
;\smallskip$

\item $R_{i,j}O=O$ unless $j\leq p-1,$ $o_{i,j+1}=0,$ $o_{i,j}o_{k,j}>0$ for
some $k\neq i,$ $o_{i,j}>o_{\ell,j+1}$ for $1\leq\ell<i,$ and $o_{\ell
,j+1}>o_{i,j}$ whenever $o_{\ell,j+1}>0$ and $j<\ell\leq q,$ in which case
$R_{i,j}O$ is obtained from $O$ by transposing $o_{i,j}$ and $o_{i,j+1}%
.\smallskip$
\end{enumerate}

\begin{definition}
\label{config}A matrix $F\in\mathcal{O}$ is a \underline{configuration matrix}
if there is a step matrix $E$ and a sequence of shift operators $G_{1}%
,\ldots,G_{m}$ such that

\begin{enumerate}
\item $F=G_{m}\cdots G_{1}E;$

\item If $G_{m}\cdots G_{1}=\cdots D_{i_{2},j_{2}}\cdots D_{i_{1},j_{1}}%
\cdots,$ then $i_{1}\leq i_{2};$

\item If $G_{m}\cdots G_{1}=\cdots R_{k_{2},\ell_{2}}\cdots R_{k_{1},\ell_{1}%
}\cdots,$ then $\ell_{1}\leq\ell_{2}.\smallskip$\newline When this occurs, we
say that $F$ is \underline{derived from} $E$ and refer to the pairing
$c\left(  F\right)  \otimes r\left(  F\right)  $ as a \underline{complementary
pair (CP) related to} $c\left(  E\right)  \otimes r\left(  E\right)  $.
\end{enumerate}
\end{definition}

Let $\mathcal{C}$ denote the set of configuration matrices. For $F=\left(
f_{i,j}\right)  \in\mathcal{C}$ with column face $U_{1}|\cdots|U_{p}$ and row
face $V_{q}|\cdots|V_{1},$ choose proper subsets $N_{i}=\left\{  f_{i,n_{1}%
}<\cdots\right.  $ $\left.  <f_{i,n_{k}}\right\vert \left.  \max
V_{i+1}<f_{i,n_{1}}\right\}  \subset V_{i}$ and $M_{j}=\left\{  f_{m_{1}%
,j}<\cdots<f_{m_{\ell},j}\left\vert \max U_{j+1}<f_{m_{1},j}\right.  \right\}
$ $\subset U_{j}$ and define
\[
D_{N_{i}}^{i}F=D_{i,n_{k}}\cdots D_{i,n_{1}}F\text{ \ and \ }R_{M_{j}}%
^{j}F=R_{m_{\ell},j}\cdots R_{m_{1},j}F.
\]
We often suppress the superscript when it is clear from context. The fact that
$D_{i,j+1}R_{i,j}F=R_{i+1,j}D_{i,j}F$ wherever both maps in the composition
act non-trivially, gives the following useful reformulation of Definition
\ref{config}:

\begin{proposition}
A matrix $F\in\mathcal{O}$ with $c\left(  F\right)  =U_{1}|\cdots|U_{p}$ and
$r\left(  F\right)  =V_{q}|\cdots|V_{1}$ is a configuration matrix if and only
if there exists $E\in\mathcal{E}$ and proper subsets $M_{j}\subset U_{j}$ and
$N_{i}\subset V_{i}$ such that
\[
F=D_{N_{q-1}}\cdots D_{N_{1}}R_{M_{p-1}}\cdots R_{M_{1}}E.
\]

\end{proposition}

\begin{example}
\label{ex4}Four configuration matrices $F$ can be derived from the step
matrix\vspace{0.1in}
\[%
\begin{tabular}
[c]{rlll}%
$E=$ & $%
\begin{tabular}
[c]{|l|l|l|}\hline
& 2 & 3\\\hline
1 & 5 & \\\hline
4 &  & \\\hline
\end{tabular}
\ \ \ \ \ :$ &  & \\
&  &  & \\
$D_{\varnothing}D_{\varnothing}R_{\varnothing}R_{\varnothing}E=$ & $%
\begin{tabular}
[c]{|l|l|l|}\hline
& 2 & 3\\\hline
1 & 5 & \\\hline
4 &  & \\\hline
\end{tabular}
\ \ \ \ \ $ & $\leftrightarrow$ & $14|25|3\otimes4|15|23,$\\
&  &  & \\
$D_{\varnothing}D_{\varnothing}R_{5}R_{\varnothing}E=$ & $%
\begin{tabular}
[c]{|l|l|l|}\hline
& 2 & 3\\\hline
1 &  & 5\\\hline
4 &  & \\\hline
\end{tabular}
\ \ \ \ \ $ & $\leftrightarrow$ & $14|2|35\otimes4|15|23,$\\
&  &  & \\
$D_{5}D_{\varnothing}R_{\varnothing}R_{\varnothing}E=$ & $%
\begin{tabular}
[c]{|l|l|l|}\hline
& 2 & 3\\\hline
1 &  & \\\hline
4 & 5 & \\\hline
\end{tabular}
\ \ \ \ \ $ & $\leftrightarrow$ & $14|25|3\otimes45|1|23,$\\
&  &  & \\
$D_{5}D_{\varnothing}R_{5}R_{\varnothing}E=$ & $%
\begin{tabular}
[c]{|l|l|l|}\hline
& 2 & 3\\\hline
1 &  & \\\hline
4 &  & 5\\\hline
\end{tabular}
\ \ \ \ \ $ & $\leftrightarrow$ & $14|2|35\otimes45|1|23$.
\end{tabular}
\ \ \ \ \ \ \ \
\]
Up to sign, the CP's
\[
c\left(  F\right)  \otimes r\left(  F\right)  =(14|2|35+14|25|3)\otimes
(4|15|23+45|1|23)
\]
are components of $\Delta_{P}\left(  \underline{5}\right)  $.
\end{example}

Let us associate formal \textquotedblleft configuration
signs\textquotedblright\ to configuration matrices. The signs we introduce
here can be derived by induction on dimension given that $P_{2}=I$ and
$\Delta_{P}$ is a chain map. \textit{Henceforth we assume that all blocks in a
partition are increasingly ordered. }First note that a face $u=U_{1}%
|\cdots|U_{p}\subset P_{n+1}$ is an $\left(  n-p+1\right)  $-face of $p-1$
faces in dimension $n-p+2.$ Thus there are $\left(  p-1\right)  !$ ways to
produce $u$ by successively inserting bars into $\underline{n+1},$ each of
which has an associated sign. Of these, we need the right-most and left-most
insertion procedures.

When each $x\in\underline{n+1}$ has degree $1,$ the sign of a permutation
$\sigma\in S_{n+1}$ is the Koszul sign that arises from the action of
$\sigma.$ Thus, if $\sigma$ transposes adjacent subsets $U,V\subset
\underline{n+1}$ for example, then $sgn\left(  \sigma\right)  =\left(
-1\right)  ^{\#U\#V}.$ For $u=U_{1}|\cdots|U_{p}\subset P_{n+1},$ denote the
sign of the permutation $\underline{n+1}\rightarrow U_{1}\cup\cdots\cup U_{p}$
by $psgn\left(  u\right)  $; note that $\sigma$ is an unshuffle of
$\underline{n}$ when $p=2,$ in which case we denote $psgn\left(  u\right)  =$
\textit{shuff}$\left(  U_{1};U_{2}\right)  $. Let $m_{i}=\#U_{i}-1$ and
identify $u$ with the Cartesian product $P_{m_{1}+1}\times\cdots\times
P_{m_{p}+1};$ then
\[
C_{n-p+1}\left(  u\right)  =C_{m_{1}}\left(  U_{1}\right)  \otimes
\cdots\otimes C_{m_{p}}\left(  U_{p}\right)  .
\]
Finally, think of the symbol $|$ as an operator with degree $-1$ that acts by
sliding in from the left; then%
\[
|(U\otimes V)=\left(  -1\right)  ^{\#U}U|V.
\]

\begin{definition}
Given a partition $M|N$ of $\underline{n+1},$ define \underline{face operators
with respect to} \underline{$M$ and $N$}, $d_{M},d^{N}:C_{n}\left(
P_{n+1}\right)  \rightarrow C_{n-1}\left(  P_{n+1}\right)  $ by%
\[
d_{M}\left(  \underline{n+1}\right)  =d^{N}\left(  \underline{n+1}\right)
=\left(  -1\right)  ^{\#M}\text{ shuff}\left(  M;N\right)  \text{ }M|N.
\]
For $u=U_{1}|\cdots|U_{p}\subset P_{n+1}$ and non-empty $M\subset U_{k},$
define the \underline{face operator with} \underline{respect to $M$},
$d_{M}^{k}:C_{n-p+1}\left(  u\right)  \rightarrow C_{n-p}\left(  u\right)  ,$
by%
\[
d_{M}^{k}(u)=(1^{\otimes k-1}\otimes d_{M}\otimes1^{\otimes p-k})(u);
\]
for $v=V_{q}|\cdots|V_{1}\subset P_{n+1}$ and non-empty $N\subset V_{k},$
define the \underline{face operator with} \underline{respect to $N$},
$d_{k}^{N}:C_{n-q+1}\left(  v\right)  \rightarrow C_{n-q}\left(  v\right)  ,$
by%
\[
d_{k}^{N}(v)=(1^{\otimes q-k}\otimes d^{N}\otimes1^{\otimes k-1})(v).
\]

\end{definition}

\noindent Then
\[
d_{M}^{k}\left(  u\right)  =\epsilon\left(  M\right)  \text{ }U_{1}%
|\cdots|M|U_{k}\setminus M|\cdots|U_{p},
\]
\
\[
d_{k}^{N}\left(  v\right)  =\epsilon\left(  N\right)  \text{ }V_{q}%
|\cdots|V_{k}\setminus N|N|\cdots|V_{1},
\]
where
\[
\epsilon\left(  M\right)  =\left(  -1\right)  ^{m_{1}+\cdots+m_{k-1}%
+\#M}\text{\textit{shuff}}\left(  M;U_{k}\setminus M\right)  \text{ \ and
\ }m_{i}=\#U_{i}-1,
\]
\
\[
\epsilon\left(  N\right)  =\left(  -1\right)  ^{n_{q}+\cdots+n_{k+1}+\#\left(
V_{k}\setminus N\right)  }\text{\textit{shuff}}\left(  V_{k}\setminus
N;N\right)  \text{ \ and \ }n_{i}=\#V_{i}-1.\vspace*{0.05in}%
\]

Face operators give rise to boundary operators $\partial:C_{n-p+1}\left(
u\right)  \rightarrow C_{n-p}\left(  u\right)  $ and $\partial:C_{n-q+1}%
\left(  v\right)  \rightarrow C_{n-q}\left(  v\right)  $ in the standard way:
\[
\partial\left(  u\right)  =\sum_{\substack{1\leq k\leq p \\M\subset U_{k}%
}}\epsilon(M)\,d_{M}^{k}\left(  U_{1}|\cdots|U_{p}\right)
\]
and similarly for $\partial\left(  v\right)  $; in either case,
\begin{equation}
\partial\left(  \underline{n+1}\right)  =\sum_{\substack{M,N\subset
\underline{n+1} \\N=\underline{n+1}\setminus M}}\left(  -1\right)
^{\#M}\text{ \textit{shuff}}(M;N)\text{ }M|N. \label{milgram}%
\end{equation}
The sign coefficients in \ref{milgram} were given by R. J. Milgram in
\cite{Milgram}. Thus, two types of signs appear when $d_{M}^{k}$ is applied to
$U_{1}|\cdots|U_{p}$ : First, Koszul's sign appears when $d_{M}$ passes
$U_{1}\otimes\cdots\otimes U_{k-1}$ and second, Milgram's sign appears when
$d_{M}$ is applied to $U_{k}.$

A \emph{partitioning procedure} is a composition of the form%
\[
d_{M_{p-1}}^{k_{p-1}}\cdots d_{M_{2}}^{k_{2}}d_{M_{1}}.
\]
For example, a partition $u=U_{1}|\cdots|U_{p}$ of $\underline{n+1}$ can be
obtained from the \emph{right-most }partitioning procedure by setting
$M_{0}=\underline{n+1},$ $M_{i}=M_{i-1}\setminus U_{p-i+1}$ and $k_{i}=1$ for
$1\leq i\leq p-1;$ then%
\[
d_{M_{p-1}}^{1}\cdots d_{M_{2}}^{1}d_{M_{1}}\left(  \underline{n+1}\right)
=sgn_{1}\left(  u\right)  \text{ }U_{1}|\cdots|U_{p},
\]
\ where\
\[
sgn_{1}\left(  u\right)  =\left(  -1\right)  ^{\epsilon_{1}}psgn\left(
u\right)  \text{\ \ and\ \ }\epsilon_{1}=\sum\nolimits_{i=1}^{p-1}%
i\cdot\#U_{p-i}.
\]
Note that when $v=V_{q}|\cdots|V_{1}$ we have $\epsilon_{1}=\sum
\nolimits_{i=1}^{q-1}i\cdot\#V_{i+1}$. Alternatively, $u$ can be obtained from
the \emph{left-most }partitioning procedure by setting $M_{i}=U_{i}$ and
$k_{i}=i$ for $1\leq i\leq p-1;$ then%
\[
d_{U_{p-1}}^{p-1}\cdots d_{U_{2}}^{2}d_{U_{1}}\left(  \underline{n+1}\right)
=sgn_{2}\left(  u\right)  \text{ }U_{1}|\cdots|U_{p},
\]
\ where%
\[
sgn_{2}\left(  u\right)  =\left(  -1\right)  ^{\epsilon_{2}}psgn\left(
u\right)  \ \ \text{and }\ \epsilon_{2}=\epsilon_{1}+\tbinom{p-1}{2}.
\]

Let $rsgn(U_{i})$ denote the sign of the order-reversing permutation on
$U_{i},$ then%
\[
rsgn(U_{i})=\left(  -1\right)  ^{\frac{1}{2}\left(  \#U_{i}\right)  \left(
\#U_{i}-1\right)  };
\]
define%
\[
rsgn(u)=\prod_{i=1}^{p}rsgn(U_{i})=(-1)^{\frac{1}{2}\left[  (\#U_{1}%
)^{2}+\cdots+(\#U_{p})^{2}-(n+1)\right]  }.
\]

\begin{definition}
If $F\in\mathcal{C}^{q\times p}$ is derived from $E\in\mathcal{E},$ the
\underline{configuration sign of $F$} is defined to be%
\[
csgn(F)=(-1)^{\binom{q}{2}}\text{ }rsgn(c(E))\cdot sgn_{1}r(F)\cdot
sgn_{2}c(E)\cdot sgn_{2}c(F).
\]

\end{definition}

\noindent In particular, for $F=E\in\mathcal{E}^{q\times p}$ we have
\[
csgn(E)=(-1)^{\binom{q}{2}}\text{ }rsgn(c(E))\cdot sgn_{1}r(E).
\]

Signs that arise from the action of shift operators are now determined. For
$x\in\mathbb{Z}$ and $Y\subseteq\mathbb{Z},$ denote the lower and upper cuts
of $Y$ at $x$ by $[Y,x)=\left\{  y\in Y\text{ }|\text{ }y<x\right\}  $ and
$(x,Y]=\left\{  y\in Y\text{ }|\text{ }y>x\right\}  ,$ respectively.

\begin{proposition}
If $F=\left(  f_{i,j}\right)  \in\mathcal{C}$, $c\left(  F\right)
=U_{1}|\cdots|U_{p}$ and $r\left(  F\right)  =V_{q}|\cdots|V_{1}$, then%
\begin{align*}
csgn(D_{i,j}F)\cdot csgn\left(  F\right)   &  =-(-1)^{\#(f_{i+1,j}^{\prime
},V_{i+1}^{\prime}]\cup\lbrack V_{i},f_{i,j})},\vspace*{0.1in}\\
csgn(R_{i,j}F)\cdot csgn\left(  F\right)   &  =-(-1)^{\#(f_{i,j},U_{j}%
]\cup\lbrack U_{j+1}^{\prime},f_{i,j+1}^{\prime})},\vspace*{0.1in}%
\end{align*}
where $F^{\prime}=\left(  f_{i,j}^{\prime}\right)  $ is the image of $F,$
$U_{1}^{\prime}|\cdots|U_{p}^{\prime}=$ $c\left(  F^{\prime}\right)  $ and
$V_{q}^{\prime}|\cdots|V_{1}^{\prime}=r\left(  F^{\prime}\right)
.\vspace*{0.1in}$
\end{proposition}

\begin{proof}
Note that $c\left(  F\right)  =c\left(  D_{i,j}F\right)  $ and $r\left(
F\right)  =r\left(  R_{i,j}F\right)  .$ Then for example,
\begin{align*}
csgn(D_{i,j}F)\cdot csgn(F)  &  =(-1)^{\binom{q}{2}}\text{\thinspace
}rsgn(c(E))\!\cdot\!sgn_{1}r(D_{i,j}F)\!\cdot\!sgn_{2}c(E)\!\cdot
\!sgn_{2}c(D_{i,j}F)\\
&  \cdot\text{ }(-1)^{\binom{q}{2}}\text{ }rsgn(c(E))\cdot sgn_{1}r\left(
F\right)  \cdot sgn_{2}c(E)\cdot sgn_{2}c(F)\\
&  =sgn_{1}r(F)\cdot sgn_{1}r(D_{i,j}F)\cdot sgn_{2}c(F)\cdot sgn_{2}%
c(D_{i,j}F)\\
&  =psgn\left(  r\left(  F\right)  \right)  \cdot psgn\left(  r\left(
D_{i,j}F\right)  \right)  =-sgn\left(  \sigma\right)  ,
\end{align*}
where $\sigma$ is the permutation $V_{q}\cup\cdots\cup V_{1}\mapsto V_{q}%
\cup\cdots V_{i+1}^{\prime}\cup V_{i}^{\prime}\cdots\cup V_{1}.$
\end{proof}

The configuration signs of \textquotedblleft edge matrices,\textquotedblright%
\ which appear in our subsequent discussion of permutahedral sets, have a
particularly nice form.

\begin{definition}
$E\in\mathcal{E}$ is an \underline{edge matrix} if $e_{1,1}=1$.
\end{definition}

Let $\Gamma$ denote the set of all edge matrices. With one possible exception,
all blocks in the column and row face of an edge matrix consist of singleton
sets. Thus if $E\in\Gamma^{q\times p},$%
\[
c\left(  E\right)  \otimes r\left(  E\right)  =A|a_{2}|\cdots|a_{p}\otimes
b_{q}|\cdots|b_{2}|B,
\]
where $A=\left\{  1<b_{2}<\cdots<b_{q}\right\}  $ and $B=\left\{
1<a_{2}<\cdots<a_{p}\right\}  .$ Since $c\left(  E\right)  $ and $r\left(
E\right)  $ meet at the cubical vertex $b_{q}|\cdots|b_{2}|1|a_{2}%
|\cdots|a_{p}$ of $P_{p+q-1},$ there is a canonical bijection%
\[
\Gamma\leftrightarrow\left\{  \text{\textit{cubical vertices of} }P=\sqcup
P_{n+1}\right\}  .
\]
The proof of the following proposition is now immediate:

\begin{proposition}
If $E$ is an edge matrix and $b_{q}|\cdots|b_{2}|1|a_{2}|\cdots|a_{p}$ is the
corresponding cubical vertex, then%
\[
csgn\left(  E\right)  =\text{shuff}\left(  b_{2},\ldots,b_{q};a_{2}%
,\ldots,a_{p}\right)  .
\]

\end{proposition}

We are ready to define a diagonal on $C_{\ast}\left(  P_{n+1}\right)  .$

\begin{definition}
For each $n\geq0,$ define $\Delta_{P}$ on the top dimensional face
$\underline{n+1}\in C_{n}\left(  P_{n+1}\right)  $ by%
\begin{equation}
\Delta_{P}\left(  \underline{n+1}\right)  =\sum\limits_{\substack{F\in
\mathcal{C}^{q\times n-q+2}\\1\leq q\leq n+1}}csgn(F)\,c\left(  F\right)
\otimes r\left(  F\right)  ; \label{diagonal}%
\end{equation}
extend $\Delta_{P}$ to proper faces $u=U_{1}|\cdots|U_{p}\in C_{n-p+1}\left(
u\right)  =C_{n_{1}}\left(  U_{1}\right)  \otimes\cdots\otimes C_{n_{p}%
}\left(  U_{p}\right)  ,$ $n_{i}=\#U_{i}-1,$ via the standard comultiplicative extension.
\end{definition}

\begin{example}
\label{ex5}On $P_{3},$ all but two configuration matrices are step
matrices:\vspace*{0.05in}%
\[%
\begin{tabular}
[c]{ccccccc}%
$%
\begin{tabular}
[c]{|l|l|l|}\hline
1 & 2 & 3\\\hline
\end{tabular}
\ \ \ $ &  &
\begin{tabular}
[c]{|l|}\hline
1\\\hline
2\\\hline
3\\\hline
\end{tabular}
&  & $%
\begin{tabular}
[c]{|l|l|}\hline
1 & 2\\\hline
3 & \\\hline
\end{tabular}
\ \ \ $ & $\longrightarrow R_{3}\longrightarrow$ & $%
\begin{tabular}
[c]{|l|l|}\hline
1 & 2\\\hline
& 3\\\hline
\end{tabular}
\ \ \ $\\
&  &  &  &  &  & \\%
\begin{tabular}
[c]{|l|l|}\hline
& 2\\\hline
1 & 3\\\hline
\end{tabular}
&  & $%
\begin{tabular}
[c]{|l|l|}\hline
& 1\\\hline
2 & 3\\\hline
\end{tabular}
\ \ \ $ &  & $\text{%
\begin{tabular}
[c]{|l|l|}\hline
1 & 3\\\hline
2 & \\\hline
\end{tabular}
}$ & $\longrightarrow D_{3}\longrightarrow$ & $%
\begin{tabular}
[c]{|l|l|}\hline
1 & \\\hline
2 & 3\\\hline
\end{tabular}
\ \ \ $%
\end{tabular}
\ \ \ \vspace{0.1in}%
\]
Consequently,%
\[%
\begin{array}
[c]{lllll}%
\Delta_{P}\left(  \underline{3}\right)  & = & 1|2|3\otimes{123} & + &
12{3}\otimes3|2|1\\
&  &  &  & \\
& - & 1|23\otimes13|2 & + & 2|13\otimes23|1\\
&  &  &  & \\
& - & 13|2\otimes3|12 & + & 12|3\otimes2|13\\
&  &  &  & \\
& - & 1|23\otimes3|12 & + & 12|3\otimes23|1.
\end{array}
\]

\end{example}

There is a computational shortcut worth mentioning. Since $F\in\mathcal{C}$ if
and only if $F^{T}\in\mathcal{C}$, we only need to derive half of the
configuration matrices.

\begin{definition}
For $F\in\mathcal{C}$, define the \underline{transpose} of $c\left(  F\right)
\otimes r\left(  F\right)  $ to be%
\[
\left[  c\left(  F\right)  \otimes r\left(  F\right)  \right]  ^{T}=c\left(
F^{T}\right)  \otimes r\left(  F^{T}\right)  .
\]

\end{definition}

\begin{example}
Refer to \textit{Example \ref{ex5} and note} that e\textit{ach component in
the right-hand column is the transpose of the component to its left. On
}$P_{4}$ we have:%
\[%
\begin{array}
[c]{ll}%
\Delta_{P}\left(  \underline{4}\right)  & =1234\otimes4|3|2|1\\
& +123|4\otimes\left(  3|2|14+3|24|1+34|2|1\right) \\
& -12|34\otimes\left(  2|14|3+24|1|3\right) \\
& +1|234\otimes14|3|2\\
& -23|14\otimes\left(  3|24|1+34|2|1\right) \\
& +13|24\otimes\left(  3|14|2+34|1|2\right) \\
& +\left(  13|24+1|234-14|23+134|2\right)  \otimes4|3|12\\
& -\left(  12|34+124|3\right)  \otimes\left(  4|2|13+4|23|1\right) \\
& +3|124\otimes34|2|1\\
& -2|134\otimes(24|3|1+4|23|1)\\
& +24|13\otimes4|23|1\\
& +\left(  1|234-14|23\right)  \otimes4|13|2\\
& \pm\text{\textit{(all transposes of the above).}}%
\end{array}
\]

\end{example}

We conclude this section with a proof of the fact that $\Delta_{P}$ is a chain
map. First note that%
\begin{align*}
\Delta_{P}\partial\left(  \underline{n+1}\right)   &  =\sum\pm\Delta
_{P}\left(  M\right)  |\Delta_{P}\left(  N\right) \\
&  =\sum\pm\left(  u_{i}\otimes v_{j}\right)  |\left(  u^{k}\otimes v^{\ell
}\right)  =\sum\pm u_{i}|u^{k}\otimes v_{j}|v^{\ell},
\end{align*}
where $u_{i}\otimes v_{j}=c\left(  F_{j\times i}\right)  \otimes r\left(
F_{j\times i}\right)  ,$ $u^{k}\otimes v^{\ell}=c\left(  F^{\ell\times
k}\right)  \otimes r\left(  F^{\ell\times k}\right)  \ $and $F_{j\times i}$
and $F^{\ell\times k}$ range over all configurations matrices with entries
from $M$ and $N,$ respectively. Although $u_{i}|u^{k}\otimes v_{j}|v^{\ell}$
is not a CP, there is the associated block matrix%
\begin{equation}%
\begin{tabular}
[c]{|c|c|}\hline
$0$ & $\overset{\text{\mathstrut}}{\underset{\text{\mathstrut}}{F^{\ell\times
k}}}$\\\hline
$\overset{\text{\mathstrut}}{\underset{\text{\mathstrut}}{F_{j\times i}}}$ &
$0$\\\hline
\end{tabular}
\ . \label{block}%
\end{equation}
Thus the components of $\Delta_{P}\partial\left(  \underline{n+1}\right)  $
lie in one-to-one correspondence with all such block matrices. Let
$a_{i}\otimes b_{j}=A_{1}|\cdots|A_{i}\otimes B_{j}|\cdots|B_{1}$ and
$a^{k}\otimes b^{\ell}=A^{1}|\cdots|A^{k}\otimes B^{\ell}|\cdots|B^{1}$ be the
SCP's related to $u_{i}\otimes v_{j}$ and $u^{k}\otimes v^{\ell}.$ Denoting a
column (or row) by its set of non-zero entries, the step matrices%
\[%
\begin{tabular}
[c]{|c|c|c|}\hline
$\overset{\text{\mathstrut}}{\underset{\text{\mathstrut}}{A_{1}}}$ & $\cdots$
& $\overset{\text{\mathstrut}}{\underset{\text{\mathstrut}}{A_{i}}}$\\\hline
\end{tabular}
\ \ =\ \
\begin{tabular}
[c]{|c|}\hline
$\overset{\text{\mathstrut}}{\underset{\text{\mathstrut}}{B_{1}}}$\\\hline
$\underset{}{\vdots}$\\\hline
$\overset{\text{\mathstrut}}{\underset{\text{\mathstrut}}{B_{j}}}$\\\hline
\end{tabular}
\ \ \text{and }\
\begin{tabular}
[c]{|c|c|c|}\hline
$\overset{\text{\mathstrut}}{\underset{\text{\mathstrut}}{A^{1}}}$ & $\cdots$
& $\overset{\text{\mathstrut}}{\underset{\text{\mathstrut}}{A^{k}}}$\\\hline
\end{tabular}
\ \ =\ \
\begin{tabular}
[c]{|c|}\hline
$\overset{\text{\mathstrut}}{\underset{\text{\mathstrut}}{B^{1}}}$\\\hline
$\underset{}{\vdots}$\\\hline
$\overset{\text{\mathstrut}}{\underset{\text{\mathstrut}}{B^{\ell}}}$\\\hline
\end{tabular}
\]
involve elements of $M$ and $N,$ respectively, and the block matrix associated
with the pairing $a_{i}|a^{k}\otimes b_{j}|b^{\ell}$ is%
\[%
\begin{tabular}
[c]{|ccc|ccc|}\hline
& $0$ &  & $\overset{\text{\mathstrut}}{\underset{\text{\mathstrut}}{A^{1}}}$
& \multicolumn{1}{|c}{$\cdots$} & \multicolumn{1}{|c|}{$\overset
{\text{\mathstrut}}{\underset{\text{\mathstrut}}{A^{k}}}$}\\\hline
$\overset{\text{\mathstrut}}{\underset{\text{\mathstrut}}{A_{1}}}$ &
\multicolumn{1}{|c}{$\cdots$} & \multicolumn{1}{|c|}{$\overset
{\text{\mathstrut}}{\underset{\text{\mathstrut}}{A_{i}}}$} &  & $0$ & \\\hline
\end{tabular}
\ \ =\ \
\begin{tabular}
[c]{|c|c|}\hline
& $\overset{\text{\mathstrut}}{\underset{\text{\mathstrut}}{B^{1}}}%
$\\\cline{2-2}%
$0$ & $\underset{}{\vdots}$\\\cline{2-2}
& $\overset{\text{\mathstrut}}{\underset{\text{\mathstrut}}{B^{\ell}}}%
$\\\hline
$\overset{\text{\mathstrut}}{\underset{\text{\mathstrut}}{B_{1}}}$ &
\\\cline{1-1}%
$\underset{}{\vdots}$ & $0$\\\cline{1-1}%
$\overset{\text{\mathstrut}}{\underset{\text{\mathstrut}}{B_{j}}}$ & \\\hline
\end{tabular}
\ \ .
\]
Our main result combines the statements in Lemmas \ref{Lemma1} and
\ref{Lemma2} below:

\begin{theorem}
\label{main}The cellular boundary map $\partial:C_{\ast}\left(  P_{n+1}%
\right)  \rightarrow C_{\ast}\left(  P_{n+1}\right)  $ is a $\Delta_{P}
$-coderivation for all $n\geq1.$
\end{theorem}

\begin{corollary}
$\left(  C_{\ast}\left(  P_{n+1}\right)  ,\Delta_{P},\partial\right)  $ is a
DG coalgebra and the cellular projection $\rho_{n+1}:P_{n+1}\rightarrow I^{n}$
induces a DG coalgebra map%
\[
(\rho_{n+1})_{\ast}:C_{\ast}(P_{n+1})\rightarrow C_{\ast}(I^{n}).
\]

\end{corollary}

\begin{lemma}
\label{Lemma1}Each non-zero component $\left(  u_{i}\otimes v_{j}\right)
|\left(  u^{k}\otimes v^{\ell}\right)  $ of $\Delta_{P}\partial\left(
\underline{n+1}\right)  $ is a non-zero component of$\ \left(  1\otimes
\partial+\partial\otimes1\right)  \Delta_{P}\left(  \underline{n+1}\right)
.$\vspace{0.1in}
\end{lemma}

\begin{proof}
Consider a component $\left(  u_{i}\otimes v_{j}\right)  |\left(  u^{k}\otimes
v^{\ell}\right)  $ of $\Delta_{P}\partial\left(  \underline{n+1}\right)  ,$
where $u_{i}\otimes v_{j}=U_{1}|\cdots|U_{i}\otimes V_{j}|\cdots|V_{1}$ is a
CP of partitions of $M=U_{1}\cup\cdots\cup U_{i}$ and $u^{k}\otimes v^{\ell
}=U^{1}|\cdots|U^{k}\otimes V^{\ell}|\cdots|V^{1}$ is a CP of partitions of
$N=\underline{n+1}\diagdown M.$ The related SCP's $a_{i}\otimes b_{j}%
=A_{1}|\cdots|A_{i}\otimes B_{j}|\cdots|B_{1}$ and $a^{k}\otimes b^{\ell
}=A^{1}|\cdots|A^{k}\otimes B^{\ell}|\cdots|B^{1}$ give the component $\left(
a_{i}\otimes b_{j}\right)  |\left(  a^{k}\otimes b^{\ell}\right)  $ of
$\Delta_{P}\partial\left(  \underline{n+1}\right)  .$ Let $E=\left(
e_{i,j}\right)  $ be the block matrix associated with $a_{i}|a^{k}\otimes
b_{j}|b^{\ell}$. There are two cases:\vspace{0.1in}\newline\underline
{\textit{Case 1}}:\textit{\ }$e_{\ell+1,i}>e_{\ell,i+1}.$\vspace
{0.1in}\newline Then $\min U_{i}=\min A_{i}>\max A^{1}\geq\max U^{1}\ $and the
CP%
\[
u\otimes v=U_{1}|\cdots|U_{i-1}|U^{1}\cup U_{i}|U^{2}|\cdots|U^{k}\otimes
v_{j}|v^{\ell}%
\]
is a component of $\Delta_{P}\left(  \underline{n+1}\right)  $ with associated
configuration matrix%
\[
F\text{ }=\text{ }%
\begin{tabular}
[c]{|rc|c|cc|}\hline
$\qquad0\hspace*{-0.1in}$ &  & $\underset{\text{\mathstrut}}{\overset
{\text{\mathstrut}}{U^{1}}}$ & $\cdots$ & \multicolumn{1}{|c|}{$\underset
{\text{\mathstrut}}{\overset{\text{\mathstrut}}{U^{k}}}$}\\\hline
\multicolumn{1}{|c}{$\underset{\text{\mathstrut}}{\overset{\text{\mathstrut}%
}{U_{1}}}$} & \multicolumn{1}{|c|}{$\cdots$} & $\underset{\text{\mathstrut}%
}{\overset{\text{\mathstrut}}{U_{i}}}$ &  & \multicolumn{1}{l|}{$\hspace
*{-0.1in}0\qquad$}\\\hline
\end{tabular}
\text{ }=\text{ }%
\begin{tabular}
[c]{|c|}\hline
$\underset{\text{\mathstrut}}{\overset{\text{\mathstrut}}{v^{\ell}}}$\\\hline
$\overset{\text{\mathstrut}}{\underset{\text{\mathstrut}}{v_{j}}}$\\\hline
\end{tabular}
\ .\
\]
It follows that $u_{i}|u^{k}\otimes v=d_{U_{i}}^{i}\left(  u\right)  \otimes
v$ is a component of $\left(  1\otimes\partial+\partial\otimes1\right)
\Delta_{P}\left(  \underline{n+1}\right)  .$ To check signs, we verify that
the product of expressions (I) through (VI) below is $1.$ Let $V_{q}^{\prime
}|...|V_{1}^{\prime}=v_{j}|v^{\ell}$ and note that $u\otimes v=c\left(
F\right)  \otimes r\left(  F\right)  $ is related to the SCP $a\otimes
b=A_{1}|\cdots|A_{i-1}|A^{1}\cup A_{i}|A^{2}|\cdots|A^{k}\otimes b_{j}%
|b^{\ell}=c\left(  E\right)  \otimes r\left(  E\right)  .$

\begin{enumerate}
\item[I.] $csgn(F)=$ I$_{1}$ $\cdot$ I$_{2}$ $\cdot$ I$_{3}$ $\cdot$ I$_{4}$
$\cdot$ I$_{5}=(-1)^{\tbinom{q}{2}}\cdot\left[  sgn_{2}u\cdot sgn_{2}a\right]
\cdot rsgn(a)\cdot(-1)^{\epsilon_{1}^{\prime}}\cdot psgn(v),$ where
$\epsilon_{1}^{\prime}=\sum\nolimits_{i=1}^{q-1}i\cdot\#V_{i+1}^{\prime
}.\vspace{0.05in}$

\item[II.] $sgn(d_{U_{i}}^{i}\left(  u\right)  )=$ II$_{1}$ $\cdot$
II$_{2}=(-1)^{\#M+i+1}\cdot(-1)^{\#U_{i}\#U^{1}},$ where the shuffle sign
II$_{2}$ follows by assumption.$\vspace{0.05in}$

\item[III.] $sgn\left(  d_{M}\left(  \underline{n+1}\right)  \right)  =$
III$_{1} $ $\cdot$ III$_{2}=(-1)^{\#M}$ $\cdot$ shuff$(M;N).\vspace{0.05in}$

\item[IV.] $csgn(F_{j\times i})=$ IV$_{1}$ $\cdot$ IV$_{2}$ $\cdot$ IV$_{3} $
$\cdot$ IV$_{4}$ $\cdot$ IV$_{5}=(-1)^{\tbinom{j}{2}}\cdot\lbrack sgn_{2}%
u_{i}\cdot sgn_{2}a_{i}]\cdot rsgn(a_{i})\cdot(-1)^{\epsilon_{1}}\cdot
psgn(v_{j}),$ where $\epsilon_{1}=\sum\nolimits_{i=1}^{j-1}i\cdot
\#V_{i+1}.\vspace{0.05in}$

\item[V.] $csgn(F^{\ell\times k})=$ V$_{1}$ $\cdot$ V$_{2}$ $\cdot$ V$_{3} $
$\cdot$ V$_{4}$ $\cdot$ V$_{5}=(-1)^{\tbinom{\ell}{2}}\cdot\lbrack
sgn_{2}u^{k}\cdot sgn_{2}a^{k}]\cdot rsgn(a^{k})\cdot(-1)^{\epsilon^{1}}\cdot
psgn(v^{\ell}),$ where $\epsilon^{1}=\sum\nolimits_{i=1}^{\ell-1}%
i\cdot\#V^{i+1}.\vspace{0.05in}$

\item[VI.] $\left(  -1\right)  ^{\dim u^{k}\dim v_{j}}=(-1)^{(\ell-1)(i-1)}$
($u_{i}\otimes v_{j}$ is a component of $\Delta_{P}\left(  M\right)  ;$ hence
$\dim\left(  u_{i}\otimes v_{j}\right)  =\#M-1$ and $\dim v_{j}=\#M-1-\dim
u_{i}=i-1$).
\end{enumerate}

\noindent Then by straightforward calculation,

\begin{enumerate}
\item[(1)] I$_{5}$ $\cdot$ III$_{2}$ $\cdot$ IV$_{5}$ $\cdot$ V$_{5}%
=1;\vspace{0.05in}$

\item[(2)] I$_{2}$ $\cdot$ IV$_{2}$ $\cdot$ V$_{2}=$ I$_{3}$ $\cdot$ II$_{2}$
$\cdot$ IV$_{3}$ $\cdot$ V$_{3}=\left(  -1\right)  ^{\#A_{i}\#A^{1}%
+\#U_{i}\#U^{1}};\vspace{0.05in}$

\item[(3)] I$_{1}$ $\cdot$ IV$_{1}$ $\cdot$ V$_{1}$ $=$ (I$_{4}$ $\cdot$
IV$_{4}$ $\cdot$ V$_{4}$) $\cdot$ (II$_{1}$ $\cdot$ III$_{1}$) $\cdot$ VI
$=\left(  -1\right)  ^{j\ell}\vspace{0.05in}$\newline($\#M=i+j-1$ since
$v_{j}=r\left(  F_{j\times i}\right)  $)$.$
\end{enumerate}

\noindent\underline{\textit{Case 2}}:\textit{\ }$e_{\ell+1,i}<e_{\ell,i+1}%
.$\vspace{0.1in}\newline Then $\max\left(  V_{1}\right)  \leq\max\left(
B_{1}\right)  <\min\left(  B^{\ell}\right)  =\min\left(  V^{\ell}\right)
\ $and the CP%
\[
u\otimes v=u_{i}|u^{k}\otimes V_{j}|\cdots|V_{1}\cup V^{\ell}|\cdots|V^{1}%
\]
is a component of $\Delta_{P}\left(  \underline{n+1}\right)  $ with associated
configuration matrix%
\[
F\text{ \ }=\text{ \ }%
\begin{tabular}
[c]{|c|c|}\hline
& $\underset{\text{\mathstrut}}{\overset{\text{\mathstrut}}{V^{1}}}%
$\\\cline{2-2}%
$%
\begin{array}
[c]{c}%
0\\
\text{\mathstrut}\\
\text{\mathstrut}%
\end{array}
$ & $\vdots$\\\hline
$\underset{\text{\mathstrut}}{\overset{\text{\mathstrut}}{V_{1}}}$ &
$\underset{\text{\mathstrut}}{\overset{\text{\mathstrut}}{V^{\ell}}}$\\\hline
$\vdots$ & $%
\begin{array}
[c]{c}%
\text{\mathstrut}\\
\text{\mathstrut}\\
0
\end{array}
$\\\cline{1-1}%
$\underset{\text{\mathstrut}}{\overset{\text{\mathstrut}}{V_{j}}}$ & \\\hline
\end{tabular}
\text{ \ }=\text{ \ }%
\begin{tabular}
[c]{|l|l|}\hline
$\ \underset{\text{\mathstrut}}{\overset{\text{\mathstrut}}{u_{i}}}$ \  &
$\ \underset{\text{\mathstrut}}{\overset{\text{\mathstrut}}{u^{k}}}$
\ \\\hline
\end{tabular}
\]
It follows that $u_{i}|u^{k}\otimes v_{j}|v^{\ell}=u\otimes d_{\ell}^{V_{1}%
}\left(  v\right)  $ is a component of $\left(  1\otimes\partial
+\partial\otimes1\right)  \Delta_{P}\left(  \underline{n+1}\right)  $. The
sign check is similar to the one in \textit{Case 1} above and is left to the
reader.\vspace{0.1in}
\end{proof}

\begin{lemma}
\label{Lemma2}Each non-zero component $d_{M}^{k}(u)\otimes v$ or $u\otimes
d_{\ell}^{N}\left(  v\right)  $ of $\left(  1\otimes\partial+\partial
\otimes1\right)  $ $\Delta_{P}\left(  \underline{n+1}\right)  $ is a non-zero
component of $\Delta_{P}\partial\left(  \underline{n+1}\right)  .$%
\vspace{0.1in}
\end{lemma}

\begin{proof}
For simplicity we work with $\mathbb{Z}_{2}$ coefficients; sign checks with
$\mathbb{Z}$ coefficients are straightforward calculations and left to the
reader. Given an SCP $a\otimes b=c\left(  E\right)  \otimes r\left(  E\right)
=A_{1}|\cdots A_{p}\otimes B_{q}|\cdots|B_{1}$ of partitions of $\underline
{n+1},$ let $u\otimes v=c\left(  F\right)  \otimes r\left(  F\right)
=U_{1}|\cdots|U_{p}\otimes V_{q}|\cdots|V_{1}$ be a related CP. Then there
exist $M_{j}\subset A_{j}$ and $N_{i}\subset B_{i}$ with $\min M_{j}>\max
A_{j+1}$ and $\min N_{i}>\max B_{i+1}$ such that%
\begin{equation}
F=D_{N_{q-1}}\cdots D_{N_{1}}R_{M_{p-1}}\cdots R_{M_{1}}E. \label{f}%
\end{equation}
Then $u\otimes v$ is a non-zero component of $\Delta_{P}\left(  \underline
{n+1}\right)  .$ For each proper $M\subset U_{k}$, we prove that the component
$d_{M}^{k}(u)\otimes v$ of $\left(  1\otimes\partial+\partial\otimes1\right)
\Delta_{P}\left(  \underline{n+1}\right)  $ is a non-zero component of
$\Delta_{P}\partial\left(  \underline{n+1}\right)  $ if and only if the
following conditions hold:\vspace{0.1in}

\begin{enumerate}
\item[(1)] $m=\min M\in A_{k};$\vspace{0.1in}

\item[(2)] $\left(  m,M\right]  =\left(  m,A_{k}\cup M_{k-1}\right]  ;$%
\vspace{0.1in}

\item[(3)] $m\in B_{r}$ implies $N_{r-1}=\varnothing.$\vspace{0.1in}
\end{enumerate}

\noindent The dual statement for $u\otimes d_{\ell}^{N}\left(  v\right)  $
with $N\subset V_{\ell}$ and is also true; the proof follows by "mirror
symmetry." Suppose conditions (1) - (3) hold. Set $M_{0}=M_{p}=\varnothing;$
then clearly, $U_{i}=\left(  A_{i}\cup M_{i-1}\right)  \setminus M_{i}$ for
$1\leq i\leq p,$ and $M_{k-1}\subseteq M$ by conditions (1) and (2). Thus
$U_{1}\cup\cdots\cup U_{k-1}\cup M=A_{1}\cup\cdots\cup A_{k-1}\cup M$ and it
follows that $d_{M}^{k}(u)\otimes v$ is the non-zero component
\[
\Delta_{P}\left(  A_{1}\cup\cdots\cup A_{k-1}\cup M\text{ }|\text{ }%
A_{k}\setminus M\cup A_{k+1}\cup\cdots\cup A_{p}\right)
\]
of $\Delta_{P}\partial\left(  \underline{n+1}\right)  .$ Conversely, if
conditions (1) - (3) fail to hold, we prove that there exists a unique CP
$\bar{u}\otimes\bar{v}\neq u\otimes v$ such that $u\otimes v+\bar{u}%
\otimes\bar{v}\in\ker\left(  \partial\otimes1+1\otimes\partial\right)  .$
\vspace{0.1in}\newline For existence, we consider all possible cases.\vspace
{0.1in}\newline\underline{\textit{Case 1}}: Assume $\left(  1\right)
^{\prime}:m\notin A_{k}.$\vspace{0.1in}\newline Let
\[
\bar{u}=U_{1}|\cdots|U_{k-1}\cup M|U_{k}\setminus M|\cdots|U_{p};
\]
then
\[
d_{U_{k-1}}^{k-1}(\bar{u})\otimes v=d_{M}^{k}(u)\otimes v.\vspace{0.1in}%
\]
Now $M\subset M_{k-1}$ since $m\in M_{k-1}$; hence $\bar{u}\otimes v$ may be
obtained by replacing $R_{M_{k-1}}$ with $R_{M_{k-1}\setminus M}$ in (\ref{f})
and $\bar{u}\otimes v$ is a CP related to $a\otimes b$.\vspace{0.1in}%
\newline\underline{\textit{Case 2}}: Assume $\left(  1\right)  \wedge\left(
2\right)  ^{\prime}:m\in A_{k}$ and $\left(  m,M\right]  \subset\left(
m,A_{k}\cup M_{k-1}\right]  .$\vspace{0.1in}\newline Let%
\[
\mu=\min\left(  m,A_{k}\cup M_{k-1}\right]  \setminus M\text{ \ and
\ }L=\left[  A_{k},m\right)  \cup\mu.
\]
Note that $\mu\in A_{i}\ $for some$\ 1\leq i\leq k.$\vspace{0.1in}%
\newline\underline{\textit{Subcase 2A}}: \textit{Assume }$\min L>\max
A_{k+1},$ $k<p.$\vspace{0.1in}\newline Let
\[
\bar{u}=U_{1}|\cdots|M|\left(  U_{k}\setminus M\right)  \cup U_{k+1}%
|\cdots|U_{p};
\]
then
\[
d_{U_{k}\setminus M}^{k+1}\left(  \bar{u}\right)  \otimes v=d_{M}^{k}\left(
u\right)  \otimes v.\vspace{0.1in}%
\]
Note that $\min A_{k}=m$ since $\min L>\max A_{k+1}>\min A_{k}.$ Thus $L=\mu.$
Now, $\min M_{k}>\max A_{k+1}$ by (\ref{f}) and $\min U_{k}\setminus
M=\min\left[  \left(  A_{k}\cup M_{k-1}\right)  \setminus M_{k}\right]
\setminus M\geq\min\left(  A_{k}\cup M_{k-1}\right)  \setminus M=\min\left(
m,A_{k}\cup M_{k-1}\right]  \setminus M=\mu=\min L>\max A_{k+1}$ so that $\min
M_{k}\cup\left(  U_{k}\setminus M\right)  >\max A_{k+1}.$ Hence $\bar
{u}\otimes v$ can be obtained by replacing $R_{M_{k}}$ with $R_{M_{k}%
\cup\left(  U_{k}\setminus M\right)  }\ $in (\ref{f}) and $\bar{u}\otimes v$
is a CP related to $a\otimes b.$\vspace{0.1in}\newline\underline
{\textit{Subcase 2B}}: $\min L<\max A_{k+1}$ with $k\leq p$. \vspace
{0.1in}\newline\underline{\textit{Subcase 2B1}}: \textit{Assume} $\min
A_{i-1}>\max A_{i}\setminus\mu$ \textit{with }$\mu\in A_{i}$ and $1<i\leq
k.\vspace{0.1in}\newline$When $i=k$ let
\[
\bar{u}=U_{1}|\cdots|U_{k-1}\cup M|U_{k}\setminus M|\cdots|U_{p};
\]
and when $1<i<k,$ let%
\[
\bar{u}=U_{1}|\cdots|U_{i-1}\cup U_{i}|\cdots|M|U_{k}\setminus M|\cdots
|U_{p}.
\]
Then for all $i\leq k$,%
\[
d_{M}^{k}\left(  u\right)  \otimes v=d_{U_{i-1}}^{i-1}\left(  \bar{u}\right)
\otimes v.\vspace{0.1in}%
\]
When $i=k$, $\min A_{k-1}\cup(A_{k}\cap M)\leq\min A_{k}\cap M<\mu=\max
A_{k}=\max A_{k}\setminus M$ so that%
\[
\bar{a}\otimes\bar{b}=A_{1}|\cdots|A_{k-1}\cup(A_{k}\cap M)|A_{k}\setminus
M|\cdots|A_{p}\otimes b
\]
is an SCP; let $\bar{E}$ be the associated step matrix and let
\[
\bar{F}=D_{N_{q-1}}\cdots D_{N_{1}}R_{M_{p-1}}\cdots R_{M_{k}}R_{M_{k-1}%
\setminus M}\cdots R_{M_{1}}\bar{E}.
\]
$\vspace{0.1in}\newline$When $i<k,$ we have $\mu=\max A_{i}>\max A_{k}\geq\max
A_{k}\cap M$ so that $\min A_{k}\cap M<\max L;$ furthermore, $\max A_{k}=\max
A_{k}\cap M$ by the minimality of $\mu$ so that $\min A_{k-1}<\max A_{k}\cap
M.$ And finally, $\min L<\max A_{k+1}$ by assumption 2B. Thus%
\begin{align*}
\bar{a}\otimes\bar{b}  &  =A_{1}|\cdots|A_{i-1}\cup A_{i}\setminus\mu
|\cdots|A_{k}\cap M|L|\cdots|A_{p}\hspace*{1in}\\
&  \hspace*{1.5in}\otimes B_{q}|\cdots|B_{r+1}|B_{r-1}|\cdots|B_{j}\cup
\mu|\cdots|B_{1}%
\end{align*}
is an SCP; let $\bar{E}$ be the associated step matrix. Note that $U_{i-1}\cup
U_{i}=$\linebreak$\left(  A_{i-1}\cup M_{i-2}\cup A_{i}\setminus\mu\right)
\setminus\left(  M_{i}\setminus\mu\right)  $ and $\mu\in M_{j}$ for $i\leq
j\leq k-1.$ Let%
\begin{align*}
\bar{F}  &  =D_{N_{q-1}}\cdots D_{N_{r-1}\cup\mu}\cdots D_{N_{j}\cup\mu}\cdots
D_{N_{1}}\\
&  \hspace*{0.25in}R_{M_{p-1}}\cdots R_{M_{k}}R_{\left(  M_{k-1}\setminus
\mu\right)  \setminus M}^{k-1}R_{M_{k-1}\setminus\mu}^{k-2}\cdots
R_{M_{i}\setminus\mu}^{i-1}\cdots R_{M_{1}}\bar{E},
\end{align*}
where $\mu\in B_{r},\bar{B}_{j}.$ Then for all $i\leq k$, $\bar{u}\otimes
v=c\left(  \bar{F}\right)  \otimes r\left(  \bar{F}\right)  $ is a CP related
to $\bar{a}\otimes\bar{b}.$\vspace{0.1in}\newline\underline{\textit{Subcase
2B2}}: \textit{Assume} $\min A_{i-1}<\max A_{i}\setminus\mu$ \textit{with
}$\mu\in A_{i}$ and $1<i\leq k.$\vspace{0.1in}\newline Let
\[
\bar{u}\otimes\bar{v}=U_{1}|\cdots|M|U_{k}\setminus M|\cdots|U_{p}\otimes
V_{q}|\cdots|V_{r}\cup V_{r-1}|\cdots|V_{1},
\]
where $\mu\in B_{r},\bar{B}_{j}.$ Then
\[
d_{M}^{k}(u)\otimes v=\bar{u}\otimes d_{r-1}^{V_{r-1}}(\bar{v}).\vspace{0.1in}%
\]
When $i=k$, $\max L=\mu\in A_{k}$ so that $\min A_{k-1}<\max A_{k}\setminus
\mu=\max A_{k}\setminus L$. Furthermore, $\min A_{k}\setminus L=m<\mu=\max L;$
and finally, $\min L<\max A_{k+1}$ by assumption 2B. Thus
\[
\bar{a}\otimes\bar{b}=A_{1}|\cdots|M|L|\cdots|A_{p}\otimes B_{q}%
|\cdots|B_{r+1}|B_{r-1}|\cdots|B_{j}\cup\mu|\cdots|B_{1}%
\]
is an SCP; let $\bar{E}$ be the associated step matrix. Since $\min\left(
\mu,A_{k}\cup M_{k-1}\right]  \setminus M>\mu=\max L,$ the operator
$R_{\left(  \mu,A_{k}\cup M_{k-1}\right]  \setminus M}^{k}$ is defined. Note
that $M_{k}\subset L\cup\left(  \mu,A_{k}\cup M_{k-1}\right]  \setminus M$ and
let%
\begin{align*}
\bar{F}  &  =D_{N_{q-1}}^{q-2}\cdots D_{N_{r}}^{r-1}D_{N_{r-2}\cup\mu}\cdots
D_{N_{j}\cup\mu}\cdots D_{N_{1}}\\
&  \hspace*{0.5in}R_{M_{p-1}}^{p}\cdots R_{M_{k}}^{k+1}R_{\left(  \mu
,A_{k}\cup M_{k-1}\right]  \setminus M}^{k}\cdots R_{M_{1}}\bar{E}.
\end{align*}
$\vspace{0.1in}$\newline When $1<i<k$ we have $\min A_{i-1}<\max
A_{i}\setminus\mu$ by assumption 2B2, and $\min A_{i}\setminus\mu<\max
A_{i+1}$ since $\mu\in A_{i}\cap M_{k-1}$ implies $\mu>\min A_{i}$. Next,
$\min A_{k-1}<\max A_{k}=\max A_{k}\cap M$ since $\max A_{k}<\mu\in A_{i},$
and $\min A_{k}\cap M<\mu=\max L.$ Finally, $\min L<\max A_{k+1}$ by
assumption 2B. Thus
\begin{align*}
\bar{a}\otimes\bar{b}  &  =A_{1}|\cdots|A_{i}\setminus\mu|\cdots|A_{k}\cap
M|L|\cdots|A_{p}\hspace*{1in}\\
&  \hspace*{1.5in}\otimes B_{q}|\cdots|B_{r+1}|B_{r-1}|\cdots|B_{j}\cup
\mu|\cdots|B_{1}%
\end{align*}
is an SCP; let $\bar{E}$ be the associated step matrix. Since $\min
M_{k-1}=\min M_{k-1}\setminus M=\mu>\max A_{k},$ both $R_{M_{k-1}\setminus\mu
}$ and $R_{(M_{k-1}\setminus\mu)\setminus M}$ are defined, so let
\begin{align*}
\bar{F}  &  =D_{N_{q-1}}^{q-2}\cdots D_{N_{r}}^{r-1}D_{N_{r-2}\cup\mu}\cdots
D_{N_{j}\cup\mu}\cdots D_{N_{1}}\\
&  \hspace*{0.25in}R_{M_{p-1}}^{p}\cdots R_{M_{k}}^{k+1}R_{\left(
M_{k-1}\setminus\mu\right)  \setminus M}^{k}R_{M_{k-1}\setminus\mu}\cdots
R_{M_{i}\setminus\mu}\cdots R_{M_{1}}\bar{E}.
\end{align*}
Then for all $i\leq k$, $\bar{u}\otimes v=c\left(  \bar{F}\right)  \otimes
r\left(  \bar{F}\right)  $ is a CP related to $\bar{a}\otimes\bar{b}$.
\vspace{0.1in}\newline\underline{\textit{Case 3}}: Assume $\left(  1\right)
\wedge\left(  2\right)  \wedge\left(  3\right)  ^{\prime}:m\in A_{k}\cap
B_{r},$ $\left(  m,M\right]  =\left(  m,A_{k}\cup M_{k-1}\right]  $ and
$N_{r-1}\neq\varnothing.$\vspace{0.1in}\newline Note that $M_{k}\subset\left(
A_{k}\cup M_{k-1}\right)  \setminus M=\left[  A_{k},m\right)  $ by conditions
(1) and (2) so that $U_{k}\setminus M=\left[  A_{k},m\right)  \setminus M_{k}%
$. Let $\nu=\min N_{r-1};$ then $\nu\in B_{i}\cap A_{j}$ for some$\,1\leq
i\leq r-1$ and
\[
j=k+\#[B_{i},\nu)+\sum_{s=i+1}^{r-1}\left(  \#B_{s}-1\right)  .
\]
\underline{\textit{Subcase 3A}}: \textit{Assume} $A_{j}=\nu.$\vspace
{0.1in}\newline In subcases 3A1 and 3A2, $\bar{u}$ is defined so that
\[
d_{\left[  A_{k},m\right)  }^{k+1}(\bar{u})\otimes v=d_{M}^{k}(u)\otimes v.
\]
\underline{\textit{Subcase 3A1}}: $j=k+1.$\vspace{0.1in}\newline Let%
\[
\bar{u}=U_{1}|\cdots|M|\left(  U_{k}\setminus M\right)  \cup U_{k+1}%
|U_{k+2}|\cdots|U_{p}.
\]
But $\nu>m$ since $\nu\in A_{k+1}\cap N_{r-1},$ consequently $M_{k}%
=\varnothing$ so that $U_{k}\setminus M=\left[  A_{k},m\right)  $ and
$U_{k+1}=A_{k+1}=\nu;$ thus $M_{k+1}=\varnothing.$ Clearly%
\begin{align*}
\bar{a}\otimes\bar{b}  &  =A_{1}|\cdots|A_{k}\cap M|\left[  A_{k},m\right)
\cup\nu|A_{k+2}|\cdots|A_{p}\hspace*{1in}\\
&  \hspace*{1.5in}\otimes\ B_{q}|\cdots|B_{r}\cup\nu|\cdots|B_{i}\setminus
\nu|\cdots|B_{1}%
\end{align*}
is an SCP; let $\bar{E}$ be the associated step matrix and let%
\[
\bar{F}=D_{N_{q-1}}\cdots D_{N_{r-1}\setminus\nu}\cdots D_{N_{i}\setminus\nu
}\cdots D_{N_{1}}R_{M_{p-1}}\cdots R_{\varnothing}^{k+1}R_{\varnothing}%
^{k}\cdots R_{M_{1}}\bar{E};
\]
then $\bar{u}\otimes v=c\left(  \bar{F}\right)  \otimes r\left(  \bar
{F}\right)  $ is a CP related to $\bar{a}\otimes\bar{b}.$\vspace
{0.1in}\newline\underline{\textit{Subcase 3A2}}: $j>k+1.$\vspace
{0.1in}\newline Let%
\[
\bar{u}=U_{1}|\cdots|M|U_{k}\setminus M|\cdots|U_{j-1}\cup U_{j}%
|U_{j+1}|\cdots|U_{p}.
\]
Again, $\nu>m$ implies that $M_{j-1}=\varnothing$ and $U_{j}=A_{j}=\nu$.
Clearly%
\begin{align*}
\bar{a}\otimes\bar{b}  &  =A_{1}|\cdots|A_{k}\cap M|\left[  A_{k},m\right)
\cup\nu|\cdots|A_{j-1}|A_{j+1}|\cdots|A_{p}\hspace*{0.5in}\\
&  \hspace*{1.5in}\otimes\ B_{q}|\cdots|B_{r}\cup\nu|\cdots|B_{i}\setminus
\nu|\cdots|B_{1}%
\end{align*}
is an SCP; let $\bar{E}$ be the associated step matrix and let
\begin{align*}
\bar{F}  &  =D_{N_{q-1}}\cdots D_{N_{r-1}\setminus\nu}\cdots D_{N_{i}%
\setminus\nu}\cdots D_{N_{1}}\\
&  \hspace*{0.25in}R_{M_{p-1}}\cdots R_{\varnothing}^{j}R_{M_{j-2}\cup\nu
}^{j-1}\cdots R_{M_{k+1}\cup\nu}^{k+2}R_{M_{k}\cup\nu}^{k+1}R_{\varnothing
}^{k}\cdots R_{M_{1}}\bar{E};
\end{align*}
then $\bar{u}\otimes v=c\left(  \bar{F}\right)  \otimes r\left(  \bar
{F}\right)  $ is a CP related to $\bar{a}\otimes\bar{b}.$ \vspace
{0.1in}\newline\underline{\textit{Subcase 3B}}: \textit{Assume} $A_{j}\neq
\nu.$ \vspace{0.1in}\newline Note that $i>1$ by assumption and let
\[
\bar{u}\otimes\bar{v}=U_{1}|\cdots|M|U_{k}\setminus M|\cdots|U_{p}\otimes
V_{q}|\cdots|V_{i}\cup V_{i-1}|\cdots|V_{1};
\]
then
\[
d_{M}^{k}(u)\otimes v=\bar{u}\otimes d_{i-1}^{V_{i-1}}\left(  \bar{v}\right)
.\vspace{0.1in}%
\]
Note that $\nu>m$ implies $M_{j-1}=\varnothing$ and $U_{j}=A_{j}\setminus
M_{j}$. Clearly%
\begin{align*}
\bar{a}\otimes\bar{b}  &  =A_{1}|\cdots|A_{k}\cap M|\left[  A_{k},m\right)
\cup\nu|\cdots|A_{j-1}|A_{j}\setminus\nu|A_{j+1}|\cdots|A_{p}\\
&  \hspace*{1.5in}\otimes\ B_{q}|\cdots|B_{r}\cup\nu|\cdots|(B_{i}\cup
B_{i-1})\setminus\nu|\cdots|B_{1}%
\end{align*}
is a SCP; let $\bar{E}$ be the associated step matrix and let
\begin{align*}
\bar{F}  &  =D_{N_{q-1}}\cdots D_{N_{r-1}\setminus\nu}\cdots D_{N_{i}%
\setminus\nu}\cdots D_{N_{1}}\\
&  \hspace*{0.25in}R_{M_{p-1}}^{p}\cdots R_{M_{j}}^{j+1}R_{M_{j-1}\cup\nu}%
^{j}\cdots R_{M_{k}\cup\nu}^{k+1}R_{\varnothing}^{k}\cdots R_{M_{1}}\bar{E};
\end{align*}
then $\bar{u}\otimes\bar{v}=c\left(  \bar{F}\right)  \otimes r\left(  \bar
{F}\right)  $ is a CP related to $\bar{a}\otimes\bar{b}.$\smallskip

For uniqueness of each pair $\bar{u}\otimes\bar{v}$ constructed above, note
the transformations $R$ and $D$ fix minimal elements, i.e., if $\bar{u}%
\otimes\bar{v}=R(\bar{a})\otimes D(\bar{b}),$ then necessarily $\min\bar
{U}_{i}=\min\bar{A}_{i}$ and $\min\bar{V}_{i}=\min\bar{B}_{i}$ for all $i;$ in
particular, if $R(\bar{a})=\tilde{R}(a^{\prime})$ or $D(\bar{b})=\tilde
{D}(b^{\prime})$ then $\min\bar{A}_{i}=\min{A}_{i}^{\prime}$ or $\min\bar
{B}_{i}=\min{B}_{i}^{\prime}.$ Consequently, for $d_{M}^{k}(u)\otimes v$ or
$u\otimes d_{\ell}^{N}(v)$ in the cases above, there is exactly one way to
construct a step matrix $\bar{E}$ so that $\bar{a}$ is step increasing and
$\bar{b}$ is step decreasing (it is straightforward to check that a
construction with distinct ${u}\otimes{v,}$ $\bar{u}\otimes\bar{v},$ and
$u^{\prime}\otimes v^{\prime}$ would contradict the necessary condition above
either for ${a}$ and ${a}^{\prime}$ or for ${b}$ and ${b}^{\prime}$). This
completes the proof.
\end{proof}

\section{Permutahedral Sets}

This section introduces the notion of a permutahedral set $\mathcal{Z}$, which
is a combinatorial object generated by permutahedra and equipped with
appropriate face and degeneracy operators. We construct the generating
category $\mathbf{P\ }$and show how to lift the diagonal on the permutahedra
$P$ constructed above to a diagonal on $\mathcal{Z}$. Naturally occurring
examples of permutahedral sets include the double cobar construction, i.e.,
Adams' cobar construction \cite{Adams} on the cobar with coassociative
coproduct \cite{Baues1}, \cite{CM}, \cite{KS1} (see Subsection \ref{sub4}
below). Permutahedral sets are distinguished from simplicial or cubical sets
by their higher order structure relations. While our construction of
$\mathbf{P\ }$follows the analogous (but not equivalent) construction for
polyhedral sets given by D.W. Jones in \cite{Jones}, there is no mention of
structure relations in \cite{Jones}.

\subsection{Singular Permutahedral Sets}

By way of motivation we begin with constructions of two singular permutahedral
sets--our universal examples. Whereas the first emphasizes coface and
codegeneracy operators, the second emphasizes cellular chains and is
appropriate for homology theory. We begin by constructing the various maps we
need to define singular coface and codegeneracy operators.

Fix a positive integer $n.$ For $0\leq p\leq n,$ let
\[
\underline{p}=\left\{
\begin{array}
[c]{cc}%
\varnothing, & p=0\\
\left\{  1,\ldots,p\right\}  , & 1\leq p\leq n
\end{array}
\right.  \text{ \ and \ }\overline{p}=\left\{
\begin{array}
[c]{cc}%
\varnothing, & p=0\\
\left\{  n-p+1,\ldots,n\right\}  , & 1\leq p\leq n;
\end{array}
\right.
\]
then $\underline{p}$ and $\overline{p}$ contain the first and last $p$
elements of $\underline{n},$ respectively; note that $\underline{p}%
\cap\overline{q}=\left\{  p\right\}  $ whenever $p+q=n+1$. Given integers
$r,s\in\underline{n}$ such that $r+s=n+1,$ there is a canonical projection
$\Delta_{r,s}:P_{n}\rightarrow P_{r}\times P_{s}$ whose restriction to a
vertex $v=a_{1}|\cdots|a_{n}\in P_{n}$ is given by
\[
\Delta_{r,s}(v)=b_{1}|\cdots|b_{r}\times c_{1}|\cdots|c_{s},
\]
where $\left(  b_{1},\ldots,b_{r};c_{1},\ldots,c_{k-1},c_{k+1},\ldots
,c_{s}\right)  \,$is the unshuffle of $\left(  a_{1},\ldots,a_{n}\right)  $
with $b_{i}\in\underline{r},$ $c_{j}\in\overline{s},$ $c_{k}=r.$ For example,
$\Delta_{2,3}(2|4|1|3)=2|1\times2|4|3$ and $\Delta_{3,2}(2|4|1|3)=2|1|3\times
4|3.$ Since the image of the vertices of a cell of $P_{n}$ uniquely determines
a cell in $P_{r}\times P_{s}$ the map $\Delta_{r,s}$ is well-defined and
cellular. Furthermore, the restriction of $\Delta_{r,s}$ to an $(n-k)$-cell
$A_{1}|\cdots|A_{k}\subset P_{n}$ is given by
\[
\Delta_{r,s}\left(  A_{1}|\cdots|A_{k}\right)  =\left\{
\begin{array}
[c]{ll}%
\underline{r}\times\left(  A_{1}|\cdots|A_{i}\setminus\underline{r-1}\text{
}|\cdots|A_{k}\right)  , & \hspace*{-0.3in}\text{if }\underline{r}\subseteq
A_{i},\text{ some }i,\\
& \\
\left(  A_{1}|\cdots|A_{j}\setminus\overline{s-1}\text{ }|\cdots|A_{k}\right)
\times\overline{s}, & \hspace*{-0.3in}\text{if }\overline{s}\subseteq
A_{j},\text{ some }j,\\
& \\%
\begin{array}
[c]{l}%
\hspace*{-0.06in}\left(  A_{1}\setminus\overline{s-1}\text{ }|\cdots
|A_{k}\setminus\overline{s-1}\right) \\
\hspace*{0.5in}\times\left(  A_{1}\setminus\underline{r-1}\text{ }%
|\cdots|A_{k}\setminus\underline{r-1}\right)  ,
\end{array}
& \text{otherwise.}%
\end{array}
\right.
\]
Note that $\Delta_{r,s}$ acts homeomorphically in the first two cases and
degeneratively in the third when $1<k<n$. When $n=3$ for example,
$\Delta_{2,2}$ maps the edge $1|23$ onto the edge $1|2\times23$ and the edge
$13|2$ onto the vertex $1|2\times3|2$ (see Figure 3).

\hspace*{0.3in}\setlength{\unitlength}{0.0005in}\begin{picture}
(2975,1000) \thicklines \put(1,239){\line( 1,0){1800}}
\put(1801,239){\line( 0,-1){1800}} \put(1801,-1561){\line(-1,
0){1800}} \put(1,-1561){\line( 0,1){1800}}
\put(1,239){\makebox(0,0){$\bullet$}}
\put(1,-661){\makebox(0,0){$\bullet$}}
\put(1,-1561){\makebox(0,0){$\bullet$}}
\put(1801,239){\makebox(0,0){$\bullet$}}
\put(1801,-661){\makebox(0,0){$\bullet$}}
\put(1801,-1561){\makebox(0,0){$\bullet$}}
\put(900,-680){\makebox(0,0){$123$}}
\put(-490,-1260){\makebox(0,0){${1|23}$}}
\put(950,530){\makebox(0,0){${3|12}$}}
\put(-490,-111){\makebox(0,0){${13|2}$}}
\put(2300,-111){\makebox(0,0){${23|1}$}}
\put(2300,-1260){\makebox(0,0){${2|13}$}}
\put(950,-1890){\makebox(0,0){${12|3}$}}
\put(3200,-400){\makebox(0,0){$\Delta_{2,2}$}}
\put(2700,-661){\vector(1,0){1000}}
\put(4700,-661){\makebox(0,0){${1|2} \times 23$}}
\put(8325,-661){\makebox(0,0){${2|1} \times 23$}}
\put(6500,530){\makebox(0,0){$12 \times {3|2}$}}
\put(6500,-1890){\makebox(0,0){$12 \times {2|3}$}}
\put(5601,239){\line( 1,0){1800}} \put(7401,239){\line(
0,-1){1800}} \put(7401,-1561){\line(-1, 0){1800}}
\put(5601,-1561){\line( 0,1){1800}}
\put(5601,239){\makebox(0,0){$\bullet$}}
\put(5601,-1561){\makebox(0,0){$\bullet$}}
\put(7401,239){\makebox(0,0){$\bullet$}}
\put(7401,-1561){\makebox(0,0){$\bullet$}}
\put(6500,-661){\makebox(0,0){$12 \times 23$}}
\end{picture}$\vspace{1.1in}$

\begin{center}
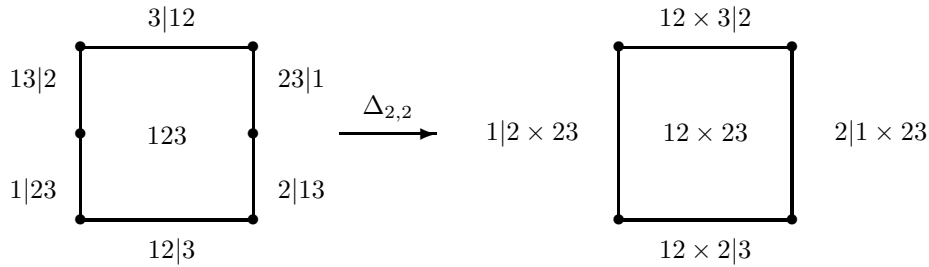
Figure 3: The projection $\Delta_{2,2}:P_{3}\rightarrow I^{2}.\vspace{0.2in} $
\end{center}

Now identify the set $U=\left\{  u_{1}<\cdots<u_{n}\right\}  $ with $P_{n}$
and the ordered partitions of $U$ with the faces of $P_{n}$ in the obvious
way. Then $\left(  \Delta_{r,s}\times1\right)  \circ\Delta_{r+s-1,t}=\left(
1\times\Delta_{s,t}\right)  \circ\Delta_{r,s+t-1}$ whenever $r+s+t=n+2$ so
that $\Delta_{\ast,\ast}$ acts coassociatively with respect to Cartesian
product. It follows that each $k$-tuple $\left(  n_{1},\ldots,n_{k}\right)
\in\mathbb{N}^{k}$ with $k\geq2$ and $n_{1}+\cdots+n_{k}=n+k-1$ uniquely
determines a cellular projection $\Delta_{n_{1}\cdots n_{k}}:P_{n}\rightarrow
P_{n_{1}}\times\cdots\times P_{n_{k}} $ given by the composition%
\[
\Delta_{n_{1}\cdots n_{k}}=\left(  \Delta_{n_{1},n_{2}}\times1^{\times
k-2}\right)  \circ\cdots\circ\left(  \Delta_{n_{\left(  k-2\right)
}-k+3,n_{k-1}}\times1\right)  \circ\Delta_{n_{\left(  k-1\right)  }-k+2,n_{k}%
},
\]
where $n_{\left(  q\right)  }=n_{1}+\cdots+n_{q};$ and in particular,
\begin{equation}
\Delta_{n_{1}\cdots n_{k}}\left(  \underline{n}\right)  =\underline{n_{1}%
}\times\underline{n_{\left(  2\right)  }-1}\setminus\underline{n_{1}-1}%
\times\cdots\times\underline{n_{\left(  k\right)  }-\left(  k-1\right)
}\setminus\underline{n_{\left(  k-1\right)  }-\left(  k-1\right)  }.
\label{proj}%
\end{equation}
Note that formula \ref{proj} with $k=n-1$ and $n_{i}=2$ for all $i$ defines a
projection $\rho_{n}:P_{n}\rightarrow I^{n-1}$%
\[
\rho_{n}\left(  \underline{n}\right)  =\Delta_{2\cdots2}\left(  \underline
{n}\right)  =12\times23\times\cdots\times\left\{  n-1,n\right\}
\]
(see Figure 4) acting on a vertex $u=u_{1}|\cdots|u_{n}$ as follows: For each
$i\in\underline{n-1},$ let $\left\{  u_{j},u_{k}\text{ }|\text{ }j<k\right\}
=\left\{  u_{1},\ldots,u_{n}\right\}  \cap\left\{  i,i+1\right\}  $ and set
$v_{i}=u_{j},$ $v_{i+1}=u_{k};$ then $\rho_{n}(u)=v_{1}|v_{2}\times
\cdots\times v_{n-1}|v_{n}.$\vspace{0.2in}

\hspace*{0.4in}\setlength{\unitlength}{0.00015in}\begin{picture}
(7500,7500) \thicklines
\put(6000,4800){\line( 0,-1){4800}}
\put(6000,4800){\makebox(0,0){$\bullet$}}
\put(6000,2400){\makebox(0,0){$\bullet$}}
\put(6000,0){\makebox(0,0){$\bullet$}}
\put(6000,3600){\makebox(0,0){$\bullet$}}
\put(6000,0){\line( 1, 0){4800}}
\put(10800,0){\line( 0, 1){4800}}
\put(10800,4800){\line(-1, 0){4800}}
\put(10800,4800){\makebox(0,0){$\bullet$}}
\put(10800,3600){\makebox(0,0){$\bullet$}}
\put(10800,2400){\makebox(0,0){$\bullet$}}
\put(10800,0){\makebox(0,0){$\bullet$}}
\put(6000,2400){\line( 1, 0){4800}}
\put(3000,6800){\line( 0,-1){4800}}
\put(3000,6800){\makebox(0,0){$\bullet$}}
\put(3000,5600){\makebox(0,0){$\bullet$}}
\put(3000,4400){\makebox(0,0){$\bullet$}}
\put(3000,2000){\makebox(0,0){$\bullet$}}
\put(3000,6800){\line( 1, 0){4800}}
\put(3000,5600){\line( 1, 0){1200}}
\put(5000,5600){\line(1,0){2800}}
\put(3000,2000){\line( 1, 0){1200}}
\put(4800,2000){\line(1,0){900}}
\put(6300,2000){\line( 1, 0){1500}}
\put(7800,5000){\line( 0,1){1700}}
\put(7800,2600){\line(0,1){2050}}
\put(7800,2000){\line(0,1){200}}
\put(7800,6800){\makebox(0,0){$\bullet$}}
\put(7800,5600){\makebox(0,0){$\bullet$}}
\put(7800,4400){\makebox(0,0){$\bullet$}}
\put(7800,2000){\makebox(0,0){$\bullet$}}
\put(3000,2000){\line( 3,-2){3000}}
\put(6000,4800){\line(-3,2){3000}}
\put(4500,5800){\line( 0,-1){4800}}
\put(4500,5800){\makebox(0,0){$\bullet$}}
\put(4500,4600){\makebox(0,0){$\bullet$}}
\put(4500,3400){\makebox(0,0){$\bullet$}}
\put(4500,1000){\makebox(0,0){$\bullet$}}
\put(6000,3600){\line(-3, 2){1500}}
\put(4500,3400){\line(-3, 2){1500}}
\put(9300,3400){\line(-3, 2){1500}}
\put(9300,5800){\line( 0,-1){800}}
\put(9300,4600){\line(0,-1){2000}}
\put(9300,2180){\line(0,-1){1200}}
\put(9370,5800){\makebox(0,0){$\bullet$ }}
\put(9300,3400){\makebox(0,0){$\bullet$}}
\put(9300,1000){\makebox(0,0){$\bullet$}}
\put(9300,4600){\makebox(0,0){$\bullet$}}
\put(10800,3600){\line(-3, 2){1500}}
\put(7800,2000){\line(3, -2){3000}}
\put(7800,6800){\line( 3,-2){3000}}
\put(17000,4800){\line( 0,-1){4800}}
\put(17000,4800){\makebox(0,0){$\bullet$}}
\put(17000,2400){\makebox(0,0){$\bullet$}}
\put(17000,0){\makebox(0,0){$\bullet$}}
\put(17000,0){\line( 1, 0){4800}}
\put(21800,0){\line( 0, 1){4800}}
\put(21800,4800){\line(-1, 0){4800}}
\put(21800,4800){\makebox(0,0){$\bullet$}}
\put(21800,2400){\makebox(0,0){$\bullet$}}
\put(21800,0){\makebox(0,0){$\bullet$}}
\put(17000,2400){\line( 1, 0){4800}}
\put(14000,6800){\line( 0,-1){4800}}
\put(14000,6800){\makebox(0,0){$\bullet$}}
\put(14000,4400){\makebox(0,0){$\bullet$}}
\put(14000,2000){\makebox(0,0){$\bullet$}}
\put(14000,6800){\line( 1, 0){4800}}
\put(14000,4400){\line( 1, 0){2600}}
\put(17300,4400){\line( 1,0){1500}}
\put(14000,2000){\line( 1, 0){2600}}
\put(17300,2000){\line( 1,0){1500}}
\put(18800,5000){\line( 0,1){1700}}
\put(18800,2600){\line(0,1){2050}}
\put(18800,2000){\line( 0,1){200}}
\put(18800,6800){\makebox(0,0){$\bullet$}}
\put(18800,4400){\makebox(0,0){$\bullet$}}
\put(18800,2000){\makebox(0,0){$\bullet$}}
\put(14000,2000){\line( 3,-2){3000}}
\put(17000,4800){\line(-3,2){3000}}
\put(18800,2000){\line(3, -2){3000}}
\put(18800,6800){\line( 3,-2){3000}}
\put(6000,-11200){\line( 0,-1){4800}}
\put(6000,-11200){\makebox(0,0){$\bullet$}}
\put(6000,-16000){\makebox(0,0){$\bullet$}}
\put(6000,-16000){\line( 1, 0){4800}}
\put(10800,-16000){\line( 0, 1){4800}}
\put(10800,-11200){\line(-1, 0){4800}}
\put(10800,-11200){\makebox(0,0){$\bullet$}}
\put(10800,-16000){\makebox(0,0){$\bullet$}}
\put(3000,-9200){\line( 0,-1){4800}}
\put(3000,-9200){\makebox(0,0){$\bullet$}}
\put(3000,-14000){\makebox(0,0){$\bullet$}}
\put(3000,-9200){\line( 1, 0){4800}}
\put(3000,-14000){\line( 1, 0){1200}}
\put(4800,-14000){\line(1,0){900}}
\put(6300,-14000){\line( 1, 0){1500}}
\put(7800,-11000){\line( 0,1){1700}}
\put(7800,-14000){\line(0,1){2600}}
\put(7800,-9200){\makebox(0,0){$\bullet$}}
\put(7800,-14000){\makebox(0,0){$\bullet$}}
\put(3000,-14000){\line( 3,-2){3000}}
\put(6000,-11200){\line(-3,2){3000}}
\put(4500,-10200){\line( 0,-1){4800}}
\put(4500,-10200){\makebox(0,0){$\bullet$}}
\put(4500,-15050){\makebox(0,0){$\bullet$}}
\put(9300,-10200){\line( 0,-1){800}}
\put(9300,-11400){\line(0,-1){3600}}
\put(9370,-10200){\makebox(0,0){$\bullet$ }}
\put(9300,-15020){\makebox(0,0){$\bullet$}}
\put(7800,-14000){\line(3, -2){3000}}
\put(7800,-9200){\line( 3,-2){3000}}
\put(17000,-11200){\line( 0,-1){4800}}
\put(17000,-11200){\makebox(0,0){$\bullet$}}
\put(17000,-16000){\makebox(0,0){$\bullet$}}
\put(17000,-16000){\line( 1, 0){4800}}
\put(21800,-16000){\line( 0, 1){4800}}
\put(21800,-11200){\line(-1, 0){4800}}
\put(21800,-11200){\makebox(0,0){$\bullet$}}
\put(21800,-16000){\makebox(0,0){$\bullet$}}
\put(14000,-9200){\line( 0,-1){4800}}
\put(14000,-9200){\makebox(0,0){$\bullet$}}
\put(14000,-14000){\makebox(0,0){$\bullet$}}
\put(14000,-9200){\line( 1, 0){4800}}
\put(14000,-14000){\line( 1, 0){2600}}
\put(17300,-14000){\line( 1,0){1500}}
\put(18800,-11000){\line( 0,1){1700}}
\put(18800,-14000){\line(0,1){2650}}
\put(18800,-9200){\makebox(0,0){$\bullet$}}
\put(18800,-14000){\makebox(0,0){$\bullet$}}
\put(14000,-14000){\line( 3,-2){3000}}
\put(14000,-9200){\line( 3,-2){3000}}
\put(18800,-14000){\line(3, -2){3000}}
\put(18800,-9200){\line( 3,-2){3000}}
\put(5750,-4300){\makebox(0,0){$\Delta_{3,2}$}}
\put(7700,-3400){\vector(0,-1){2000}}
\put(7700,-2700){\makebox(0,0){$1234$}}
\put(7700,-6300){\makebox(0,0){$123\times34$}}
\put(9500,-2700){\vector(1,0){4000}}
\put(9900,-6300){\vector(1,0){3200}}
\put(16000,-2700){\makebox(0,0){$12\times234$}}
\put(16000,-6300){\makebox(0,0){$12\times23\times34$}}
\put(16000,-3400){\vector(0,-1){2000}}
\put(18500,-4300){\makebox(0,0){$1\times\Delta_{2,2}$}}
\put(11500,-7600){\makebox(0,0){$\Delta_{2,2}\times1$}}
\put(11500,-1500){\makebox(0,0){$\Delta_{2,3}$}}
\end{picture}$\vspace{2.6in}$

\begin{center}
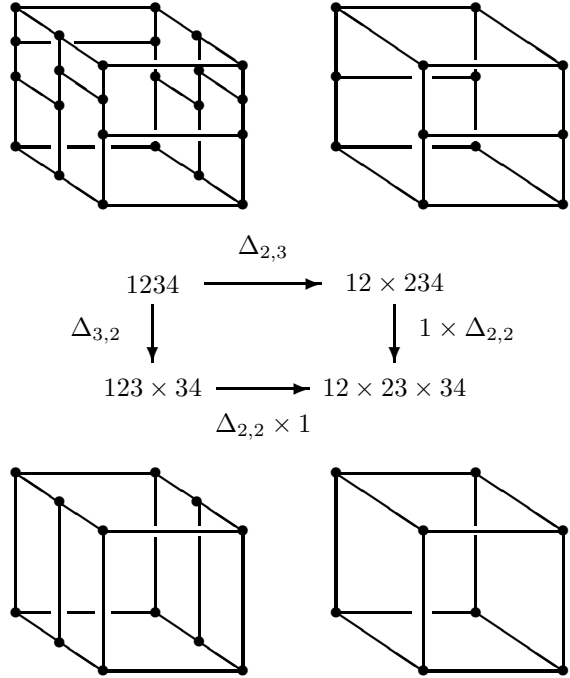
Figure 4: The projection $\rho_{4}:P_{4}\rightarrow I^{3}.\vspace{0.2in}$
\end{center}

Now choose a (non-cellular) homeomorphism $\gamma_{n}:I^{n-1}\rightarrow
P_{n}$ whose restriction to a vertex $v=v_{1}|v_{2}\times\cdots\times
v_{n-1}|v_{n}$ can be expressed inductively as follows: Set $A_{2}=v_{1}%
|v_{2};$ if $A_{k-1}$ has been obtained from $v_{1}|v_{2}\times\cdots\times
v_{k-2}|v_{k-1},$ set%
\[
A_{k}=\left\{
\begin{array}
[c]{ll}%
A_{k-1}|k, & \text{if }v_{k}=k,\\
k|A_{k-1}, & \text{otherwise.}%
\end{array}
\right.
\]
For example, $\gamma_{4}\left(  2|1\times3|2\times3|4\right)  =3|2|1|4.$ Then
$\gamma_{n}$ sends the vertices of $I^{n-1}$ to cubical vertices of $P_{n}$
and the vertices of $P_{n}$ fixed by $\gamma_{n}\rho_{n}$ are exactly its
cubical vertices. Given a codimension 1 face $A|B\subset P_{n},$ index the
elements of $A$ and $B$ as follows: If $n\in A,$ write $A=\left\{
a_{1}<\cdots<a_{m}\right\}  $ and $B=\left\{  b_{1}<\cdots<b_{\ell}\right\}
;$ if $n\in B,$ write $A=\left\{  a_{1}<\cdots<a_{\ell}\right\}  $ and
$B=\left\{  b_{1}<\cdots<b_{m}\right\}  .$ Then $A|B$ uniquely embeds in
$P_{n}$ as the subcomplex%
\[
P_{\ell}\times P_{m}=\left\{
\begin{array}
[c]{ll}%
a_{1}|\cdots|a_{m}|B\times A|b_{1}|\cdots|b_{\ell}, & \text{if }n\in A\\
A|b_{1}|\cdots|b_{m}\times a_{1}|\cdots|a_{\ell}|B, & \text{if }n\in B.
\end{array}
\right.
\]
For example, $14|23$ embeds in $P_{4}$ as $1|4|23\times14|2|3$. Let
$\iota_{A|B}:A|B\hookrightarrow P_{\ell}\times P_{m}$ denote this embedding
and let $h_{A|B}=\iota_{A|B}^{-1};$ then $h_{A|B}:P_{\ell}\times
P_{m}\rightarrow A|B$ is an orientation preserving homeomorphism. Also define
the cellular projection
\[
\phi_{A|B}:P_{n}\rightarrow P_{\ell}\times P_{m}=\left\{
\begin{array}
[c]{ll}%
b_{1}\cdots b_{\ell}\times a_{1}\cdots a_{m}, & \text{if }n\in A\\
a_{1}\cdots a_{\ell}\times b_{1}\cdots b_{m}, & \text{if }n\in B
\end{array}
\right.
\]
on a vertex $c=c_{1}|\cdots|c_{n}$ by $\phi_{A|B}\left(  c\right)
=u_{1}|\cdots|u_{\ell}\times v_{1}|\cdots|v_{m},$ where $(u_{1},\ldots
,u_{\ell};$ $v_{1},\ldots,v_{m})\,$is the unshuffle of $\left(  c_{1}%
,\ldots,c_{n}\right)  $ with $u_{i}\in B,$ $v_{j}\in A$ when $n\in A$ or with
$u_{i}\in A,$ $v_{j}\in B$ when $n\in B.$ Note that unlike $\Delta_{r,s},$ the
projection $\phi_{A|B}$ always degenerates on the top cell; furthermore,
$\phi_{A|B}\circ h_{A|B}=\phi_{B|A}\circ h_{A|B}=1$. We note that when $A$ or
$B$ is a singleton set, the projection $\phi_{A|B}$ was defined by R.J.
Milgram in \cite{Milgram}.

The \textit{singular codegeneracy operator associated with }$A|B$ is the map
$\beta_{A|B}:P_{n}\rightarrow P_{n-1}$ given by the composition
\[
P_{n}\overset{\phi_{A|B}}{\longrightarrow}P_{\ell}\times P_{m}\overset
{\rho_{\ell}\times\rho_{m}}{\longrightarrow}I^{\ell-1}\times I^{m-1}%
=I^{n-2}\overset{\gamma_{n-1}}{\longrightarrow}P_{n-1};
\]
the \textit{singular coface operator associated with }$A|B$ is the map
$\delta_{A|B}:P_{n-1}\rightarrow P_{n}$ given by the composition%
\[
P_{n-1}\overset{\rho_{n-1}}{\longrightarrow}I^{n-2}=I^{\ell-1}\times
I^{m-1}\overset{\gamma_{\ell}\times\gamma_{m}}{\longrightarrow}P_{\ell}\times
P_{m}\overset{h_{A|B}}{\longrightarrow}A|B\overset{i}{\hookrightarrow}P_{n}.
\]
Unlike the simplicial or cubical case, $\delta_{A|B}$ need not be injective.
We shall often abuse notation and write $h_{A|B}:P_{\ell}\times P_{m}%
\rightarrow P_{n}$ when we mean $i\circ h_{A|B}.$

We are ready to define our first universal example. For future reference and
to emphasize the fact that our definition depends only on positive integers,
let $\left(  n_{1},\ldots,n_{k}\right)  \in\mathbb{N}^{k}$ such that
$n_{\left(  k\right)  }=n$ and denote%
\[
{\mathcal{P}}_{n_{1}\cdots n_{k}}\left(  n\right)
=\{\text{\textit{Partitions} }A_{1}|\cdots|A_{k}\text{ \textit{of }}%
\underline{n}\text{ }|\text{ }\#A_{i}=n_{i}\}.
\]

\begin{definition}
Let $Y$ be a topological space. The \underline{singular permutahedral set of
$Y$} consists of the singular set
\[
{Sing}_{\ast}^{P}Y=\bigcup\limits_{n\geq1}\left[  {Sing}_{n}^{P}Y=\left\{
\text{Continuous maps }P_{n}{\rightarrow}Y\right\}  \right]
\]
together with singular face and degeneracy operators
\[
d_{A|B}:Sing_{n}^{P}Y\rightarrow Sing_{n-1}^{P}Y\text{ \ and \ }\varrho
_{A|B}:Sing_{n-1}^{P}Y\rightarrow Sing_{n}^{P}Y
\]
defined respectively for each $n\geq2$ and $A|B\in{\mathcal{P}}_{\ast\ast
}\left(  n\right)  $ as the pullback along $\delta_{A|B}$ and $\beta_{A|B}$,
i.e., for $f\in Sing_{n}^{P}Y$ and $g\in Sing_{n-1}^{P}Y,$
\[
d_{A|B}(f)=f\circ\delta_{A|B}\text{ \ and \ }\varrho_{A|B}(g)=g\circ
\beta_{A|B}.
\]
\vspace*{0.1in}
\end{definition}

\hspace{1.3in}\setlength{\unitlength}{0.0002in}\begin{picture}
(0,0) \thicklines
\put(1500,-600){\vector(3,-1){8000}}
\put(11200,-1000){\vector(0,-1){1500}}
\put(4600,0){\makebox(0,0){$\delta_{A|B}:P_{n-1} \rightarrow I^{n-1} \rightarrow P_{\ell} \times P_m \rightarrow A|B \hookrightarrow P_n$}}
\put(11200,-3500){\makebox(0,0){$Y$}}
\put(12000,-1500){\makebox(0,0){$f$}}
\put(4000,-2800){\makebox(0,0){$d_{A|B}(f)$}}
\end{picture}\vspace*{0.8in}

\begin{center}
Figure 5:\ The singular face operator associated with $A|B$.\vspace{0.2in}
\end{center}

Although coface operators $\delta_{A|B}:P_{n-1}\rightarrow P_{n}$ need not be
inclusions, the top cell of $P_{n-1}$ is always non-degenerate (c.f.
Definition \ref{degen}); however, the top cell of $P_{n-2}$ may degenerate
under quadratic compositions $\delta_{A|B}\delta_{C|D}:P_{n-2}\rightarrow
P_{n}$ . For example, $\delta_{12|34}\delta_{13|2}:P_{2}\rightarrow P_{4}$ is
a constant map, since $\delta_{12|34}:P_{3}\rightarrow P_{2}\times
P_{2}\hookrightarrow P_{4}$\ sends the edge\thinspace$13|2$ to the vertex
$1|2\times3|2$.

\begin{definition}
\label{admissible} \textit{A quadratic composition of face operators }%
$d_{C|D}d_{A|B}$ \underline{acts on $P_{n}$} \textit{if the top cell of
}$P_{n-2}$\textit{\ is non-degenerate under the composition}
\[
\delta_{A|B}\delta_{C|D}:P_{n-2}\rightarrow P_{n}.
\]

\end{definition}

\noindent Theorem \ref{faces} below gives the conditions under which a
quadratic composition acts on $P_{n}$. For comparison, quadratic compositions
of simplicial or cubical face operators always act on the simplex or cube.
When $d_{C|D}d_{A|B}$ acts on $P_{n},$ we assign the label $d_{C|D}d_{A|B}$ to
the codimension 2 face $\delta_{A|B}\delta_{C|D}\left(  \underline{n}\right)
$. The various paths of descent from the top cell to a cell in codimension 2
gives rise to relations among compositions of face and degeneracy operators
(see Figure 6). \vspace{0.2in}

\setlength{\unitlength}{0.0005in}\begin{picture} (2975,1000)
\thicklines\put(3601,239){\line( 1, 0){1800}}
\put(5401,239){\line( 0,-1){1800}} \put(5401,-1561){\line(-1,
0){1800}} \put(3601,-1561){\line( 0, 1){1800}}
\put(3601,239){\makebox(0,0){$\bullet$}}
\put(3601,-661){\makebox(0,0){$\bullet$}}
\put(3601,-1561){\makebox(0,0){$\bullet$}}
\put(5401,239){\makebox(0,0){$\bullet$}} \put
(5401,-661){\makebox(0,0){$\bullet$}}
\put(5401,-1561){\makebox(0,0){$\bullet$}}
\put(4500,-680){\makebox(0,0){$123$}}
\put(2000,-1861){\makebox(0,0){$d_{1|2}d_{12|3}=d_{1|2}d_{1|23}$}}
\put(2000,-699){\makebox(0,0){$d_{1|2}d_{13|2}=d_{2|1}d_{1|23}$}}
\put(2000,464){\makebox(0,0){$d_{2|1}d_{13|2}=d_{1|2}d_{3|12}$}}
\put(7000,-1861){\makebox(0,0){$d_{1|2}d_{2|13}=d_{2|1}d_{12|3}$}}
\put(7000,-699){\makebox(0,0){$d_{2|1}d_{2|13}=d_{1|2}d_{23|1}$}}
\put(7000,464){\makebox(0,0){$d_{2|1}d_{23|1}=d_{2|1}d_{3|12}$}}
\put(3110,-1260){\makebox(0,0){$d_{1|23}$}}
\put(4550,530){\makebox(0,0){$d_{3|12}$}}
\put(3110,-111){\makebox(0,0){$d_{13|2}$}}
\put(5900,-111){\makebox(0,0){$d_{23|1}$}}
\put(5900,-1260){\makebox(0,0){$d_{2|13}$}}
\put(4550,-1890){\makebox(0,0){$d_{12|3}$}}
\end{picture}\vspace{1.1in}

\begin{center}
Figure 6: Quadratic relations on the vertices of $P_{3}.$\vspace{0.2in}
\end{center}

It is interesting to note that singular permutahedral sets have higher order
structure relations, an example of which appears below in Figure 7 (see also
(\ref{hrelation})). This distinguishes permutahedral sets from simplicial or
cubical sets in which relations are strictly quadratic. Our second universal
example, called a \textquotedblleft singular multipermutahedral
set,\textquotedblright\ specifies a singular permutahedral set by restricting
to maps $f=\bar{f}\circ\Delta_{n_{1}\cdots n_{k}}$ for some continuous
$\bar{f}:P_{n_{1}}\times\cdots\times P_{n_{k}}\rightarrow Y$. Face and
degeneracy operators satisfy those relations above in which $\Delta
_{n_{1}\cdots n_{k}}$ plays no essential role.

\vspace*{0.7in} \hspace*{0.09in}%
\setlength{\unitlength}{0.000236in}\begin{picture} (0,0)
\thicklines
\put(800,0){\vector(1,0){1500}}
\put(4200,0){\vector(1,0){6000}}
\put(12200,0){\vector(1,0){6000}}
\put(0,800){\line(0,1){1000}}
\put(19200,1750){\vector(0,-1){1000}}
\put(0,1750){\line(1,0){19225}}
\put(19200,-600){\vector(0,-1){4000}}
\put(18400,-5500){\vector(-1,0){3000}}
\put(18800,-600){\vector(-1,-1){4000}}
\put(11400,-600){\vector(2,-3){2650}}
\put(4000,-500){\vector(2,-1){9200}}
\put(0,-5500){\vector(1,0){13200}}
\put(0,-800){\line(0,-1){4720}}
\put(0,0){\makebox(0,0){$P_2$}}
\put(3200,0){\makebox(0,0){$P_1$}}
\put(7400,900){\makebox(0,0){$\delta_{1|2}$}}
\put(11100,0){\makebox(0,0){$P_2$}}
\put(15300,900){\makebox(0,0){$\delta_{1|23}$}}
\put(19200,0){\makebox(0,0){$P_3$}}
\put(9621,2600){\makebox(0,0){$\delta_{13|2}$}}
\put(19200,0){\makebox(0,0){$P_3$}}
\put(18300,-3100){\makebox(0,0){$\delta_{12|34}$}}
\put(19200,-5500){\makebox(0,0){$P_4$}}
\put(14300,-5500){\makebox(0,0){$Y$}}
\put(1500,900){\makebox(0,0){$\beta_{2|1}$}}
\put(17000,-4800){\makebox(0,0){$f$}}
\put(10000,-1700){\makebox(0,0){$d_{1|23}d_{12|34}(f)$}}
\put(4800,-2700){\makebox(0,0){$d_{2|1}d_{1|23}d_{12|34}(f)$}}
\put(16000,-1700){\makebox(0,0){$d_{12|34}(f)$}}
\put(5500,-4800){\makebox(0,0){$\varrho_{2|1}d_{2|1}d_{1|23}d_{12|34}(f)=d_{13|2}d_{12|34}(f)$}}
\end{picture}\vspace*{1.4in}

\begin{center}
Figure 7: A quartic relation in ${Sing}_{\ast}^{P}Y$.\vspace*{0.2in}
\end{center}

Once again, fix a positive integer $n,$ but this time consider $\left(
n_{1},\ldots,n_{k}\right)  \in\left(  \mathbb{N}\cup0\right)  ^{k}$ with
$n_{(k)}=n-1$ and the projection $\Delta_{n_{1}+1\cdots n_{k}+1}%
:P_{n}\rightarrow P_{n_{1}+1}\times\cdots\times P_{n_{k}+1}$ with $\Delta
_{n}:P_{n}\rightarrow P_{n}$ defined to be the identity. Given a topological
space $Y,$ let%
\[
{Sing}^{n_{1}\cdots n_{k}}Y=\left\{  \bar{f}\circ\Delta_{n_{1}+1\cdots
n_{k}+1}:P_{n}\rightarrow Y\text{ }|\text{\ }\bar{f}\ \text{\textit{is
continuous}}\right\}  ;
\]
define $f,f^{\prime}\in{Sing}^{n_{1}\cdots n_{k}}Y$ to be equivalent if there
exists $g:P_{n_{1}+1}\times\cdots\times P_{n_{i-1}+1}\times P_{1}\times
P_{n_{i+1}+1}\times\cdots\times P_{n_{k}+1}\rightarrow Y$ for some $i<k$ such
that
\[
f=g\circ(1^{\times i-1}\times\phi_{\underline{n_{i}+1}|n_{i}+1}\times1^{\times
k-i-1})\circ\Delta_{n_{1}+1\cdots n_{i-1}+1,n_{i}+2,n_{i+2}+1\cdots n_{k}+1}%
\]
and
\[
f^{\prime}=g\circ(1^{\times i}\times\phi_{1|\underline{n_{i+2}+1}\setminus
1}\times1^{\times k-i-2})\circ\Delta_{n_{1}+1\cdots n_{i}+1,n_{i+2}%
+2,n_{i+3}+1\cdots n_{k}+1},
\]
in which case we write $f\sim f^{\prime}.$ The geometry of the cube motivates
this equivalence; the degeneracies in the product of cubical sets implies the
identification (c.f. \cite{Kan} or the definition of the cubical set functor
$\mathbf{\Omega}X$ in the Appendix).

Define the singular set
\[
{Sing}_{n}^{M}Y=\bigcup_{\substack{\left(  n_{1},\ldots,n_{k}\right)
\in\left(  \mathbb{N}\cup0\right)  ^{k}\\n_{(k)}=n-1}}\left.  {Sing}%
^{n_{1}\cdots n_{k}}Y\right/  \sim.
\]
Singular face and degeneracy operators
\[
d_{A|B}:Sing_{n}^{M}Y\rightarrow Sing_{n-1}^{M}Y\text{ \ and \ }\varrho
_{A|B}:Sing_{n-1}^{M}Y\rightarrow Sing_{n}^{M}Y
\]
are defined piece-wise for each $n\geq2$ and $A|B\in{\mathcal{P}}_{\ast,\ast
}\left(  n\right)  ,$ depending on the form of $A|B.$ More precisely, for each
pair of integers $\left(  p_{i},q_{i}\right)  ,$ $1\leq i\leq k,$ with%
\[
p_{i}=1+\sum_{j=1}^{i-1}n_{j}\text{ \ and \ }q_{i}=1+\sum_{j=i+1}^{k}%
n_{j},\text{ \ let}%
\]%
\[
\mathcal{Q}_{p_{i},q_{i}}\left(  n\right)  =\left\{  U|V\in{\mathcal{P}}%
_{\ast,\ast}\left(  n\right)  \text{ }|\text{ }\left(  \underline{p_{i}%
}\subseteq U\text{ or }\underline{p_{i}}\subseteq V\right)  \text{ and
}\left(  \overline{q_{i}}\subseteq U\text{ or }\overline{q_{i}}\subseteq
V\right)  \right\}  ;\vspace*{0.1in}%
\]
in particular, when $r+s=n+1,$ set $k=2,$ $p_{1}=q_{2}=1,$ $p_{2}=r$ and
$q_{1}=s,$ then%
\[
\mathcal{Q}_{r,1}\left(  n\right)  =\left\{  U|V\in{\mathcal{P}}_{\ast,\ast
}\left(  n\right)  \text{ }|\text{ }\underline{r}\subseteq U\text{ or
}\underline{r}\subseteq V\right\}  \text{ and}%
\]%
\[
\mathcal{Q}_{1,s}\left(  n\right)  =\left\{  U|V\in{\mathcal{P}}_{\ast,\ast
}\left(  n\right)  \text{ }|\text{ }\overline{s}\subseteq U\text{ or
}\overline{s}\subseteq V\right\}  .
\]
Since we identify $\underline{r}|\overline{s}\subset P_{n+1}$ with
$P_{r}\times P_{s}=\Delta_{r,s}\left(  P_{n}\right)  ,$ it follows that
$A|B\in\mathcal{Q}_{p_{i},q_{i}}\left(  n\right)  $ for some $i$ if and only
if $\delta_{A|B}\delta_{\underline{r}|\overline{s}}:P_{n-1}\rightarrow
P_{n+1}$ is non-degenerate; consequently we consider cases $A|B\in
\mathcal{Q}_{p_{i},q_{i}}\left(  n\right)  $ for some $i$ and $A|B\notin
\mathcal{Q}_{p_{i},q_{i}}\left(  n\right)  $ for all $i.$

Since our definitions of $d_{A|B}$ and $\varrho_{A|B}$ are independent in the
first case and interdependent in the second, we define both operators
simultaneously. But first we need some notation: Given an increasingly ordered
set $M=\{m_{1}<\cdots<m_{k}\}\subset{\mathbb{N}},$ let $I_{M}:M\rightarrow
\underline{\#M}$ denote the \emph{indexing map} $m_{i}\mapsto i$ and let
$M+z=\{m_{i}+z\}$ denote \emph{translation by} $z\in\mathbb{Z}$. Of course,
$M-z$ and $M+z$ are left and right translations when $z>0;$ we adopt the
convention that translation takes preference over set operations.

Assume $A|B\in{\mathcal{Q}}_{p_{i},q_{i}}(n)$ for some $i,$ and let
\[
C_{i}=\left\{  p_{i},p_{i}+1,...,p_{i}+n_{i}\right\}  ;
\]%
\[
A_{i}=\left(  C_{i}\cap A\right)  -n_{(i-1)},\text{ \ }B_{i}=\left(  C_{i}\cap
B\right)  -n_{(i-1)};
\]%
\begin{equation}
n_{i}^{\prime}=\#(A\cap C_{i})-1,\text{ \ }\,n_{i}^{\prime\prime}=\#(B\cap
C_{i})-1. \label{formal}%
\end{equation}
For example, $n=6,$ $n_{1}=3$ and $n_{2}=2$ determines the projection
$\Delta_{4,3}:P_{6}\rightarrow1234\times456$ and pairs $\left(  p_{1}%
,q_{1}\right)  =\left(  1,3\right)  $ and $\left(  p_{2},q_{2}\right)
=\left(  4,1\right)  .$ Thus $A|B=1234|56\in{\mathcal{Q}}_{3,2}(6)$ and the
composition $\delta_{\underline{4}|\overline{3}}\delta_{A|B}:P_{5}\rightarrow
P_{7}$ is non-degenerate. Furthermore, $C_{2}=456,$ $A_{2}=\left(
456\cap1234\right)  -3=1,$ $B_{2}=23,$ $n_{i}^{\prime}=0,$ $n_{i}%
^{\prime\prime}=1$ and we may think of $d_{A|B}$ acting on $1234\times456$ as
$1\times d_{1|23}.$

\vspace*{0.2in} \hspace*{0.7in}%
\setlength{\unitlength}{0.00017in}\begin{picture}
(0,0) \thicklines
\put(2000,-500){\vector(3,-2){7000}}
\put(2900,-3700){\vector(3,-1){5800}}
\put(2000,-11900){\vector(3,2){7000}}
\put(800,-1000){\vector(0,-1){1500}}
\put(800,-4000){\vector(0,-1){4500}}
\put(2900,-8700){\vector(3,1){5800}}
\put(800,-11400){\vector(0,1){1500}}
\put(800,0){\makebox(0,0){$P_{n-1}$}}
\put(800,-12200){\makebox(0,0){$P_{n}$}}
\put(-2000,-1600){\makebox(0,0){$_{\Delta_{n_1^{\prime}+1,n_1^{\prime
\prime}+1,n_2+1}}$}}
\put(-950,-6200){\makebox(0,0){$_{h_{A_{1}|B_{1}}\times 1}$}}
\put(800,-3200){\makebox(0,0){$_{P_{n_1^{\prime}+1} \times
P_{n_1^{\prime \prime }+1}\times P_{n_2+1}}$}}
\put(800,-9200){\makebox(0,0){$_{ P_{n_1+1}\times P_{n_2+1}}$}}
\put(9600,-6200){\makebox(0,0){$Y$}}
\put(-1300,-10800){\makebox(0,0){$_{\Delta_{n_1+1,n_2+1}}$}}
\put(5600,-7200){\makebox(0,0){$\bar{f}$}}
\put(5600,-10300){\makebox(0,0){$f$}}
\put(6400,-1900){\makebox(0,0){$d_{A|B}(f)$}}
\put(5600,-5450){\makebox(0,0){$\tilde{f}$}}
\put(17200,-500){\vector(-3,-2){7000}}
\put(16300,-3700){\vector(-3,-1){5800}}
\put(17200,-11900){\vector(-3,2){7000}}
\put(18400,-1000){\vector(0,-1){1500}}
\put(18400,-4000){\vector(0,-1){4500}}
\put(16300,-8700){\vector(-3,1){5800}}
\put(18400,-11400){\vector(0,1){1500}}
\put(18400,0){\makebox(0,0){$P_{n}$}}
\put(18400,-12200){\makebox(0,0){$P_{n-1}$}}
\put(21600,-10800){\makebox(0,0){$_{\Delta_{n_1^{\prime}+1,n_1^{\prime\prime}+1,n_2+1}}$}}
\put(20200,-6200){\makebox(0,0){$_{\phi_{A_{1}|B_{1}} \times 1}$}}
\put(18400,-9200){\makebox(0,0){$_{P_{n_1^{\prime}+1} \times
P_{n_1^{\prime \prime }+1}\times P_{n_2+1}}$}}
\put(18400,-3200){\makebox(0,0){$_{ P_{n_1+1}\times P_{n_2+1}}$}}
\put(20700,-1600){\makebox(0,0){$_{\Delta_{n_1+1,n_2+1}}$}}
\put(13600,-5350){\makebox(0,0){$\tilde{g}$}}
\put(13600,-10300){\makebox(0,0){$g$}}
\put(12600,-1900){\makebox(0,0){$\varrho_{A|B}(g)$}}
\put(13600,-7000){\makebox(0,0){$\bar{g}$}}
\end{picture}\vspace*{2.2in}

\begin{center}
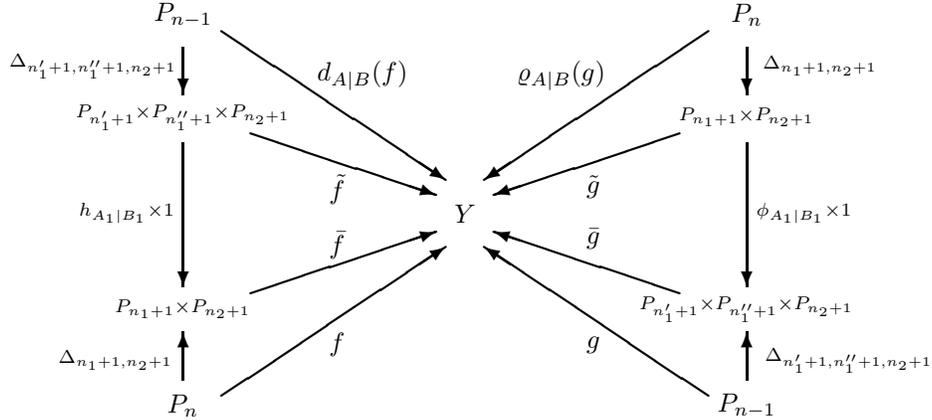
Figure 8:\ Face and degeneracy operators when $i=1$ and $k=2.$\vspace*{0.2in}
\end{center}

For $f=\bar{f}\circ\Delta_{n_{1}+1\cdots n_{k}+1}\in Sing_{n}^{M}Y,$ let
$\tilde{f}=\bar{f}\circ(1^{\times i-1}\times h_{A_{i}|B_{i}}\times1^{\times
k-i})$ and define
\[
d_{A|B}(f)=\tilde{f}\circ\Delta_{n_{1}+1\cdots n_{i}^{\prime}+1,n_{i}%
^{\prime\prime}+1\cdots n_{k}+1}.
\]
Dually, note that $n_{i}^{\prime}+n_{i}^{\prime\prime}=n_{i}-1$ implies the
sum of coordinates $(n_{1},\ldots,n_{i-1},n_{i}^{\prime},n_{i}^{\prime\prime
},$ $n_{i+1},\ldots,n_{k})\in\left(  \mathbb{N}\cup0\right)  ^{k+1}$ is $n-2.$
So for $g=\bar{g}\circ\Delta_{n_{1}+1\cdots n_{i}^{\prime}+1,n_{i}%
^{\prime\prime}+1\cdots n_{k}+1}\in Sing_{n-1}^{M}Y,$ let $\tilde{g}=\bar
{g}\circ(1^{\times i-1}\times\phi_{A_{i}|B_{i}}\times1^{\times k-i})$ and
define%
\[
\varrho_{A|B}(g)=\tilde{g}\circ\Delta_{n_{1}+1\cdots n_{k}+1}%
\]
(see Figure 8).

On the other hand, assume that $A|B\notin{\mathcal{Q}}_{p_{i},q_{i}}(n)$ for
all $i$ and define $d_{A|B}$ inductively as follows: When $k=2,$ set
$r=n_{1}+1,$ $s=n_{2}+1$ and let%
\[%
\begin{array}
[c]{ll}%
K|L= & \left\{
\begin{array}
[c]{ll}%
(\underline{r}\cap A)\cup\overline{s}\,|\,\underline{r}\cap B, & r\in A\\
\underline{r}\cap A\,|\,(\underline{r}\cap B)\cup\overline{s}, & r\in B
\end{array}
\right. \\
& \\
M|N= & \left\{
\begin{array}
[c]{lll}%
(\overline{s}\cap A)-1\,|\,\underline{n-1}\setminus(\overline{s}\cap A)-1, &
r\in B & \\
\underline{n-1}\setminus(\overline{s}\cap B)-\#L\,|\,(\overline{s}\cap
B)-\#L, & r\in A, & n\in A\\
I_{\underline{n}\setminus L}(A)\,|\,\underline{n-1}\setminus I_{\underline
{n}\setminus L}(A), & r\in A, & n\in B
\end{array}
\right. \\
& \\
C|D= & \left\{
\begin{array}
[c]{lll}%
I_{\underline{n}\setminus B}(\underline{r}\cap A)\,|\,\underline{n-1}\setminus
I_{\underline{n}\setminus B}(\underline{r}\cap A), & r\in B, & n\in B\\
I_{\underline{n}\setminus A}(\overline{s}\cap B)\,|\,\underline{n-1}\setminus
I_{\underline{n}\setminus A}(\overline{s}\cap B), & r\in A, & n\in B\\
&  & \\
\underline{n-1}\setminus I_{\underline{n}\setminus B}(\overline{s}\cap
A)\,|\,I_{\underline{n}\setminus B}(\overline{s}\cap A), & r\in B, & n\in A\\
\underline{n-1}\setminus I_{\underline{n}\setminus A}(\underline{r}\cap
B)\,|\,I_{\underline{n}\setminus A}(\underline{r}\cap B), & r\in A, & n\in A.
\end{array}
\right.
\end{array}
\]
Then define
\begin{equation}
d_{A|B}=\varrho_{C|D}d_{M|N}d_{K|L}. \label{multip}%
\end{equation}

\begin{remark}
This definition makes sense since $K|L\in{\mathcal{Q}}_{p_{1},q_{1}}%
(n),\,$\ $M|N\in{\mathcal{Q}}_{p_{3},q_{3}}(n-1),$ $C|D\in{\mathcal{Q}}%
_{p_{1},q_{1}}(n-1)\,$with either $r,n\in B$ or $r,n\in A$\thinspace and
$C|D\in{\mathcal{Q}}_{p_{3},q_{3}}(n-1)\,$with either $r\in B,$ $n\in A$ or
$r\in A,$ $n\in B$. Of course, ${\mathcal{Q}}_{\ast\ast}(n-1)$ is considered
with respect to the decomposition $n-2=m_{1}+m_{2}+m_{3}$ fixed after the
action of $d_{K|L}(\underline{r}\times\underline{s})$.
\end{remark}

If $k=3,$ consider the pair $(r,s)=(n_{1}+1,n-n_{1}),$ then $(r_{1}%
,s_{1})=(n_{2}+1,n-n_{1}-n_{2}-1)$ for $A_{1}|B_{1}=I_{\underline{n}%
\setminus\underline{r}}(\overline{s}\cap A)|I_{\underline{n}\setminus
\underline{r}}(\overline{s}\cap B)\in{\mathcal{P}}_{p_{1},q_{1}}(n-r),$ and so
on. Now dualize and use the same formulas above to define the degeneracy
operator $\varrho_{A|B}$.

\begin{definition}
Let $Y$ be a topological space. The \underline{singular multipermutahedral
set} \underline{of $Y$} consists of the singular set $Sing_{\ast}^{M}Y$
together with the singular face and degeneracy operators
\[
d_{A|B}:Sing_{n}^{M}Y\rightarrow Sing_{n-1}^{M}Y\text{ \ and \ }\varrho
_{A|B}:Sing_{n-1}^{M}Y\rightarrow Sing_{n}^{M}Y
\]
defined respectively for each $n\geq2$ and $A|B\in{\mathcal{P}}_{\ast\ast
}\left(  n\right)  $.
\end{definition}

\begin{remark}
The operator $d_{A|B}$ defined in (\ref{multip}) applied to $d_{U|V}$ for
some$\ U|V\in{\mathcal{P}}_{r,s}(n+1)$ yields the higher order structural
relation
\begin{equation}
d_{A|B}d_{U|V}=\varrho_{C|D}d_{M|N}d_{K|L}d_{U|V} \label{hrelation}%
\end{equation}
discussed in our first universal example.
\end{remark}

Now $Sing_{\ast}^{M}Y$ determines the singular (co)homology of a space $Y$ in
the following way: Let $R$ be a commutative ring with identity. For $n\geq1,$
let $C_{n-1}(Sing^{M}Y)$ denote the $R$-module generated by $Sing_{n}^{M}Y$
and form the \textquotedblleft chain complex\textquotedblright%
\[
(C_{\ast}(Sing^{M}Y),d)=\bigoplus_{\substack{n_{(k)}=n-1 \\n\geq1}%
}(C_{n-1}(Sing^{n_{1}\cdots n_{k}}Y),d_{n_{1}\cdots n_{k}}),
\]
where
\[
d_{n_{1}\cdots n_{k}}=\sum_{\substack{A|B\in\bigcup_{i=1}^{k}\mathcal{Q}%
_{p_{i},q_{i}}\left(  n\right)  }}\,-(-1)^{n_{(i-1)}+n_{i}^{\prime}}\text{
\textit{ shuff}}(C_{i}\cap A;C_{i}\cap B)\text{ }d_{A|B}.
\]
Refer to the example in Figure 7 and note that for $f\in C_{4}(Sing^{M}Y)$
with\linebreak\ $d_{13|2}d_{12|34}\left(  f\right)  \neq0$, the component
$d_{13|2}d_{12|34}\left(  f\right)  $ of $d^{2}\left(  f\right)  \in
C_{2}(Sing^{M}Y)$ is not cancelled and $d^{2}\neq0$. Hence $d$ is not a
differential. To remedy this, form the quotient
\[
C_{\ast}^{\diamondsuit}(Y)=C_{\ast}\left(  Sing^{M}Y\right)  /DGN,
\]
where $DGN$ is the submodule generated by the degeneracies, and obtain the
\emph{singular permutahedral chain complex }$\left(  C_{\ast}^{\diamondsuit
}(Y),d\right)  $. Because the signs in $d$ are determined by the index $i$,
which is missing in our first universal example, we are unable to use our
first example to define a chain complex with signs. However, we could use it
to define a unoriented theory with $\mathbb{Z}_{2}$-coefficients.

The singular homology of $Y$ is recovered from the composition
\[
C_{\ast}(SingY)\rightarrow C_{\ast}(Sing^{I}Y)\rightarrow C_{\ast}%
(Sing^{M}Y)\rightarrow C_{\ast}^{\diamondsuit}(Y)
\]
arising from the canonical cellular projections
\[
P_{n+1}\rightarrow I^{n}\rightarrow\Delta^{n}.
\]
Since this composition is a chain map, there is a natural isomorphism
\[
H_{\ast}(Y)\approx H_{\ast}^{\diamondsuit}(Y)=H_{\ast}(C_{\ast}^{\diamondsuit
}(Y),d).
\]
The fact that our diagonal on $P$ and the A-W diagonal on simplices commute
with projections allows us to recover the singular cohomology ring of $Y$ as
well. Finally, we remark that a cellular projection $f$ between polytopes
induces a chain map between corresponding singular chain complexes whenever
chains on the target are normalized. Here $C_{\ast}(SingY)$ and $C_{\ast
}(Sing^{I}Y)$ are non-normalized and the induced map $f^{\ast}$ is not a chain
map; but fortunately $d^{2}=0$ does not depend $df^{\ast}=f^{\ast}d$.

\subsection{Abstract Permutahedral Sets}

We begin by constructing a generating category $\mathbf{P}$ for permutahedral
sets similar to that of finite ordered sets and monotonic maps for simplicial
sets. The objects of $\mathbf{P}$ are the sets $n!=S_{n}$ of permutations of
$\underline{n},$ $n\geq1.$ But before we can define the morphisms we need some
preliminaries. First note that when $P_{n}$ is identified with its vertices
$n!$, the maps $\rho_{n}$ and $\gamma_{n}$ defined above become
\[
\rho_{n}:n!\rightarrow2!^{n-1}\ \ \text{and}\ \ \gamma_{n}:2!^{n-1}\rightarrow
n!.
\]
Given a non-empty increasingly ordered set $M=\left\{  m_{1}<\cdots
<m_{k}\right\}  \subset\mathbb{N},$ let $M!$ denote the set of all
permutations of $M$ and let $J_{M}:M!\rightarrow k!$ be the map defined for
$a=\left(  m_{\sigma\left(  1\right)  },...,m_{\sigma\left(  k\right)
}\right)  \in M!$ by $J_{M}(a)=\sigma.$ For $n,m\in\mathbb{N}$ and partitions
$A_{1}|\cdots|A_{k}\in{\mathcal{P}}_{n_{1}\cdots n_{k}}(n)$ and $B_{1}%
|\cdots|B_{\ell}\in{\mathcal{P}}_{m_{1}\cdots m_{\ell}}(m)$ with
$n-k=m-\ell=\varkappa,$ define the morphism
\[
f_{A_{1}|\cdots|A_{k}}^{B_{1}|\cdots|B_{\ell}}:m!\rightarrow n!
\]
by the composition
\[
m!\overset{sh_{B}}{\rightarrow}\prod_{j=1}^{\ell}B_{j}\overset{\sigma_{\max}%
}{\rightarrow}\prod_{r=1}^{\ell}B_{j_{r}}\overset{J_{B}}{\rightarrow}%
\prod_{j=r}^{\ell}m_{j_{r}}!\overset{\rho_{\ast}}{\rightarrow}2!^{\varkappa
}\overset{\gamma_{\ast}}{\rightarrow}\prod_{s=1}^{k}n_{i_{s}}!\overset
{J_{A}^{-1}}{\rightarrow}\prod_{s=1}^{k}A_{i_{s}}\overset{\sigma_{\max}^{-1}%
}{\rightarrow}\prod_{i=1}^{k}A_{i}\overset{{\iota_{A}}}{{\rightarrow}}n!
\]
where $sh_{B}$ is a surjection defined for $b=\{b_{1},...,b_{m}\}\in m!$ by
\[
sh_{B}(b)=(b_{1,1},..,b_{m_{1},1};...;b_{1,\ell},..,b_{m_{\ell},\ell}),\,
\]
in which the right-hand side is the unshuffle of $b$ with $b_{r,t}\in
B_{t},\,1\leq r\leq m_{t},\,1\leq t\leq\ell;$ $\sigma_{\max}\in S_{\ell}$ is a
permutation defined by $j_{r}=\sigma_{\max}(r),\,$ $\max B_{j_{r}}=\max
(B_{1}\cup B_{2}\cup\cdots\cup B_{j_{r}});$ $J_{B}=\prod_{r=1}^{\ell
}J_{B_{j_{r}}};$ $\rho_{\ast}=\prod_{r=1}^{\ell}\rho_{j_{r}}$ and
$\gamma_{\ast}=\prod_{s=1}^{k}\gamma_{i_{s}};$ finally, $\iota_{A}$ is the
inclusion. It is easy to see that
\[
f_{A_{1}|\cdots|A_{k}}^{B_{1}|\cdots|B_{\ell}}=f_{A_{1}|\cdots|A_{k}%
}^{\underline{\varkappa+1}}\circ f_{\underline{\varkappa+1}}^{B_{1}%
|\cdots|B_{\ell}}\ \ \text{and}\ \ f_{\underline{n}}^{\underline{n}}%
=\gamma_{n}\circ\rho_{n}.
\]
In particular, the maps $f_{A|B}^{\underline{n-1}}:(n-1)!\rightarrow n!$ and
$f_{\underline{n-1}}^{A|B}:n!\rightarrow\left(  n-1\right)  !$ are generator
morphisms denoted by $\delta_{A|B}$ and $\beta_{A|B},$ respectively (see
Theorem \ref{highrelations} below, the statement of which requires some new
set operations).

\begin{definition}
Given non-empty disjoint subsets $A,B,U\subset\underline{n+1}$ with $A\cup
B\subseteq U,$ define the \underline{lower and upper disjoint unions} (with
respect to $U$) by\vspace{0.1in}\newline%
\begin{tabular}
[c]{lll}%
$\hspace*{0.3in}$ & $A\underline{\sqcup}B=$ & $\left\{
\begin{array}
[c]{ll}%
I_{U\diagdown A}\left(  B\right)  +\#A-1, & \text{if }\min B>\min\left(
U\diagdown A\right) \\
I_{U\diagdown A}\left(  B\right)  +\#A-1\cup\underline{\#A}, & \text{if }\min
B=\min\left(  U\diagdown A\right)
\end{array}
\right.  $\\
and &  & \\
& $A\overline{\sqcup}B=$ & $\left\{
\begin{array}
[c]{ll}%
I_{U\diagdown B}\left(  A\right)  , & \text{if }\max A<\max\left(  U\diagdown
B\right) \\
I_{U\diagdown B}\left(  A\right)  \cup\overline{\#B}-1, & \text{if }\max
A=\max\left(  U\diagdown B\right)  \text{.}%
\end{array}
\right.  $%
\end{tabular}
\vspace{0.1in}\newline If either $A$ or $B$ is empty, define $A\underline
{\sqcup}B=A\overline{\sqcup}B=A\cup B$. Furthermore, given non-empty disjoint
subsets $A,B_{1},\ldots,B_{k}\subset\underline{n+1}$ with $k\geq1,$ set
$U=A\cup B_{1}\cup\cdots\cup B_{k}$ and define
\[
A\square(B_{1}|\cdots|B_{k})=(B_{1}|\cdots|B_{k})\square A=\left\{
\begin{array}
[c]{ll}%
A\underline{\sqcup}B_{1}|\cdots|A\underline{\sqcup}B_{k}, & \text{if }\max
A<\max U\\
B_{1}\overline{\sqcup}A|\cdots|B_{k}\overline{\sqcup}A, & \text{if }\max
A=\max U.
\end{array}
\right.
\]

\end{definition}

\noindent Note that if $A|B$ is a partition of $\underline{n+1}$, then
\[
A\underline{\sqcup}B=A\overline{\sqcup}B=\underline{n}.
\]
Given a partition $A_{1}|\cdots|A_{k+1}$ of $\underline{n},$ define $A_{1}%
^{1}|\cdots|A_{k+1}^{1}=A_{1}^{1}|\cdots|A_{1}^{k+1}=$\linebreak$A_{1}%
|\cdots|A_{k+1};$ inductively, given $A_{1}^{i}|\cdots|A_{k-i+2}^{i}$ the
partition of $\underline{n-i+1},\,1\leq i<k,$ let
\[
A_{1}^{i+1}|\cdots|A_{k-i+1}^{i+1}=A_{1}^{i}\square(A_{2}^{i}|\cdots
|A_{k-i+2}^{i})
\]
be the partition of $\underline{n-i};$ and given $A_{i}^{1}|\cdots
|A_{i}^{k-i+2}$ the partition of $\underline{n-i+1},\,1\leq i<k,$ let
\[
A_{i+1}^{1}|\cdots|A_{i+1}^{k-i+1}=(A_{i}^{1}|\cdots|A_{i}^{k-i+1})\square
A_{i}^{k-i+2}%
\]
be the partition of $\underline{n-i}.$

\begin{theorem}
\label{highrelations}For $A_{1}|\cdots|A_{k+1}\in{\mathcal{P}}_{n_{1}\cdots
n_{k+1}}(n),$ $2\leq k\leq n,$ the map $f_{A_{1}|\cdots|A_{k+1}}%
^{\underline{n-k}}:(n-k)!\rightarrow n!\,$ can be expressed as a composition
of $\delta$'s two ways:
\[
f_{A_{1}|\cdots|A_{k+1}}^{\underline{n-k}}=\delta_{A_{1}^{1}\,|\,A_{2}^{1}%
\cup\cdots\cup A_{k+1}^{1}}\cdots\delta_{A_{1}^{k}\,|\,A_{2}^{k}}%
=\delta_{A_{1}^{1}\cup\cdots\cup A_{1}^{k}\,|\,A_{1}^{k+1}}\cdots\delta
_{A_{k}^{1}\,|\,A_{k}^{2}}.
\]

\end{theorem}

\begin{proof}
The proof is straightforward and omitted.
\end{proof}

\noindent There is also the dual set of relations among the $\beta$'s.

\begin{example}
Theorem \ref{highrelations} defines structure relations among the $\delta$'s,
the first of which is
\begin{equation}
\delta_{A|B\cup C}\,\delta_{A\square(B|C)}=\delta_{A\cup B|C}\,\delta
_{(A|B)\square C} \label{coquadrel}%
\end{equation}
when $k=2.$ In particular, let $A|B|C=12|345|678$. Since $A\underline{\sqcup
}B=\{1234\},$ $A\underline{\sqcup}C=\left\{  567\right\}  ,$ $A\overline
{\sqcup}C=\{12\}$ and $B\overline{\sqcup}C=\left\{  34567\right\}  ,$ we
obtain the following quadratic relation on $12|345|678$:
\[
\delta_{12|345678}\delta_{1234|567}=\delta_{12345|678}\delta_{12|34567};
\]
similarly, on $345|12|678$ we have
\[
\delta_{345|12678}\delta_{1234|567}=\delta_{12345|678}\delta_{34567|12}.
\]

\end{example}

\begin{theorem}
\label{faces}Let $A|B\in\mathcal{P}_{p,q}\left(  n+1\right)  $ and
$C|D\in\mathcal{P}_{\ast\ast}\left(  n\right)  .\ $Then $\delta_{A|B}%
\delta_{C|D}$ coincides with a map $f_{X|Y|Z}^{\underline{n-1}}%
:(n-1)!\rightarrow(n+1)!$ if and only if%
\begin{equation}
C|D\in\left\{
\begin{array}
[c]{ll}%
\mathcal{Q}_{q,1}(n)\cup\mathcal{Q}_{1,p}\left(  n\right)  , & \text{if
}n+1\in A\\
& \\
\mathcal{Q}_{p,1}(n)\cup\mathcal{Q}_{1,q}(n), & \text{if }n+1\in B.
\end{array}
\right.  \label{conditions}%
\end{equation}

\end{theorem}

\begin{proof}
If $\delta_{A|B}\delta_{C|D}$ coincides with $f_{X|Y|Z}^{\underline{n-1}},$
then according to relation (\ref{quadrel}) we have either
\[
{A|B}=X{|}Y{\cup Z}\ \text{and}\ {C|D}=X\square(Y|Z)
\]
or
\[
{A|B}={X\cup Y|Z}\ \text{and}\ {C|D}=(X|Y)\square Z.
\]
Hence there are two cases.\vspace*{0.1in}\newline\underline{\textit{Case 1:}%
}\ $A|B=X|Y\cup Z.$\vspace*{0.1in}\newline\underline{\textit{Subcase 1a:}%
}\ Assume $n+1\in A.$ If $\max Y=\max(Y\cup Z),$ then $\overline{p}\subseteq
Y\overline{\sqcup}X;$ otherwise $\max\left(  Y\cup Z\right)  =\max Z$ and
$\overline{p}\subseteq Z\overline{\sqcup}X.\ \ $In either case,
$C|D=Y\overline{\sqcup}X|Z\overline{\sqcup}X\in\mathcal{Q}_{1,p}\left(
n\right)  $.\vspace*{0.1in}\newline\underline{\textit{Subcase 1b:}}\ Assume
$n+1\in B.$ If $\min Y=\min(Y\cup Z),$ then $\underline{p}\subseteq
X\underline{\sqcup}Y;$ otherwise $\min(Y\cup Z)=\min Z$ and $\underline
{p}\subseteq X\underline{\sqcup}Z.\ $In either case, {$C$}${|}${$D$}
${=}X\underline{\sqcup}Y|X\underline{\sqcup}Z\in\mathcal{Q}_{p,1}\left(
n\right)  .$ \vspace*{0.1in}\newline\underline{\textit{Case 2:}}\ $A|B=X\cup
Y|Z.$\vspace*{0.1in}\newline\underline{\textit{Subcase 2a:}}\ Assume $n+1\in
A.$ If $\min X=\min(X\cup Y),$ then $\underline{q}\subseteq Z\underline
{\sqcup}X;$ otherwise $\min(X\cup Y)=\min Y$ and $\underline{q}\subseteq
Z\underline{\sqcup}Y.\ $In either case, {$C$}${|}${$D$} ${=}Z\underline
{\sqcup}X|Z\underline{\sqcup}Y\in\mathcal{Q}_{q,1}\left(  n\right)  .$%
\vspace*{0.1in}\newline\underline{\textit{Subcase 2b:}}\ Assume $n+1\in B.$ If
$\max X=\max(X\cup Y),$ then $\overline{q}\subseteq X\overline{\sqcup}Z;$
otherwise $\max\left(  X\cup Y\right)  =\max Y$ and $\overline{q}\subseteq
Y\overline{\sqcup}Z.\ \ $In either case, $C|D=X\overline{\sqcup}%
Z|Y\overline{\sqcup}Z\in\mathcal{Q}_{1,q}\left(  n\right)  $.\vspace
*{0.1in}\newline Conversely, given $A|B\in\mathcal{P}_{p,q}(n+1)$ and $C|D$
satisfying conditions (\ref{conditions}) above, let
\[
\lbrack A|B;\,C|D]=\left\{
\begin{array}
[c]{ll}%
A\,|\,I_{B}^{-1}\left(  \overline{q}\cap C-p+1\right)  \,|\,I_{B}^{-1}\left(
\overline{q}\cap D-p+1\right)  , & C|D\in\mathcal{Q}_{p,1}\left(  n\right) \\
& \\
I_{A}^{-1}\left(  \underline{p}\cap C\right)  \,|\,I_{A}^{-1}\left(
\underline{p}\cap D\right)  \,|\,B, & C|D\in\mathcal{Q}_{1,q}\left(  n\right)
\end{array}
\right.
\]
and
\[
\lbrack A|B;\,C|D]=\left\{
\begin{array}
[c]{ll}%
A\,|\,I_{B}^{-1}\left(  \underline{q}\cap C\right)  \,|\,I_{B}^{-1}\left(
\underline{q}\cap D\right)  , & C|D\in\mathcal{Q}_{1,p}\left(  n\right) \\
& \\
I_{A}^{-1}\left(  \overline{p}\cap C-q+1\right)  \,|\,I_{A}^{-1}\left(
\overline{p}\cap D-q+1\right)  \,|\,B, & C|D\in\mathcal{Q}_{q,1}\left(
n\right)  .
\end{array}
\right.
\]
A straightforward calculation shows that
\[
\left[  X|Y\cup Z;\text{ }X\square(Y|Z)\right]  =X|Y|Z=[{X\cup Y|Z};\text{
}(X|Y)\square Z].
\]
Consequently, if $X|Y|Z=[A|B;\,C|D],$ either%
\[
{A|B}={X|}\text{ }{Y\cup Z}\ \text{and}\ {C|D}=X\square(Y|Z)
\]
when $C|D\in\mathcal{Q}_{p,1}(n)\cup\mathcal{Q}_{1,p}(n)$ or
\[
{A|B}={X\cup Y|Z}\ \text{and}\ {C|D}=(X|Y)\square Z
\]
when $C|D\in\mathcal{Q}_{q,1}(n)\cup\mathcal{Q}_{1,q}(n).$
\end{proof}

\noindent On the other hand, if $C|D\not \in \mathcal{Q}_{p,1}(n)\cup
\mathcal{Q}_{1,p}(n)\cup\mathcal{Q}_{q,1}(n)\cup\mathcal{Q}_{1,q}(n),$ higher
order structure relations involving both coface and codegeneracy operators appear.

\begin{definition}
Let $\mathcal{C}$ be the category of sets. A \underline{permutahedral set} is
a contravariant functor
\[
\mathcal{Z}:\mathbf{{P}\rightarrow\mathcal{C}}.
\]

\end{definition}

Thus a permutahedral set $\mathcal{Z}$ is a graded set $\mathcal{Z}%
=\{\mathcal{Z}_{n}\}_{n\geq1}\,$endowed with face and degeneracy operators
\[
d_{A|B}=\mathcal{Z}(\delta_{A|B}):\mathcal{Z}_{n}\rightarrow\mathcal{Z}%
_{n-1}\text{ \ and \ }\varrho_{M|N}=\mathcal{Z}(\beta_{M|N}):\mathcal{Z}%
_{n}\rightarrow\mathcal{Z}_{n+1}%
\]
satisfying an appropriate set of relations, which includes quadratic relations
such as
\begin{equation}
d_{A\square(B|C)}d_{A|B\cup C}\ =d_{(A|B)\square C}d_{A\cup B|C}
\label{quadrel}%
\end{equation}
induced by (\ref{coquadrel}) and higher order relations such as
\[
d_{A|B}d_{U|V}=\varrho_{C|D}d_{M|N}d_{K|L}d_{U|V}%
\]
discussed in (\ref{hrelation}).

Let us define the abstract analog of a singular multipermutahedral set, which
leads to a singular chain complex with arbitrary coefficients.

\begin{definition}
\label{derformal}For $n\geq1,\mathbb{\ }$let $X_{n}=\bigcup_{n_{\left(
k\right)  }=n-1,\text{ }n_{k}\geq0\,}X^{n_{1}\cdots n_{k}}$ and $X_{n-1}=$
\linebreak$\bigcup_{m_{\left(  \ell\right)  }=n-2,\text{ }m_{\ell}\geq
0\,}X^{m_{1}\cdots m_{\ell}}$ be filtered sets; let $A|B\in{\mathcal{Q}%
}_{p_{i},q_{i}}(n)$ for some $i.$ A map $g:X_{n}\rightarrow X_{n-1}$ acts as
an \underline{$A|B$-formal derivation} if $g|_{X^{n_{1}\cdots n_{k}}}%
:X^{n_{1}\cdots n_{k}}\rightarrow X^{n_{1}\cdots n_{i}^{\prime},n_{i}%
^{\prime\prime}\cdots n_{k}},$ where $(n_{i}^{\prime},n_{i}^{\prime\prime})$
is given by (\ref{formal}).
\end{definition}

Let $\mathcal{C}_{M}$ denote the category whose objects are positively graded
sets $X_{\ast}$ filtered by subsets $X_{n}=\bigcup_{n_{\left(  k\right)
}=n-1,\text{ }n_{k}\geq0}X^{n_{1}\cdots n_{k}}$ and whose morphisms are
filtration preserving set maps.

\begin{definition}
A \underline{multipermutahedral set} is a contravariant functor $\mathcal{Z}%
:\mathbf{{P}\rightarrow\mathcal{C}}_{M}$ such that
\[
\mathcal{Z}(\delta_{A|B}):\mathcal{Z}(n!)\rightarrow\mathcal{Z}((n-1)!)
\]
acts as an $A|B$-formal derivation for each $A|B\in{\mathcal{Q}}_{p_{i},q_{i}%
},$ all $i\geq1.$
\end{definition}

Thus a multipermutahedral set $\mathcal{Z}$ is a graded set $\left\{
\mathcal{Z}_{n}\right\}  _{n\geq1}$ with
\[
\mathcal{Z}_{n}=\bigcup_{\substack{n_{\left(  k\right)  }=n-1 \\n_{k}\geq
0}}\mathcal{Z}^{n_{1}\cdots n_{k}},
\]
$\,$together with face and degeneracy operators
\[
d_{A|B}=\mathcal{Z}(\delta_{A|B}):\mathcal{Z}_{n}\rightarrow\mathcal{Z}%
_{n-1}\text{ \ and \ }\varrho_{M|N}=\mathcal{Z}(\beta_{M|N}):\mathcal{Z}%
_{n}\rightarrow\mathcal{Z}_{n+1}%
\]
satisfying the relations of a permutahedral set and the additional requirement
that $d_{A|B}$ respect underlying multigrading. This later condition allows us
to form the chain complex of $\mathcal{Z}$ with signs mimicking the cellular
chain complex of permutahedra (see below). Note that the chain complex of a
permutahedral set is only defined with $\mathbb{Z}_{2} $-coefficients in general.

\subsection{The Cartesian product of permutahedral sets}

The objects and morphisms in the category $\mathbf{P\times P}$ are the sets
and maps
\[
n!!=\bigcup_{r+s=n}r!\times s!\text{ \ and \ }\bigcup_{f,g\in\mathbf{P}%
}f\times g:m!!\rightarrow n!!
\]
all $m,n\geq1.$ There is a functor $\Delta:\mathbf{P\rightarrow P\times P}$
defined as follows. If ${A|B}\in\mathcal{Q}_{r,1}(n)\cup\mathcal{Q}_{1,s}(n),$
define $\Delta_{r,s}(A|B)=A_{1}|B_{1}\times A_{2}|B_{2}\in r!\times s!\ $and
define $\delta_{A|B}:(n-1)!\rightarrow n!$ by
\[
\Delta(\delta_{A|B})=\delta_{A_{1}|B_{1}}\times\delta_{A_{2}|B_{2}},
\]
where $\delta_{A_{i}|B_{i}}=1$ for either $i=1$ or $i=2.$ Define $\Delta
(\beta_{A|B})$ similarly$.$ On the other hand, if $A|B\notin\mathcal{Q}%
_{r,1}(n)\cup\mathcal{Q}_{1,s}(n),$ define
\[
\Delta(\delta_{A|B})=\Delta(\delta_{K|L})\Delta(\delta_{M|N})\Delta
(\beta_{C|D}),
\]
where $K|L,\,M|N,\,C|D$ are given by the formulas in (\ref{multip}). Dually,
define $\Delta(\beta_{M|N}).$ It is easy to check that $\Delta$ is well defined.

Given multipermutahedral sets $\mathcal{Z}^{\prime},\mathcal{Z}^{\prime\prime
}:\mathbf{P}\rightarrow\mathcal{C}_{M}$, first define a functor
\[
\mathcal{Z}^{\prime}\tilde{\times}\mathcal{Z}^{\prime\prime}:\mathbf{P}%
\times\mathbf{P}\rightarrow\mathcal{C}_{M}%
\]
on an object $n!!$ by
\[
(\mathcal{Z}^{\prime}\tilde{\times}\mathcal{Z}^{\prime\prime})(n!!)=\bigcup
_{r+s=n}\mathcal{Z}^{\prime}(r!)\times\mathcal{Z}^{\prime\prime}%
(s!)\diagup{\Huge \sim},
\]
where $\left(  a,b\right)  \sim\left(  c,e\right)  $ if and only if
$a=\varrho_{\underline{r}|r+1}^{\prime}(c)$ and $e=\varrho_{1|\underline
{s+1}\setminus\underline{1}}^{\prime\prime}(b).$ On a map $h=\bigcup(f\times
g):m!!\rightarrow n!!,$
\[
(\mathcal{Z}^{\prime}\tilde{\times}\mathcal{Z}^{\prime\prime})(h):(\mathcal{Z}%
^{\prime}\tilde{\times}\mathcal{Z}^{\prime\prime})(n!!)\rightarrow
(\mathcal{Z}^{\prime}\tilde{\times}\mathcal{Z}^{\prime\prime})(m!!)
\]
is the map induced by $\bigcup(\mathcal{Z}^{\prime}(f)\times\mathcal{Z}%
^{\prime\prime}(g)).$ Now define the product $\mathcal{Z}^{\prime}%
\times\mathcal{Z}^{\prime\prime}$ to be the composition of functors
\[
\mathcal{Z}^{\prime}\times\mathcal{Z}^{\prime\prime}=\mathcal{Z}^{\prime
}\tilde{\times}\mathcal{Z}^{\prime\prime}\circ\Delta:\mathbf{P\rightarrow
\mathcal{C}}_{M}.
\]
The face operator $d_{A|B}$ on $\mathcal{Z}^{\prime}\times\mathcal{Z}%
^{\prime\prime}$ is given by
\begin{equation}
d_{A|B}(a\times b)=\left\{
\begin{array}
[c]{ll}%
d_{\underline{r}\cap A|\underline{r}\cap B}^{\prime}\left(  a\right)  \times
b, & \text{\textit{if }}A|B\in\mathcal{Q}_{1,s}\left(  n\right)  ,\\
a\times d_{(\overline{s}\cap A)-r+1\,|\,(\overline{s}\cap B)-r+1}%
^{\prime\prime}(b), & \text{\textit{if }}A|B\in\mathcal{Q}_{r,1}\left(
n\right)  ,\\
\varrho_{C|D}d_{M|N}d_{K|L}\left(  a\times b\right)  , &
\text{\textit{otherwise}},
\end{array}
\right.  \label{productd1}%
\end{equation}
with $M|N,\,K|L,\,C|D$ given by the formulas in (\ref{multip}).

\begin{example}
The canonical map $\iota:Sing^{P}X\times Sing^{P}Y\rightarrow Sing^{P}(X\times
Y)$ defined for $(f,g)\in Sing_{r}^{P}X\times Sing_{s}^{P}Y$ by
\[
\iota(f,g)=(f\times g)\circ\Delta_{r,s}%
\]
is a map of permutahedral sets. Consequently, if $X$ is a topological monoid,
the singular permutahedral complex $Sing^{P}X$ inherits a canonical monoidal structure.
\end{example}

\subsection{The diagonal on a permutahedral set}

Let $\mathcal{Z}=(\mathcal{Z}_{n},d_{A|B},\varrho_{M|N})$ be a
multipermutahedral set. The chain complex of $\mathcal{Z}$ is
\[
C_{\ast}^{\diamondsuit}({\mathcal{Z}})=C_{\ast}({\mathcal{Z}})/DGN,
\]
where $DGN$ is the submodule generated by the degeneracies;
\[
(C_{\ast}({\mathcal{Z}}),d)=\bigoplus_{\substack{n_{(k)}+1=n\\n\geq1}%
}(C_{n-1}({\mathcal{Z}}^{n_{1}\cdots n_{k}}),d_{n_{1}\cdots n_{k}})
\]
and%
\[
d_{n_{1}\cdots n_{k}}=\,\sum\limits_{_{\substack{A|B\in%
{\textstyle\bigcup_{i=1}^{k}}
\mathcal{Q}_{p_{i},q_{i}}\left(  n\right)  }}}-(-1)^{n_{(i-1)}+n_{i}^{\prime}%
}\text{ \textit{shuff}}(C_{i}\cap A;C_{i}\cap B)\text{ }d_{A|B}.
\]
The explicit diagonal
\[
\Delta:C_{\ast}^{\diamondsuit}(\mathcal{Z})\rightarrow C_{\ast}^{\diamondsuit
}(\mathcal{Z})\otimes C_{\ast}^{\diamondsuit}(\mathcal{Z})
\]
on $a\subset\mathcal{Z}_{n}$ is given by
\begin{equation}
\Delta(a)=\sum\limits_{\substack{F\in\mathcal{C}^{q\times n-q+2}\\1\leq q\leq
n+1}}csgn(F)\,\ d_{c\left(  F\right)  }(a)\otimes d_{r\left(  F\right)  }(a),
\label{pformula}%
\end{equation}
where $d_{A_{1}|...|A_{k+1}}={\mathcal{Z}}(f_{A_{1}|...|A_{k+1}}%
^{\underline{n-k}}).$

\subsection{The double cobar-construction $\Omega^{2}C_{\ast}(X)$\label{sub4}}

Given a simplicial, cubical or a permutahedral set $W$ with base point $\ast,$
let ${C}_{\ast}(W)$ denote the quotient\linebreak${C}_{\ast}(W)/C_{>0}(\ast).$
Say that $W$ is $k$\textit{-reduced} if $W_{i}$ contains exactly one element
for each $i\leq k$ and let $\Omega C$ denote the cobar construction on a
1-reduced DG coalgebra $C.$ In \cite{KS1} and \cite{KS2}, Kadeishvili and
Saneblidze construct functors from the category of 1-reduced simplicial sets
to the category of cubical sets and from the category of 1-reduced cubical
sets to the category of multipermutahedral sets (denoted by $\mathbf{\Omega}$
in either case) for which the following statements hold (c.f. \cite{CM},
\cite{Baues1}):

\begin{theorem}
\cite{KS1} Given a 1-reduced simplicial set $X,$ there is a canonical
identification isomorphism
\[
\Omega{C}_{\ast}(X)\approx C_{\ast}^{\Box}(\mathbf{\Omega}X).
\]

\end{theorem}

\begin{theorem}
\cite{KS2}\label{KS2} Given a 1-reduced cubical set $Q,$ there is a canonical
identification isomorphism
\[
\Omega{C}_{\ast}^{\Box}(Q)\approx C_{\ast}^{\diamondsuit}(\mathbf{\Omega}Q).
\]

\end{theorem}

\noindent For completeness, definitions of these two functors appear in the
appendix. Since the chain complex of any cubical set $Q$ is a DG coalgebra
with strictly coassociative coproduct, setting $Q=\mathbf{\Omega}X$ in Theorem
\ref{KS2} immediately gives:

\begin{theorem}
For a 2-reduced simplicial set $X$ there is a canonical identification
isomorphism
\[
\Omega^{2}{C}_{\ast}(X)\approx C_{\ast}^{\diamondsuit}(\mathbf{\Omega^{2}}X).
\]

\end{theorem}

Now if $X=Sing^{1}Y,$ then $\Omega{C}_{\ast}(X)$ is Adams' cobar construction
for the space $Y$ \cite{Adams}; consequently, there is a canonical (geometric)
coproduct on $\Omega^{2}C_{\ast}(Sing^{1}Y).$ We shall extend this canonical
coproduct to an \textquotedblleft$A_{\infty}$-Hopf algebra\textquotedblright%
\ structure in the sequel \cite{SU3}.

\section{The Multiplihedra and Associahedra}

The multiplihedron $J_{n}$ and the associahedron $K_{n+1}$ are cellular
projections of $P_{n}$ defined in terms of planar trees. Consequently, we
shall need to index the faces of $P_{n}$ four ways: (1) by partitions of
$\underline{n}$ (see Section 2), (2) by (planar rooted) $p$-leveled trees with
$n+1$ leaves (PLT's), (3) by parenthesized strings of $n+1$ indeterminants
with $p-1$ levels of subscripted parentheses and (4) by $\left(  p-1\right)
$-fold compositions of face operators acting on $n+1$ indeterminants. The
second and third serve as transitional intermediaries between the first and fourth.

Define a correspondence between PLT's and partitions of $\underline{n}$ as
follows: Let $T_{n+1}^{p}$ be a PLT with $n+1$ leaves, $p$-levels and root in
level $p.$ Number the leaves from left to right and assign the label $i\ $to
the node at which the branch of leaf $i$ meets the branch of leaf $i+1,$
$1\leq i\leq n$ (a node may have multiple labels). Let $U_{j}=\left\{
\text{labels assigned to }j\text{-level nodes}\right\}  $ and identify
$T_{n+1}^{p}$ with the partition $U_{1}|\cdots|U_{p}$ of $\underline{n}$ (see
Figure 9).$\ $Thus binary $n$-leveled trees parametrize the vertices of
$P_{n}$. Loday and Ronco constructed a map from $S_{n}$ to binary $n$-leveled
trees \cite{Loday}; its extension to faces of $P_{n}$ was given by Tonks
\cite{tonks}. Note that the map from PLT's to partitions defined above gives
an inverse.\vspace{0.5in}

\begin{center}
\setlength{\unitlength}{0.005in}\begin{picture}
(152,-130)(202,-90) \thicklines
\put(150,-150){\line( -1,1){150}} \put(150,-150){\line( 1,1){150}}
\put(150,0){\line( 0,-1){250}} \put(100,-100){\line( 0,1){100}}
\put (100,-50){\line( -1,1){50}} \put(250,-50){\line( 0,1){50}}
\put(250,-50){\line ( -1,1){50}}
\put(150,-150){\makebox(0,0){$\bullet$}} \put(100,-100){\makebox
(0,0){$\bullet$}} \put(100,-50){\makebox(0,0){$\bullet$}} \put
(250,-50){\makebox(0,0){$\bullet$}}
\put(178,-160){\makebox(0,0){$34$}}
\put(80,-110){\makebox(0,0){$1$}}
\put(125,-55){\makebox(0,0){$2$}}
\put(285,-55){\makebox(0,0){$56$}}
\put(475,-50){\makebox(0,0){$\longleftrightarrow$\hspace{0.3in}$256|1|34$}}
\put(475,-125){\makebox(0,0){$\longleftrightarrow$\hspace{0.3in}
$d_{(0,1)}d_{(1,1)(4,2)}$}}
\put(475,-200){\makebox(0,0){$\longleftrightarrow$\hspace{0.3in}
$\left( \left( \bullet\left( \bullet
\text{\thinspace}\bullet\right) _1 \right)_2 \bullet
\left( \bullet \bullet \bullet\right)_1 \right)
$}}
\end{picture}\vspace{0.9in}

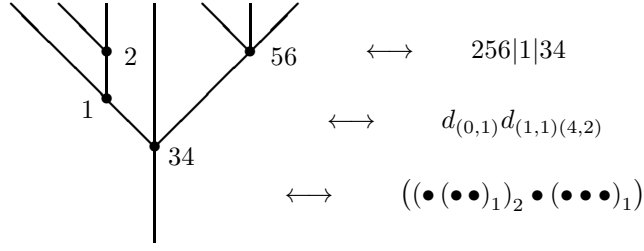
Figure 9: Various representations of the face $256|1|34.$\vspace{0.1in}
\end{center}

To define the correspondence between PLT's and subscripted parenthesizations
of $n+1$ indeterminants, begin by identifying the top cell of $P_{n}$ with the
$\left(  n+1\right)  $-leaf corolla and the (unsubscripted) parenthesized
string $\left(  x_{1}x_{2}\cdots x_{n+1}\right)  .$ Let $T_{n+1}^{p}$ be a PLT
with $p>1.$ If the branches meeting at a level 1 node contain leaves
$i,\ldots,i+k,$ enclose the corresponding indeterminants in a pair of
parentheses with subscript 1; if the branches meeting at a level 2 node
contain leaves $i_{1},\ldots,i_{k},$ enclose the corresponding indeterminants
in a pair of parentheses with subscript 2; and so on for $p-1$ steps (see
Example \ref{one}).

Compositions of face operators encode this parenthesization procedure. For
$s\geq1,$ choose $s$ pairs of indices $\left(  i_{1},\ell_{1}\right)
\cdots\left(  i_{s},\ell_{s}\right)  $ such that $0\leq i_{r}<i_{r+1}\leq n-1$
and $i_{r}+\ell_{r}+1\leq i_{r+1}.$ The face operator
\[
d_{\left(  i_{1},\ell_{1}\right)  \cdots\left(  i_{s},\ell_{s}\right)  }%
:P_{n}\rightarrow\partial P_{n}%
\]
acts on $\left(  x_{1}x_{2}\cdots x_{n+1}\right)  $ by simultaneously
inserting $s$ disjoint (non-nested) pairs of inner parentheses with subscript
$1$, the first enclosing $x_{i_{1}+1}\cdots x_{i_{1}+\ell_{1}+1},$ the second
enclosing $x_{i_{2}+1}\cdots x_{i_{2}+\ell_{2}+1},$ and so on.$\ $Thus,%
\[
d_{\left(  i_{1},\ell_{1}\right)  \cdots\left(  i_{s},\ell_{s}\right)
}\left(  x_{1}x_{2}\cdots x_{n+1}\right)  =\hspace*{2.5in}%
\]%
\[
\hspace*{1in}\left(  x_{1}\cdots\left(  x_{i_{1}+1}\cdots x_{i_{1}+\ell_{1}%
+1}\right)  _{1}\cdots\left(  x_{i_{s}+1}\cdots x_{i_{s}+\ell_{s}+1}\right)
_{1}\cdots x_{n+1}\right)  .
\]
A composition of face operators%
\begin{equation}
d_{\left(  i_{1}^{k},\ell_{1}^{k}\right)  \cdots\left(  i_{s_{k}}^{k}%
,\ell_{s_{k}}^{k}\right)  }\cdots d_{\left(  i_{1}^{1},\ell_{1}^{1}\right)
\cdots\left(  i_{s_{1}}^{1},\ell_{s_{1}}^{1}\right)  }:P_{n}\rightarrow
\partial^{k}P_{n} \label{comp}%
\end{equation}
continues this process inductively:$\ $If the $j^{th}$ operator inserted
parentheses with subscript $j$, treat each such pair and its contents as a
single indeterminant and apply the $\left(  j+1\right)  ^{st}$ as above,
inserting parentheses subscripted by $j+1.$

\label{one}Refer to Figure 9 above.$\ $The composition $d_{(0,1)}%
d_{(1,1)(4,2)}$ acts on $\left(  \bullet\bullet\bullet\bullet\bullet
\bullet\bullet\right)  $ in the following way: First, $d_{(1,1)(4,2)}$
simultaneously inserts two inner pairs of parentheses with subscript $1$:
\[
\left(  \bullet\left(  \bullet\text{\thinspace}\bullet\right)  _{1}%
\bullet\left(  \bullet\bullet\bullet\right)  _{1}\right)  .
\]
Next, $d_{(0,1)}$ inserts the single pair with subscript $2$:
\[
\left(  \left(  \bullet\left(  \bullet\text{\thinspace}\bullet\right)
_{1}\right)  _{2}\bullet\left(  \bullet\bullet\bullet\right)  _{1}\right)  .
\]
We summarize the discussion above as a proposition.

\begin{proposition}
\label{corr}The following correspondences (defined above) preserve
combinatorial structure:
\end{proposition}

\begin{tabular}
[c]{lll}%
$\left\{  \text{\textit{Faces of} }P_{n}\right\}  $ & $\leftrightarrow$ &
$\left\{  \text{\textit{Partitions of }}\underline{n}\right\}  $\\
&  & \\
& $\leftrightarrow$ & $\left\{  \text{\textit{Leveled trees with }}n+1\text{
leaves}\right\}  $\\
&  & \\
& $\leftrightarrow$ & $\left\{
\begin{tabular}
[c]{c}%
\textit{Strings of} $n+1\text{\ \textit{indeterminants}}$\\
\textit{with subscripted parentheses}%
\end{tabular}
\ \ \right\}  $\\
&  & \\
& $\leftrightarrow$ & $\left\{
\begin{array}
[c]{c}%
\text{\textit{Compositions of face operators}}\\
\text{\textit{acting on }}n+1\text{\ \textit{indeterminants.}}%
\end{array}
\right\}  $%
\end{tabular}
$\vspace{0.1in}$

Assign the identity face operator $Id$ to the top dimensional face of
$P_{n+1}$ and use the correspondences above to assign compositions of faces
operators to lower dimensional faces (see Figure 10). For faces in codimension
$1$ we have:\
\[%
\begin{tabular}
[c]{c|c}%
$\underset{\ }{\text{\textbf{Face of }}P_{n+1}}$ & \textbf{Face operator}%
\\\hline
$P_{n}\times0$ & $d_{\left(  0,n\right)  }$\\
& \\
$P_{n}\times1$ & $d_{\left(  n,1\right)  }$\\
& \\
$d_{\left(  i_{1},\ell_{1}\right)  \cdots\left(  i_{k},\ell_{k}\right)
}\times\left[  0,1-2^{\ell_{1}+\cdots+\ell_{k}-n}\right]  $ & $d_{\left(
i_{1},\ell_{1}\right)  \cdots\left(  i_{k},\ell_{k}\right)  }$\\
& \\
$d_{\left(  i_{1},\ell_{1}\right)  \cdots\left(  i_{k},\ell_{k}\right)
}\times\left[  1-2^{\ell_{1}+\cdots+\ell_{k}-n},1\right]  $ & $\left\{
\begin{array}
[c]{lc}%
d_{\left(  i_{1},\ell_{1}\right)  \cdots\left(  i_{k},\ell_{k}\right)  \left(
n,1\right)  }, & i_{k}+\ell_{k}<n\\
d_{\left(  i_{1},\ell_{1}\right)  \cdots\left(  i_{k},\ell_{k}+1\right)  }, &
i_{k}+\ell_{k}=n.
\end{array}
\right.  $%
\end{tabular}
\]
\vspace{0.1in}

\begin{center}
\setlength{\unitlength}{0.0005in}\begin{picture}
(2975,2685)(3126,-2038) \thicklines\put(3601,239){\line( 1,
0){1800}} \put(5401,239){\line( 0,-1){1800}}
\put(5401,-1561){\line(-1, 0){1800}} \put(3601,-1561){\line( 0,
1){1800}} \put(3601,239){\makebox(0,0){$\bullet$}}
\put(3601,-661){\makebox(0,0){$\bullet$}}
\put(3601,-1561){\makebox(0,0){$\bullet$}}
\put(5401,239){\makebox(0,0){$\bullet$}} \put
(5401,-661){\makebox(0,0){$\bullet$}}
\put(5401,-1561){\makebox(0,0){$\bullet$}}
\put(4500,-680){\makebox(0,0){$Id$}} \put(2800,-1861){\makebox
(0,0){$d_{(0,1)}d_{(0,1)}$}}
\put(2800,-699){\makebox(0,0){$d_{(1,1)} d_{(0,1)}$}}
\put(2800,464){\makebox(0,0){$d_{(0,1)}d_{(2,1)}$}}
\put(6200,-1861){\makebox(0,0){$d_{(0,1)}d_{(1,1)}$}} \put
(6200,-699){\makebox(0,0){$d_{(1,1)}d_{(1,1)}$}}
\put(6200,464){\makebox(0,0){$d_{(1,1)}d_{(2,1)}$}}
\put(3110,-1260){\makebox(0,0){$d_{(0,1)}$}}
\put(4550,530){\makebox(0,0){$d_{(2,1)}$}}
\put(2900,-111){\makebox(0,0){$d_{(0,1)(2,1)}$}}
\put(5900,-111){\makebox(0,0){$d_{(1,2)}$}}
\put(5900,-1260){\makebox(0,0){$d_{(1,1)}$}}
\put(4550,-1890){\makebox(0,0){$d_{(0,2)}$}}
\end{picture}\vspace{0.1in}

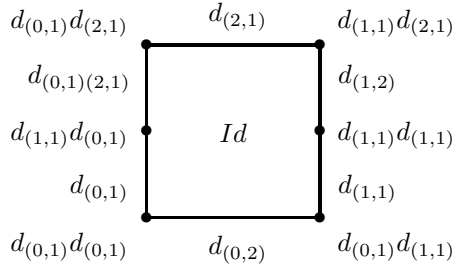
Figure 10: $P_{3}$ with face-operator labeling.\vspace{0.2in}
\end{center}

\noindent Since compositions of face operators are determined by the
correspondence between faces and partitions, we only label the codimension $1
$ faces of the related polytopes below.

The associahedra $\left\{  K_{n}\right\}  $ serve as parameter spaces for
higher homotopy associativ- ity.$\ $In his seminal papers of 1963
\cite{stasheff}, J. Stasheff constructed $K_{n}$ in the following way:$\ $Let
$K_{2}=\ast;$ if $K_{n-1}$ has been constructed, define $K_{n}$ to be the cone
on the set
\[
\bigcup\limits_{\substack{r+s=n+1 \\1\leq k\leq n-s+1}}\left(  K_{r}\times
K_{s}\right)  _{k}.
\]
Thus, $K_{n}$ is an $\left(  n-2\right)  $-dimensional convex polytope.\ 

Stasheff's motivating example of higher homotopy associativity in
\cite{stasheff} is the singular chain complex on the (Poincar\`{e}) loop space
of a connected $CW$-complex.$\ $Here associativity holds up to homotopy, the
homotopies between the various associations are homotopic, the homotopies
between these homotopies are homotopic, and so on.$\ $An abstract $A_{\infty}%
$-algebra is a DGA in which associativity behaves as in Stasheff's motivating
example.$\ $If $\varphi^{2}:A\otimes A\rightarrow A$ is the multiplication on
an $A_{\infty}$-algebra $A,$ the homotopies $\varphi^{n}:A^{\otimes
n}\rightarrow A$ are multilinear operations such that $\varphi^{3}$\ is a
chain homotopy between the associations $\left(  ab\right)  c$ and $a\left(
bc\right)  $ thought of as quadratic compositions $\varphi^{2}\left(
\varphi^{2}\otimes1\right)  $ and $\varphi^{2}\left(  1\otimes\varphi
^{2}\right)  $ in three variables, $\varphi^{4}$\ is a chain homotopy bounding
the cycle of five quadratic compositions in four variables involving
$\varphi^{2}$ and $\varphi^{3},$ and so on.$\ $Let $C_{\ast}\left(
K_{r}\right)  $ denote the cellular chains on $K_{r}.$ The natural
correspondence between faces of $K_{r}$ and the various compositions of the
$\varphi^{n}$'s in $r$ variables (modulo an appropriate equivalence) induces a
chain map $C_{\ast}\left(  K_{r}\right)  \rightarrow Hom\left(  A^{\otimes
r},A\right)  $ that determines the relations among the compositions of
$\varphi^{n}$'s.$\ $This chain map together with our diagonal on $K_{n}$ leads
to the tensor product of $A_{\infty}$-algebras (see Section 5).

Now if we disregard levels, a PLT is simply a planar rooted tree
(PRT).$\ $Quite remarkably, A. Tonks \cite{tonks} showed that $K_{n}$ is the
identification space $P_{n-1}/\sim$ in which all faces indexed by isomorphic
PRT's are identified.$\ $Since the quotient map $\theta:P_{n-1}\rightarrow
K_{n}$ is cellular, the faces of $K_{n}$ are indexed by PRT's with $n$ leaves.
The correspondence between PRT's with $n$ leaves and parenthesizations of $n$
indeterminants is simply this: Given a node $N,$ parenthesize the
indeterminants that correspond to leaves on all branches that meet at node
$N.$

\begin{example}
With one exception, all classes of faces of $P_{3}$ consist of a single
element.$\ $Elements of the exceptional class
\[
\left[  1|3|2,13|2,3|1|2\right]
\]
represent the parenthesization $\left(  \left(  \bullet\bullet\right)  \left(
\bullet\bullet\right)  \right)  .$ Whereas $1|3|2$ and $3|1|2$ insert inner
parentheses in the opposite order, the element $13|2$ inserts inner
parentheses simultaneously and represents a homotopy between $1|3|2$ and
$3|1|2$.$\ $Tonks' projection $\theta$ sends the exceptional class to the
vertex of $K_{4}$ represented by the parenthesization $\left(  \left(
\bullet\bullet\right)  \left(  \bullet\bullet\right)  \right)  .\ $The classes
of faces of $P_{4}$ with more than one element and their projections to
$K_{5}$ are:%
\[%
\begin{array}
[c]{l}%
\left[  12|4|3,124|3,4|12|3\right]  \overset{\theta}{\longrightarrow}\left(
\left(  \bullet\bullet\bullet\right)  \left(  \bullet\bullet\right)  \right)
\\
\left[  1|3|24,13|24,3|1|24\right]  \longrightarrow\left(  \left(
\bullet\bullet\right)  \left(  \bullet\bullet\right)  \bullet\right) \\
\left[  1|4|23,14|23,4|1|23\right]  \longrightarrow\left(  \left(
\bullet\bullet\right)  \bullet\left(  \bullet\bullet\right)  \right) \\
\left[  2|4|13,24|13,4|2|13\right]  \longrightarrow\left(  \bullet\left(
\bullet\bullet\right)  \left(  \bullet\bullet\right)  \right) \\
\left[  1|34|2,134|2,34|1|2\right]  \longrightarrow\left(  \left(
\bullet\bullet\right)  \left(  \bullet\bullet\bullet\right)  \right)
\end{array}
\]%
\[%
\begin{array}
[c]{l}%
\left[  1|3|2|4,13|2|4,3|1|2|4\right]  \longrightarrow\left(  \left(  \left(
\bullet\bullet\right)  \left(  \bullet\bullet\right)  \right)  \bullet\right)
\\
\left[  2|4|3|1,24|3|1,4|2|3|1\right]  \longrightarrow\left(  \bullet\left(
\left(  \bullet\bullet\right)  \left(  \bullet\bullet\right)  \right)  \right)
\\
\left[  1|2|4|3,1|24|3,1|4|2|3,14|2|3,4|1|2|3\right]  \longrightarrow\left(
\left(  \left(  \bullet\bullet\right)  \bullet\right)  \left(  \bullet
\bullet\right)  \right) \\
\left[  1|3|4|2,13|4|2,3|1|4|2,3|14|2,3|4|1|2\right]  \longrightarrow\left(
\left(  \bullet\bullet\right)  \left(  \left(  \bullet\bullet\right)
\bullet\right)  \right) \\
\left[  1|4|3|2,14|3|2,4|1|3|2,4|13|2,4|3|1|2\right]  \longrightarrow\left(
\left(  \bullet\bullet\right)  \left(  \bullet\left(  \bullet\bullet\right)
\right)  \right) \\
\left[  2|1|4|3,2|14|3,2|4|1|3,24|1|3,4|2|1|3\right]  \longrightarrow\left(
\left(  \bullet\left(  \bullet\bullet\right)  \right)  \left(  \bullet
\bullet\right)  \right)  .
\end{array}
\]
\newline Faces and edges represented by elements of the first five classes
project to edges; edges and vertices represented by elements of the next six
classes project to vertices.
\end{example}

The multiplihedra $\left\{  J_{n+1}\right\}  $, which serve as parameter
spaces for homotopy multiplicative morphisms of $A_{\infty}$-algebras, lie
between the associahedra and permutahedra (see \cite{stasheff}, \cite{Iwase}%
).\ If $f^{1}:A\rightarrow B$ is such a morphism, there is a chain homotopy
$f^{2}$ between the quadratic compositions $f^{1}\varphi_{A}^{2}$ and
$\varphi_{B}^{2}\left(  f^{1}\otimes f^{1}\right)  $ in two variables, there
is a chain homotopy $f^{3}$ bounding the cycle of the six quadratic
compositions in three variables involving $f^{1},$ $f^{2},$ $\varphi_{A}^{2},$
$\varphi_{A}^{3},$ $\varphi_{B}^{2}$ and $\varphi_{B}^{3},$ and so on.$\ $The
natural correspondence between faces of $J_{r}$ and the various compositions
of $f^{i},$ $\varphi_{A}^{j}$ and $\varphi_{B}^{k}$ in $r$ variables (modulo
an appropriate equivalence) induces a chain map $C_{\ast}\left(  J_{r}\right)
\rightarrow Hom\left(  A^{\otimes r},B\right)  $.

The multiplihedron $J_{n+1}$ can also be realized as a subdivision of the cube
$I^{n}.\ $For $n=0,1,2,$ set $J_{n+1}=P_{n+1}.$ If $J_{n}$ has been
constructed, $J_{n+1}$ is the subdivision of $J_{n}\times I$ given below and
its various $\left(  n-1\right)  $-faces are labeled as indicated:%
\[%
\begin{tabular}
[c]{c|c}%
$\underset{\ }{\text{\textbf{Face of }}J_{n+1}}$ & \textbf{Face operator}%
\\\hline
$\underset{\ }{\overset{\ }{J_{n}\times0}}$ & $d_{\left(  0,n\right)  }$\\
$\underset{\ }{J_{n}\times1}$ & $d_{\left(  n,1\right)  }$\\
$\underset{\ }{d_{\left(  i,\ell\right)  }\times I}$ & $d_{\left(
i,\ell\right)  },\hspace{0.16in}1\leq i<n-\ell$\\
$\underset{\ }{d_{\left(  i,\ell\right)  }\times\left[  0,1-2^{-i}\right]  }$
& $d_{\left(  i,\ell\right)  },\hspace{0.16in}1\leq i=n-\ell$\\
$\underset{\ }{d_{\left(  i,\ell\right)  }\times\left[  1-2^{-i},1\right]  }$
& $d_{\left(  i,\ell+1\right)  },$ $1\leq i=n-\ell$\\
$\underset{\ }{d_{\left(  0,\ell_{1}\right)  \cdots\left(  i_{k},\ell
_{k}\right)  }\times\left[  0,1-2^{k-n}\right]  }$ & $d_{\left(  0,\ell
_{1}\right)  \cdots\left(  i_{k},\ell_{k}\right)  }$\\
$d_{\left(  0,\ell_{1}\right)  \cdots\left(  i_{k},\ell_{k}\right)  }%
\times\left[  1-2^{k-n},1\right]  $ & $\left\{
\begin{array}
[c]{ll}%
d_{\left(  0,\ell_{1}\right)  \cdots\left(  i_{k},\ell_{k}\right)  \left(
n,1\right)  }, & i_{k}<n-\ell_{k}\\
d_{\left(  0,\ell_{1}\right)  \cdots\left(  i_{k},\ell_{k}+1\right)  }, &
i_{k}=n-\ell_{k}.
\end{array}
\right.  $%
\end{tabular}
\]
\vspace{0.1in}

\noindent Thus faces of $J_{n+1}$ are indexed by compositions of face
operators of the form
\begin{equation}
d_{\left(  i_{m},\ell_{m}\right)  }\cdots d_{(i_{k_{1}},\ell_{k_{1}}%
)\cdots\left(  i_{k_{s}},\ell_{k_{s}}\right)  }\cdots d_{\left(  i_{1}%
,\ell_{1}\right)  }. \label{multface}%
\end{equation}
In terms of trees and parenthesizations this says the following:$\ $Let $T$ be
a $\left(  k+1\right)  $-leveled tree with left-most branch attached at level
$p.\ $For $1\leq j<p,$ insert level $j$ parentheses one pair at a time without
regard to order as in $K_{n+2};$ next, insert all level $p$ parentheses
simultaneously as in $P_{n+1};$ finally, for $j>p,$ insert level $j$
parentheses one pair at a time without regard to order.$\ $Thus multiple lower
indices in a composition of face operators may only occur when the left-most
branch is attached above the root.$\ $This suggests the following equivalence
relation on the set of $\left(  k+1\right)  $-leveled trees with $n+2$
leaves:$\ $Let $T$ and $T^{\prime}$ be $p$-leveled trees with $n+2$ nodes
whose $p$-level meets $U_{p}$ and $U_{p}^{\prime}$ contain $1.\ $Then $T\sim
T^{\prime}$ if $T$ and $T^{\prime}$ are isomorphic as PLT's and $U_{p}%
=U_{p}^{\prime}.\ $This equivalence relation induces a cellular projection
$\pi:P_{n+1}\rightarrow J_{n+1}$ under which $J_{n+1}$ can be realized as an
identification space of $P_{n+1}.\ $Furthermore, the projection $J_{n+1}%
\rightarrow K_{n+2}$ given by identifying faces of $J_{n+1}$ indexed by
isomorphic PLT's gives the factorization $P_{n+1}\overset{\pi}{\longrightarrow
}J_{n+1}\rightarrow K_{n+2}$ of Tonks' projection.

It is interesting to note the role of the indices $\ell_{j}$ in compositions
of face operators representing the faces of $J_{n+1}$ as in (\ref{multface}%
).$\ $With one exception, each $U_{j}$ in the corresponding partition
$U_{1}|\cdots|U_{m+1}$ is a set of consecutive integers; this holds without
exception for all $U_{j}$ on $K_{n+2}$.$\ $The exceptional set $U_{p}$ is a
union of $s$ sets of consecutive integers with maximal cardinality, as is
typical of sets $U_{j}$ on $P_{n+1}$.$\ $Thus $J_{n+1}$ exhibits
characteristics of both combinatorial structures.\vspace{0.1in}

\begin{center}
\setlength{\unitlength}{0.0099in}\begin{picture} (480,-400)
\thicklines
\put(0,-120){\line( 0,-1){120}} \put(120,0){\line( 0,-1){360}}
\put (180,-120){\line( 0,-1){120}} \put(240,0){\line( 0,-1){360}}
\put (360,-120){\line( 0,-1){120}} \put(420,-120){\line(
0,-1){120}} \put (480,-120){\line( 0,-1){120}}
\put(120,0){\line( 1,0){120}} \put(0,-120){\line( 1,0){480}} \put
(180,-180){\line( 1,0){60}} \put(240,-150){\line( 1,0){120}} \put
(420,-150){\line( 1,0){60}} \put(0,-150){\line( 1,0){120}} \put
(360,-180){\line( 1,0){60}} \put(0,-240){\line( 1,0){480}} \put
(120,-360){\line( 1,0){120}}
\put(120,0){\makebox(0,0){$\bullet$}}
\put(180,0){\makebox(0,0){$\bullet$}}
\put(240,0){\makebox(0,0){$\bullet$}}
\put(0,-120){\makebox(0,0){$\bullet$}}
\put(120,-120){\makebox(0,0){$\bullet$}} \put(180,-120){\makebox
(0,0){$\bullet$}} \put(240,-120){\makebox(0,0){$\bullet$}} \put
(360,-120){\makebox(0,0){$\bullet$}}
\put(420,-120){\makebox(0,0){$\bullet$}}
\put(480,-120){\makebox(0,0){$\bullet$}} \put(240,-150){\makebox
(0,0){$\bullet$}} \put(360,-150){\makebox(0,0){$\bullet$}} \put
(180,-180){\makebox(0,0){$\bullet$}}
\put(240,-180){\makebox(0,0){$\bullet$}}
\put(420,-180){\makebox(0,0){$\bullet$}}
\put(0,-150){\makebox(0,0){$\bullet$}}
\put(120,-240){\makebox(0,0){$\bullet$}} \put(360,-240){\makebox
(0,0){$\bullet$}} \put(420,-240){\makebox(0,0){$\bullet$}} \put
(480,-150){\makebox(0,0){$\bullet$}}
\put(0,-240){\makebox(0,0){$\bullet$}}
\put(120,-150){\makebox(0,0){$\bullet$}} \put(180,-240){\makebox
(0,0){$\bullet$}} \put(240,-240){\makebox(0,0){$\bullet$}} \put
(360,-180){\makebox(0,0){$\bullet$}}
\put(420,-150){\makebox(0,0){$\bullet$}}
\put(480,-240){\makebox(0,0){$\bullet$}} \put(120,-360){\makebox
(0,0){$\bullet$}} \put(180,-360){\makebox(0,0){$\bullet$}} \put
(240,-360){\makebox(0,0){$\bullet$}}
\put(275,-100) {\makebox(0,0){$(1,1,1)$}} \put(20,-260) {\makebox
(0,0){$(0,0,0)$ }} \put(177,-55) {\makebox(0,0){$d_{(3,1)}$}}
\put(177,-295) {\makebox(0,0){$d_{(0,3)}$}} \put(55,-134)
{\makebox(0,0){$d_{(0,2)(3,1)}$}} \put(55,-180)
{\makebox(0,0){$d_{(0,2)}$}} \put(390,-200) {\makebox
(0,0){$d_{(0,1)(2,1)}$}} \put(390,-134)
{\makebox(0,0){$d_{(0,1)(2,2)}$}} \put(300,-134)
{\makebox(0,0){$d_{(2,2)}$}} \put(150,-180) {\makebox
(0,0){$d_{(1,1)}$}} \put(300,-195) {\makebox(0,0){$d_{(2,1)}$}}
\put(210,-152) {\makebox(0,0){$d_{(1,3)}$}} \put(450,-134)
{\makebox (0,0){$d_{(0,1)(3,1)}$}} \put(210,-210)
{\makebox(0,0){$d_{(1,2)}$}} \put(450,-194)
{\makebox(0,0){$d_{(0,1)}$}}
\end{picture}\vspace{3.7in}

Figure 11: $J_{4}$ as a subdivision of $J_{3}\times I.$\vspace{0.1in}
\end{center}

We realize the associahedron $K_{n+2}$ in a similar way.$\ $For $n=0,1,$ set
$K_{n+2}=P_{n+1}.\ $If $K_{n+1}$ has been constructed, let $e_{i,\epsilon}$
denote the face $(x_{1},\ldots,x_{i-1},\epsilon,x_{i+1},$ $\ldots
,x_{n})\subset I^{n},$ where $\epsilon=0,1$ and $1\leq i\leq n.$ Then
$K_{n+2}$ is the subdivision of $K_{n+1}\times I$ given below and its various
$\left(  n-1\right)  $-faces are labeled as indicated:%
\[%
\begin{tabular}
[c]{c|cc}%
$\underset{\ }{\text{\textbf{Face of }}K_{n+2}}$ & \textbf{Face operator} &
\\\hline
$e_{\ell,0}$ & $d_{\left(  0,\ell\right)  },$ & \multicolumn{1}{l}{$1\leq
\ell\leq n$}\\
$e_{n,1}$ & $d_{\left(  n,1\right)  ,}$ & \multicolumn{1}{l}{}\\
$d_{\left(  i,\ell\right)  }\times I$ & $d_{\left(  i,\ell\right)  ,}$ &
\multicolumn{1}{l}{$1\leq i<n-\ell\text{ }$}\\
&  & \\
$d_{\left(  i,\ell\right)  }\times\left[  0,1-2^{-i}\right]  $ & $d_{\left(
i,\ell\right)  ,}$ & \multicolumn{1}{l}{$1\leq i=n-\ell$}\\
&  & \\
$d_{\left(  i,\ell\right)  }\times\left[  1-2^{-i},1\right]  $ & $d_{\left(
i,\ell+1\right)  },$ & \multicolumn{1}{l}{$1\leq i=n-\ell$}%
\end{tabular}
\]
\vspace{0.1in}

\begin{center}
\setlength{\unitlength}{0.0004in}\begin{picture}
(2775,2685)(2926,-2038) \thicklines\put(3401,239){\line( 1,
0){1800}} \put(3401,239){\makebox(0,0){$\bullet$}}
\put(5201,239){\line( 0, 1){ 0}}
\put(5201,239){\makebox(0,0){$\bullet$}} \put(5201,239){\line(
0,-1){1800}} \put(5201,-661){\makebox(0,0){$\bullet$}}
\put(5201,-1561){\line(-1, 0){1800}}
\put(5201,-1561){\makebox(0,0){$\bullet$}} \put(3401,-1561){\line(
0, 1){1800}} \put(3401,-1561){\makebox(0,0){$\bullet$}}
\put(3401,239){\line( 0, 1){ 0}}
\put(4300,-680){\makebox(0,0){$Id$}}
\put(2860,-680){\makebox(0,0){$d_{(0,1)}$}}
\put(4350,600){\makebox(0,0){$d_{(2,1)}$}}
\put(5800,-211){\makebox(0,0){$d_{(1,2)}$}}
\put (5800,-1060){\makebox(0,0){$d_{(1,1)}$}}
\put(4350,-1890){\makebox (0,0){$d_{(0,2)}$}}
\end{picture}\vspace{0.1in}

Figure 12: $K_{4}$ as a subdivision of $K_{3}\times I.$ \vspace{0.3in}

\setlength{\unitlength}{0.007in}\begin{picture}
(480,0) \thicklines
\put(0,-120){\line( 0,-1){120}} \put(120,0){\line( 0,-1){360}}
\put (180,-120){\line( 0,-1){120}} \put(240,0){\line( 0,-1){360}}
\put (360,-120){\line( 0,-1){120}} \put(480,-120){\line(
0,-1){120}}
\put(120,0){\line( 1,0){120}} \put(0,-120){\line( 1,0){480}} \put
(180,-180){\line( 1,0){60}} \put(240,-150){\line( 1,0){120}} \put
(0,-240){\line( 1,0){480}} \put(120,-360){\line( 1,0){120}}
\put(120,0){\makebox(0,0){$\bullet$}}
\put(240,0){\makebox(0,0){$\bullet$}}
\put(0,-120){\makebox(0,0){$\bullet$}}
\put(120,-120){\makebox(0,0){$\bullet$}}
\put(180,-120){\makebox(0,0){$\bullet$}} \put(240,-120){\makebox
(0,0){$\bullet$}} \put(360,-120){\makebox(0,0){$\bullet$}} \put
(480,-120){\makebox(0,0){$\bullet$}}
\put(240,-150){\makebox(0,0){$\bullet$}}
\put(360,-150){\makebox(0,0){$\bullet$}} \put(180,-180){\makebox
(0,0){$\bullet$}} \put(240,-180){\makebox(0,0){$\bullet$}} \put
(120,-240){\makebox(0,0){$\bullet$}}
\put(360,-240){\makebox(0,0){$\bullet$}}
\put(0,-240){\makebox(0,0){$\bullet$}}
\put(180,-240){\makebox(0,0){$\bullet$}}
\put(240,-240){\makebox(0,0){$\bullet$}} \put(480,-240){\makebox
(0,0){$\bullet$}} \put(120,-360){\makebox(0,0){$\bullet$}} \put
(240,-360){\makebox(0,0){$\bullet$}}
\put(280,-100){\makebox(0,0){$(1,1,1)$}}
\put(20,-260){\makebox(0,0){$(0,0,0)$ }}
\put(177,-55){\makebox(0,0){$d_{(3,1)}$}} \put(177,-295){\makebox
(0,0){$d_{(0,3)}$}} \put(55,-180) {\makebox(0,0){$d_{(0,2)}$}}
\put (300,-134){\makebox(0,0){$d_{(2,2)}$}}
\put(150,-180){\makebox(0,0){$d_{(1,1)} $}} \put(300,-195)
{\makebox(0,0){$d_{(2,1)}$}} \put(210,-152){\makebox
(0,0){$d_{(1,3)}$}} \put(210,-210){\makebox(0,0){$d_{(1,2)}$}}
\put (420,-180){\makebox(0,0){$d_{(0,1)}$}}
\end{picture}\vspace{2.6in}

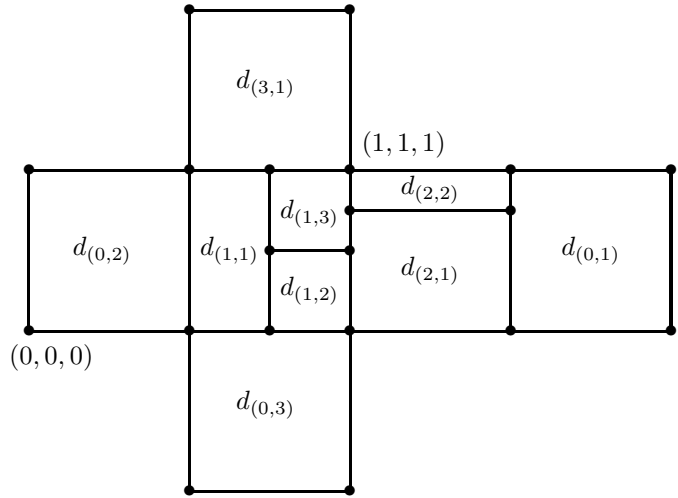
Figure 13: $K_{5}$ as a subdivision of $K_{4}\times I.$
\end{center}

\section{Diagonals on the Associahedra and Multiplihedra}

The diagonal $\Delta_{P}$ on $C_{\ast}\left(  P_{n+1}\right)  $ descends to
diagonals $\Delta_{J}$ on $C_{\ast}\left(  J_{n+1}\right)  $ and $\Delta_{K} $
on $C_{\ast}\left(  K_{n+2}\right)  $ via the cellular projections
$\pi:P_{n+1}\rightarrow J_{n+1}$ and $\theta:P_{n+1}\rightarrow K_{n+2}$
discussed in Section 2 above.$\ $This fact is an immediate consequence of
Proposition \ref{delta-w}.

\begin{definition}
\label{degen}Let $f:W\rightarrow X$ be a cellular map of CW-complexes, let
$\Delta_{W}$ be a diagonal on $C_{\ast}\left(  W\right)  $ and let $X^{\left(
r\right)  }$ denote the $r$-skeleton of $X$.$\ $A $k$-cell $e\subseteq W$ is
\underline{degenerate under $f$} if $f(e)\subseteq X^{\left(  r\right)  }$
with $r<k.\ $A component $a\otimes b$ of $\Delta_{W}$ is \underline{degenerate
under} $f$ if either $a$ or $b$ is degenerate under $f.$
\end{definition}

\noindent Let us identify the non-degenerate cells of $P_{n+1}$ under $\pi$
and $\theta.$

\begin{definition}
Let $A_{1}|\cdots|A_{p}$ be a partition of $\underline{n+1}$ with $p>1$ and
let $1\leq k<p.\ $The subset $A_{k}$ is \underline{exceptional} if for
$k<j\leq p,$ there is an element $a_{i,j}\in A_{j}\ $such that $\min
A_{k}<a_{i,j}<\max A_{k}.$
\end{definition}

\begin{proposition}
Let $a=A_{1}|\cdots|A_{p}$ be a face of $P_{n+1}$ and let
\[
d_{(i_{1}^{p-1},\ell_{1}^{p-1})\cdots(i_{s_{p-1}}^{p-1},\ell_{s_{p-1}}^{p-1}%
)}\cdots d_{(i_{1}^{1},\ell_{1}^{1})\cdots(i_{s_{1}}^{1},\ell_{s_{1}}^{1})}%
\]
be its unique representation as a composition of face operators.\vspace{0.1in}
\newline(1) The following are all equivalent:\vspace{0.1in} \newline%
\hspace*{0.25in}(1a) The face $a$ is degenerate under $\pi.$\vspace{0.1in}
\newline\hspace*{0.25in}(1b) $\min A_{j}>\min\left(  A_{j+1}\cup\cdots\cup
A_{p}\right)  $ with $A_{j}$ exceptional for some $j<p.$\vspace{0.1in}
\newline\hspace*{0.25in}(1c) $i_{1}^{k}>0$ and $s_{k}>1$ for some
$k<p.$\vspace{0.1in}\newline(2) The following are all equivalent:\vspace
{0.1in}\newline\hspace*{0.25in}(2a) The face $a$ is degenerate under $\theta
.$\vspace{0.1in}\newline\hspace*{0.25in}(2b) $A_{j}$ is exceptional for some
$j<p.$\vspace{0.1in}\newline\hspace*{0.25in}(2c) $s_{k}>1$ for some
$k<p.$\vspace{0.1in}
\end{proposition}

\begin{proof}
Obvious.
\end{proof}

\begin{example}
The subset $A_{1}=\left\{  13\right\}  $ in the partition $a=13|24$ is
exceptional and the face $a\subset P_{4}$ is degenerate under $\theta.\ $In
terms of compositions of face operators, the face $a$ corresponds
to\ $d_{\left(  0,1\right)  \left(  2,1\right)  }\left(  x_{1}\cdots
x_{5}\right)  $ with $s_{1}=2.\ $Furthermore, $a$ is also non-degenerate under
$\pi$ since $i_{1}^{1}=0$ (equivalently, $\min A_{1}<\min A_{2}$).
\end{example}

Next, we apply Tonks' projection and obtain an explicit formula for the
diagonal $\Delta_{K}$ on the associahedra.

\begin{proposition}
\label{delta-w}Let $f:W\rightarrow X$ be a surjective cellular map and let
$\Delta_{W}$ be a diagonal on $C_{\ast}\left(  W\right)  .\ $Then $\Delta_{W}
$ uniquely determines a diagonal $\Delta_{X}$ on $C_{\ast}\left(  X\right)  $
given by the non-degenerate components of $\Delta_{W}$ under $f.\ $Moreover,
$\Delta_{X}$ is the unique map that commutes the following diagram:%
\[%
\begin{array}
[c]{ccc}%
C_{\ast}\left(  W\right)  & \overset{\Delta_{W}}{\longrightarrow} & C_{\ast
}\left(  W\right)  \otimes C_{\ast}\left(  W\right) \\
f\downarrow &  & \downarrow f\otimes f\\
C_{\ast}\left(  X\right)  & \underset{\Delta_{X}}{\longrightarrow} & C_{\ast
}\left(  X\right)  \otimes C_{\ast}\left(  X\right)  .
\end{array}
\]

\end{proposition}

\begin{proof}
Obvious.
\end{proof}

In Section 2 we established correspondences between faces of $P_{n+1}$ and
PLT's with $n+2$ leaves and between faces of $K_{n+2}$ and PRT's with $n+2$
leaves. Since a PRT can be viewed as a PLT, faces of $K_{n+2}$ can be viewed
as faces of $P_{n+1}$.

\begin{definition}
\label{induced}For $n\geq0,$ let $\Delta_{P}$ be the diagonal on $C_{\ast
}\left(  P_{n+1}\right)  $ and let $\theta:P_{n+1}\rightarrow K_{n+2}$ be
Tonks' projection.$\ $View each face $e$ of the associahedron $K_{n+2}$ as a
face of $P_{n+1}$ and define $\Delta_{K}:C_{\ast}\left(  K_{n+2}\right)
\rightarrow C_{\ast}\left(  K_{n+2}\right)  \otimes C_{\ast}\left(
K_{n+2}\right)  $ by%

\[
\Delta_{K}(e)=(\theta\otimes\theta)\Delta_{P}(e).
\]

\end{definition}

\begin{corollary}
The map $\Delta_{K}$ given by Definition \ref{induced} is the diagonal on
$C_{\ast}\left(  K_{n+2}\right)  $ induced by $\Delta_{P}$.
\end{corollary}

\begin{proof}
This is an immediate application of Proposition \ref{delta-w}.
\end{proof}

Consider a CP $u\otimes v=c\left(  F\right)  \otimes r\left(  F\right)  $
related to SCP $a\otimes b=c\left(  E\right)  \otimes r\left(  E\right)  $ via
$F=D_{N_{q-1}}\cdots D_{1}R_{M_{p-1}}\cdots R_{1}E.$ Note that both factors of
$u\otimes v$ are non-degenerate under $\theta$ if and only if $b$ is
non-degenerate and each $M_{j}$ has maximal cardinality.$\ $Alternatively, if
$d^{p}\cdots d^{1}\otimes d^{q}\cdots d^{1}\left(  e^{n}\otimes e^{n}\right)
$ is a component of $\Delta_{K}\left(  e^{n}\right)  ,$ factors in the
corresponding pairing $T_{p}\otimes T_{q}$ of PRT's have $n+2$ leaves, $p+1$
and $q+1$ nodes and respective dimensions $n-p$ and $n-q.$ Hence $p+q=n$ and
$T_{p}\otimes T_{q}$ has exactly $n+2$ nodes. But if $T_{u}\otimes T_{v}$ is
the pairing of PLT's corresponding to $u\otimes v$, forgetting levels in
$T_{u}\otimes T_{v}$ gives the pairing of PRT's corresponding to
$\theta\left(  u\right)  \otimes\theta\left(  v\right)  $. Since the number of
nodes in $T_{u}\otimes T_{v}$ is at least $n+2,$ $\theta\left(  a\right)
\otimes\theta\left(  b\right)  $ is non-degenerate in $\Delta_{K}\left(
e^{n}\right)  $ if and only if the total number of nodes in $T_{u}\otimes
T_{v}$ is exactly $n+2.$

Choose a system of generators $e^{n}\in C_{n}\left(  K_{n+2}\right)  ,$
$n\geq0$.$\ $The signs in (\ref{delta-Kn}) below follow from (\ref{diagonal}).

\begin{definition}
[\cite{SU1}]\label{Delta-K}For each $n\geq0,$ define $\Delta_{K}$ on $e^{n}\in
C_{n}\left(  K_{n+2}\right)  $ by
\begin{equation}
\Delta_{K}\left(  e^{n}\right)  =\sum_{0\leq p\leq p+q=n+2}\hspace
{-0.1in}\left(  -1\right)  ^{\epsilon}\text{ }d_{(i_{p-1},\ell_{p-1})}\cdots
d_{(i_{1},\ell_{1})}\otimes d_{(i_{q-1}^{\prime},\ell_{q-1}^{\prime})}\cdots
d_{(i_{1}^{\prime},\ell_{1}^{\prime})}\left(  e^{n}\otimes e^{n}\right)  ,
\label{delta-Kn}%
\end{equation}
where
\[
\epsilon=\sum_{j=1}^{p-1}i_{j}(\ell_{j}+1)+\sum_{k=1}^{q-1}(i_{k}^{\prime
}+k+q)\ell_{k}^{\prime},
\]
and lower indices $\left(  \left(  i_{1},\ell_{1}\right)  ,\ldots,\left(
i_{p-1},\ell_{p-1}\right)  ;\left(  i_{1}^{\prime},\ell_{1}^{\prime}\right)
,\ldots,\left(  i_{q-1}^{\prime},\ell_{q-1}^{\prime}\right)  \right)  $ range
over all solutions of the following system of inequalities:\vspace{0.1in}
\begin{equation}
\hspace*{0.4in}\left\{
\begin{array}
[c]{ll}%
1\leq i_{j}^{\prime}<i_{j-1}^{\prime}\leq n+1\vspace{0.1in} & (1)\\
1\leq\ell_{j}^{\prime}\leq n+1-i_{j}^{\prime}-\ell_{\left(  j-1\right)
}^{\prime}\vspace{0.1in} & (2)\\
0\leq i_{k}\leq\min\limits_{o^{\prime}\left(  t_{k}\right)  <r<k}\left\{
i_{r},\text{ }i_{t_{k}}^{\prime}-\ell_{\left(  o^{\prime}\left(  t_{k}\right)
\right)  }\right\}  \vspace{0.1in} & (3)\\
1\leq\ell_{k}=\epsilon_{k}-i_{k}-\ell_{\left(  k-1\right)  }, & (4)
\end{array}
\right\}  _{\substack{1\leq k\leq p-1 \\1\leq j\leq q-1}} \label{system}%
\end{equation}
where
\end{definition}

$\hspace*{1.02in}\left\{  \epsilon_{1}<\cdots<\epsilon_{q-1}\right\}
=\left\{  1,\ldots,n\right\}  \setminus\left\{  i_{1}^{\prime},\ldots
,i_{q-1}^{\prime}\right\}  ;$\vspace{0.1in}

$\hspace*{1.02in}\epsilon_{0}=\ell_{0}=\ell_{0}^{\prime}=i_{p}=i_{q}^{\prime
}=0;$\vspace{0.1in}

$\hspace*{1.02in}i_{0}=i_{0}^{\prime}=\epsilon_{q}=\ell_{\left(  p\right)
}=\ell_{\left(  q\right)  }^{\prime}=n+1;$\vspace{0.1in}

$\hspace*{1.02in}\ell_{\left(  u\right)  }=\sum\nolimits_{j=0}^{u}\ell_{j}$
for $0\leq u\leq p;$\vspace{0.1in}

$\hspace*{1.02in}\ell_{\left(  u\right)  }^{\prime}=\sum\nolimits_{k=0}%
^{u}\ell_{k}^{\prime}$ for $0\leq u\leq q;$\vspace{0.1in}

$\hspace*{1.02in}t_{u}=\min\left\{  r\text{ }|\text{ }i_{r}^{\prime}%
+\ell_{\left(  r\right)  }^{\prime}-\ell_{\left(  o\left(  u\right)  \right)
}^{\prime}>{\epsilon}_{u}>i_{r}^{\prime}\right\}  ;$\vspace{0.1in}

$\hspace*{1.02in}o\left(  u\right)  =\max\left\{  r\text{ }|\text{ }%
i_{r}^{\prime}\geq\epsilon_{u}\right\}  ;$ and\vspace{0.1in}

$\hspace*{1.02in}o^{\prime}\left(  u\right)  =\max\left\{  r\text{ }|\text{
}\epsilon_{r}\leq i_{u}^{\prime}\right\}  .$\vspace{0.1in}\newline

\noindent\textit{\ Extend }$\Delta_{K}$\textit{\ to proper faces of }$K_{n+2}$
\textit{via the standard comultiplicative extension. }

\begin{theorem}
The map $\Delta_{K}$ given by Definition \ref{Delta-K} is the diagonal induced
by $\theta$.
\end{theorem}

\begin{proof}
If $v=L_{\beta}(v^{\prime})$ is non-degenerate in some component $u\otimes v$
of $\Delta_{P}$, then so is $v^{\prime},$ and we immediately obtain inequality
(1) of (\ref{system}).$\ $Next, each non-degenerate decreasing $b$ uniquely
determines an SCP $a\otimes b$.$\ $Although $a$ may be degenerate, there is a
unique non-degenerate $u=R_{M_{p-1}}\cdots R_{M_{1}}(a)$ obtained by choosing
each $M_{j}$ with maximal cardinality (the case $M_{j}=\varnothing$ for all
$j$ may nevertheless occur);$\ $then $u\otimes b$ is a non-degenerate CP
associated with $a\otimes b$ in $\Delta_{P}.\ $As a composition of face
operators, straightforward examination shows that $u$ has form $u=d_{(i_{p-1}%
,\ell_{p-1})}\cdots d_{(i_{1},\ell_{1})}(e^{n})$ and is related to
$b=d_{(i_{q-1}^{\prime},\ell_{q-1}^{\prime})}\cdots d_{(i_{1}^{\prime}%
,\ell_{1}^{\prime})}\left(  e^{n}\right)  $ by
\[
i_{k}=\min\limits_{o^{\prime}\left(  t_{k}\right)  <r<k}\left\{  i_{r},\text{
}i_{t_{k}}^{\prime}-\ell_{\left(  o^{\prime}\left(  t_{k}\right)  \right)
}\right\}  ,\ 1\leq k<p;
\]
and equality holds in (4) of (\ref{system}).$\ $Finally, let $b=L_{\beta}%
(\bar{b}).\ $As we vary $\bar{b}$ in all possible ways, each $\bar{b}$ is
non-degenerate and we obtain all possible non-degenerate CP's $\bar{u}%
\otimes\bar{b}$ associated with $\bar{a}\otimes\bar{b}$ ($\bar{u}=u$ when
$\bar{b}=b$ and $\beta=\varnothing$)$.\ $For each such $\bar{u}=d_{(i_{p-1}%
,\ell_{p-1})}\cdots d_{(i_{1},\ell_{1})}(e^{n})$ we have both inequality (3)
and equality in (4) of (\ref{system}).$\ $Hence, the theorem is
proved.$\vspace{0.15in}$
\end{proof}

\begin{example}
On $K_{4}$ we obtain:
\begin{align*}
\Delta_{K}\left(  e^{2}\right)   &  =\left\{  d_{\left(  0,1\right)
}d_{\left(  0,1\right)  }\otimes1+1\otimes d_{\left(  1,1\right)  }d_{\left(
2,1\right)  }+d_{\left(  0,2\right)  }\otimes d_{\left(  1,1\right)  }\right.
\\
&  \qquad\qquad\left.  +d_{\left(  0,2\right)  }\otimes d_{\left(  1,2\right)
}^{\mathstrut}+d_{\left(  1,1\right)  }\otimes d_{\left(  1,2\right)
}^{\mathstrut}-d_{\left(  0,1\right)  }\otimes d_{\left(  2,1\right)
}^{\mathstrut}\right\}  (e^{2}\otimes e^{2}).
\end{align*}

\end{example}

\section{Application: Tensor Products of $A_{\infty}$-(co)algebras}

In this section, we use $\Delta_{K}$ to define the tensor product of
$A_{\infty}$-(co)algebras in maximal generality.$\ $We note that a special
case was given by J. Smith \cite{Smith} for certain objects with a richer
structure than we have here.$\ $We also mention that Lada and Markl
\cite{Lada} defined an $A_{\infty}$ tensor product structure on a construct
different from the tensor product of graded modules.

We adopt the following notation and conventions: Let $R$ be a commutative ring
with unity; $R$-modules are assumed to be $\mathbb{Z}$-graded, tensor products
and $Hom$'s are defined over $R$ and all maps are $R$-module maps unless
otherwise indicated.$\ $If an $R$-module $V$ is connected, $\overline
{V}=V/V_{0}$.$\ $The symbol $1:V\rightarrow V$ denotes the identity map; the
suspension and desuspension maps $\uparrow$ and $\downarrow$ shift dimension
by $+1$ and $-1$, respectively.$\ $Define $V^{\otimes0}=R$ and $V^{\otimes
n}=V\otimes\cdots\otimes V$ with $n>0$ factors; then $TV=\oplus_{n\geq
0}V^{\otimes n}$ and $T^{a}V$ (respectively, $T^{c}V$) denotes the free tensor
algebra (respectively, cofree tensor coalgebra) of $V.\ $Given $R$-modules
$V_{1},\ldots,V_{n},$ a permutation $\sigma\in S_{n}$ induces an isomorphism
$\sigma:V_{1}\otimes\cdots\otimes V_{n}\rightarrow V_{\sigma^{-1}(1)}%
\otimes\cdots\otimes V_{\sigma^{-1}(n)}$ by $\sigma(x_{1}\cdots x_{n})=\pm$
$x_{\sigma^{-1}(1)}\cdots x_{\sigma^{-1}(n)},$ where $\pm$ is the Koszul sign.
In particular, $\sigma_{2,n}=\left(  1\text{ }3\text{ }\cdots\text{ }\left(
2n-1\right)  \text{ }2\text{ }4\text{ }\cdots\text{ }2n\right)  :\left(
A\otimes B\right)  ^{\otimes n}\rightarrow A^{\otimes n}\otimes B^{\otimes n}$
and $\sigma_{n,2}=\sigma_{2,n}^{-1}$ induce isomorphisms $\left(  \sigma
_{2,n}\right)  ^{\ast}:Hom\left(  A^{\otimes n}\otimes B^{\otimes n},A\otimes
B\right)  \rightarrow Hom\left(  \left(  A\otimes B\right)  ^{\otimes
n}\hspace*{-0.02in},A\otimes B\right)  $ and $\left(  \sigma_{n,2}\right)
_{\ast}$:$\,Hom\left(  A\otimes B,A^{\otimes n}\otimes B^{\otimes n}\right)
\rightarrow Hom\left(  A\otimes B,\right.  $ $\left.  \left(  A\otimes
B\right)  ^{\otimes n}\right)  $. The map $\iota:Hom(U,V)\otimes Hom\left(
U^{\prime},V^{\prime}\right)  \rightarrow Hom\left(  U\otimes U^{\prime
},V\otimes V^{\prime}\right)  $ is the canonical isomorphism. If $f:V^{\otimes
p}\rightarrow V^{\otimes q}$ is a map, we let $f_{i,n-p-i}=1^{\otimes
i}\otimes f\otimes1^{\otimes n-p-i}:V^{\otimes n}\rightarrow V^{\otimes
n-p+q},$ where $0\leq i\leq n-p.\ $The abbreviations \textit{DGM, DGA,}and
\textit{DGC} stand for \textit{differential graded }$R$\textit{-module, DG
}$R$-\textit{algebra }and \textit{DG }$R$-\textit{coalgebra}, respectively.

We begin with a review of $A_{\infty}$-(co)algebras paying particular
attention to the signs.$\ $Let $A$ be a connected $R$-module equipped with
operations $\{\varphi^{k}\in$\linebreak$Hom^{k-2}\left(  A^{\otimes
k},A\right)  \}_{k\geq1}.\ $For each $k$ and $n\geq1,$ linearly extend
$\varphi^{k}$ to $A^{\otimes n}$ via
\[
\sum\limits_{i=0}^{n-k}\varphi_{i,n-k-i}^{k}:A^{\otimes n}\rightarrow
A^{\otimes n-k+1},
\]
and consider the induced map of degree $-1$ given by
\[
\sum\limits_{i=0}^{n-k}\left(  \uparrow\varphi^{k}\downarrow^{\otimes
k}\right)  _{i,n-k-i}:\left(  \uparrow\overline{A}\right)  ^{\otimes
n}\rightarrow\left(  \uparrow\overline{A}\right)  ^{\otimes n-k+1}.
\]
Let $\widetilde{B}A=T^{c}\left(  \uparrow\overline{A}\right)  $ and define a
map $d_{\widetilde{B}A}:\widetilde{B}A\rightarrow\widetilde{B}A$ of degree
$-1$ by
\begin{equation}
d_{\widetilde{B}A}=\sum\limits_{\substack{1\leq k\leq n \\0\leq i\leq n-k
}}\left(  \uparrow\varphi^{k}\downarrow^{\otimes k}\right)  _{i,n-k-i}.
\label{nine}%
\end{equation}
The identities $\left(  -1\right)  ^{\left[  n/2\right]  }\uparrow^{\otimes
n}\downarrow^{\otimes n}=1^{\otimes n}$ and $\left[  n/2\right]  +\left[
\left(  n+k\right)  /2\right]  \equiv nk+\left[  k/2\right]  $ (mod 2) imply
that
\begin{equation}
d_{\widetilde{B}A}=\sum\limits_{\substack{1\leq k\leq n \\0\leq i\leq n-k
}}\left(  -1\right)  ^{\left[  \left(  n-k\right)  /2\right]  +i\left(
k+1\right)  }\uparrow^{\otimes n-k+1}\varphi_{i,n-k-i}^{k}\downarrow^{\otimes
n}. \label{d-tilde-bar}%
\end{equation}

\begin{definition}
$\left(  A,\varphi^{n}\right)  _{n\geq1}$ is an \underline{$A_{\infty}%
$-algebra} if $d_{\widetilde{B}A}^{2}=0.$
\end{definition}

\begin{proposition}
\label{relations}For each $n\geq1,$ the operations $\left\{  \varphi
^{n}\right\}  $ on an $A_{\infty}$-algebra satisfy the following quadratic
relations:
\begin{equation}
\sum_{\substack{0\leq\ell\leq n-1 \\0\leq i\leq n-\ell-1}}\left(  -1\right)
^{\ell\left(  i+1\right)  }\varphi^{n-\ell}\varphi_{i,n-\ell-1-i}^{\ell+1}=0.
\label{A-infty-alg}%
\end{equation}

\end{proposition}

\begin{proof}
For $n\geq1,$
\begin{align*}
0  &  =\sum\limits_{\substack{1\leq k\leq n \\0\leq i\leq n-k}}\left(
-1\right)  ^{\left[  \left(  n-k\right)  /2\right]  +i\left(  k+1\right)
}\uparrow\varphi^{n-k+1}\downarrow^{\otimes n-k+1}\uparrow^{\otimes
n-k+1}\varphi_{i,n-k-i}^{k}\downarrow^{\otimes n}\\
&  =\sum\limits_{\substack{1\leq k\leq n \\0\leq i\leq n-k}}\left(  -1\right)
^{n-k+i\left(  k+1\right)  }\varphi^{n-k-1}\varphi_{i,n-k-i}^{k}\\
&  =-\left(  -1\right)  ^{n}\sum\limits_{\substack{0\leq\ell\leq n-1 \\0\leq
i\leq n-\ell-1}}\left(  -1\right)  ^{\ell\left(  i+1\right)  }\varphi^{n-\ell
}\varphi_{i,n-\ell-1-i}^{\ell+1}.
\end{align*}

\end{proof}

\noindent It is easy to prove that

\begin{proposition}
If $\left(  A,\varphi^{n}\right)  _{n\geq1}$ is an $A_{\infty}$-algebra, then
$\left(  \widetilde{B}A,d_{\widetilde{B}A}\right)  $ is a DGC.
\end{proposition}

\begin{definition}
Let $\left(  A,\varphi^{n}\right)  _{n\geq1}$ be an $A_{\infty}$%
-algebra.$\ $\underline{The tilde bar construction on $A$} is the DGC $\left(
\widetilde{B}A,d_{\widetilde{B}A}\right)  .$
\end{definition}

\begin{definition}
Let $A$ and $C$ be $A_{\infty}$-algebras.$\ $A chain map $f=f^{1}:A\rightarrow
C$ is a \underline{\textit{map of }$A_{\infty}$\textit{-algebras}}%
\textit{\ }if there is a sequence of maps $\{f^{k}\in Hom^{k-1}\left(
A^{\otimes k},C\right)  \}_{k\geq2}$ such that
\[
\widetilde{f}=\sum_{n\geq1}\left(  \sum_{k\geq1}\uparrow f^{k}\downarrow
^{\otimes k}\right)  ^{\otimes n}:\widetilde{B}A\rightarrow\widetilde{B}C
\]
is a DGC map.
\end{definition}

Dually, consider a sequence of operations $\{\psi^{k}\in Hom^{k-2}\left(
A,A^{\otimes k}\right)  \}_{k\geq1}.\ $For each $k$ and $n\geq1,$ linearly
extend each $\psi^{k}$ to $A^{\otimes n}$ via
\[
\sum\limits_{i=0}^{n-1}\psi_{i,n-1-i}^{k}:A^{\otimes n}\rightarrow A^{\otimes
n+k-1},
\]
and consider the induced map of degree $-1$ given by
\[
\sum\limits_{i=0}^{n-1}\left(  \downarrow^{\otimes k}\psi^{k}\uparrow\right)
_{i,n-1-i}:\left(  \downarrow\overline{A}\right)  ^{\otimes n}\rightarrow
\left(  \downarrow\overline{A}\right)  ^{\otimes n+k-1}.
\]
Let $\widetilde{\Omega}A=T^{a}\left(  \downarrow\overline{A}\right)  $ and
define a map $d_{\widetilde{\Omega}A}:\widetilde{\Omega}A\rightarrow
\widetilde{\Omega}A$ of degree $-1$ by
\[
d_{\widetilde{\Omega}A}=\sum\limits_{\substack{n,k\geq1 \\0\leq i\leq n-1
}}\left(  \downarrow^{\otimes k}\psi^{k}\uparrow\right)  _{i,n-1-i},
\]
which can be rewritten as
\begin{equation}
d_{\widetilde{\Omega}A}=\sum\limits_{\substack{n,k\geq1 \\0\leq i\leq n-1
}}\left(  -1\right)  ^{\left[  n/2\right]  +i\left(  k+1\right)  +k\left(
n+1\right)  }\downarrow^{\otimes n+k-1}\psi_{i,n-1-i}^{k}\uparrow^{\otimes n}.
\label{d-tilde-cobar}%
\end{equation}

\begin{definition}
$\left(  A,\psi^{n}\right)  _{n\geq1}$ is an \underline{$A_{\infty}%
$-coalgebra} if $d_{\widetilde{\Omega}A}^{2}=0.$
\end{definition}

\begin{proposition}
For each $n\geq1,$ the operations $\left\{  \psi^{k}\right\}  $ on an
$A_{\infty}$-coalgebra satisfy the following quadratic relations:
\begin{equation}
\sum_{\substack{0\leq\ell\leq n-1 \\0\leq i\leq n-\ell-1}} \left(  -1\right)
^{\ell\left(  n+i+1\right)  }\psi_{i,n-\ell-1-i}^{\ell+1}\psi^{n-\ell}=0.
\label{A-infty-coalg}%
\end{equation}

\end{proposition}

\begin{proof}
The proof is similar to the proof of Proposition \ref{relations} and is omitted.
\end{proof}

\noindent Again, it is easy to prove that

\begin{proposition}
If $\left(  A,\psi^{n}\right)  _{n\geq1}$ is an $A_{\infty}$-coalgebra , then
$\left(  \widetilde{\Omega}A,d_{\widetilde{\Omega}A}\right)  $ is a DGA.
\end{proposition}

\begin{definition}
Let $\left(  A,\psi^{n}\right)  _{n\geq1}$ be an $A_{\infty}$-coalgebra.$\ $%
The \underline{tilde cobar construction on} \underline{$A$}is the DGA $\left(
\widetilde{\Omega}A,d_{\widetilde{\Omega}A}\right)  .$
\end{definition}

\begin{definition}
Let $A$ and $B$ be $A_{\infty}$-coalgebras.$\ $A chain map $g=g^{1}%
:A\rightarrow B$ is a \underline{\textit{map of }$A_{\infty}$%
\textit{-coalgebras}}\textit{\ }if there is a sequence of maps $\{g^{k}\in
Hom^{k-1}\left(  A,B^{\otimes k}\right)  \}_{k\geq2}$ such that
\begin{equation}
\widetilde{g}=\sum_{n\geq1}\left(  \sum_{k\geq1}\downarrow^{\otimes k}%
g^{k}\uparrow\right)  ^{\otimes n}:\widetilde{\Omega}A\rightarrow
\widetilde{\Omega}B, \label{dga}%
\end{equation}
is a DGA map.
\end{definition}

The structure of an $A_{\infty}$-(co)algebra is encoded by the quadratic
relations among its operations (also called \textquotedblleft higher
homotopies\textquotedblright). Although the \textquotedblleft
direction,\textquotedblright\ i.e., sign, of these higher homotopies is
arbitrary, each choice of directions determines a set of signs in the
quadratic relations, the \textquotedblleft simplest\textquotedblright\ of
which appears on the algebra side when no changes of direction are made; see
(\ref{nine}) and (\ref{A-infty-alg}) above.$\ $Interestingly, the
\textquotedblleft simplest\textquotedblright\ set of signs appear on the
coalgebra side when $\psi^{n}$ is replaced by $\left(  -1\right)  ^{\left[
\left(  n-1\right)  /2\right]  }\psi^{n},$ $n\geq1,$ i.e., the direction of
every third and fourth homotopy is reversed.$\ $The choices one makes will
depend on the application; for us the appropriate choices are as in
(\ref{A-infty-alg}) and (\ref{A-infty-coalg}).

Let $\mathcal{A}_{\infty}=\oplus_{n\geq2}C_{\ast}\left(  K_{n}\right)  $ and
let $\left(  A,\varphi^{n}\right)  _{n\geq1}$ be an $A_{\infty}$-algebra with
quadratic relations as in (\ref{A-infty-alg}).$\ $For each $n\geq2,$ associate
$e^{n-2}\in C_{n-2}\left(  K_{n}\right)  $ with the operation $\varphi^{n}$
via
\begin{equation}
e^{n-2}\mapsto\left(  -1\right)  ^{n}\varphi^{n} \label{two}%
\end{equation}
and each codimension $1$ face $d_{(i,\ell)}\left(  e^{n-2}\right)  \in
C_{n-3}\left(  K_{n}\right)  $ with the quadratic composition
\begin{equation}
d_{(i,\ell)}\left(  e^{n-2}\right)  \mapsto\varphi^{n-\ell}\varphi
_{i,n-\ell-1-i}^{\ell+1}. \label{three}%
\end{equation}
Then (\ref{two}) and (\ref{three}) induce a chain map
\begin{equation}
\zeta_{A}:\mathcal{A}_{\infty}\longrightarrow\oplus_{n\geq2}Hom^{\ast}\left(
A^{\otimes n},A\right)  \label{alg}%
\end{equation}
representing the $A_{\infty}$-algebra structure on $A$.$\ $Dually, if $\left(
A,\psi^{n}\right)  _{n\geq1}$ is an $A_{\infty}$-coalge- bra with quadratic
relations as in (\ref{A-infty-coalg}), the associations
\[
e^{n-2}\mapsto\psi^{n}\text{$\ $and$\ $}d_{(i,\ell)}\left(  e^{n-2}\right)
\mapsto\psi_{i,n-\ell-1-i}^{\ell+1}\psi^{n-\ell}%
\]
induce a chain map
\begin{equation}
\xi_{A}:\mathcal{A}_{\infty}\longrightarrow\oplus_{n\geq2}Hom^{\ast}\left(
A,A^{\otimes n}\right)  \label{coalg}%
\end{equation}
representing the $A_{\infty}$-coalgebra structure on $A.\ $The definition of
the tensor product is now immediate:

\begin{definition}
\label{tensor product}The \underline{tensor product of $A_{\infty}$-algebras}
$\left(  A,\zeta_{A}\right)  $ and $\left(  B,\zeta_{B}\right)  $ is given by%
\[
\left(  A,\zeta_{A}\right)  \otimes\left(  B,\zeta_{B}\right)  =\left(
A\otimes B,\zeta_{A\otimes B}\right)  ,
\]
where $\zeta_{A\otimes B}$ is the sum of the compositions
\[%
\begin{array}
[c]{ccc}%
C_{\ast}(K_{n}) & \overset{{\LARGE \zeta}_{A\otimes B}}{\longrightarrow} &
Hom((A\otimes B)^{\otimes n},A\otimes B)\\
&  & \\
_{\strut\strut{\LARGE \Delta}_{K}}\text{$\ $}\downarrow\text{$\ $\ } &  &
\text{$\ \ \ \ \ \ $}\uparrow\text{$\ (\sigma_{2,n})^{\ast}\iota$}\\
&  & \\
C_{\ast}(K_{n})\otimes C_{\ast}(K_{n}) & \underset{{\LARGE \zeta}%
_{A}{\LARGE \otimes\zeta}_{B}}{\longrightarrow} & Hom(A^{\otimes n},A)\otimes
Hom\left(  B^{\otimes n},B\right)
\end{array}
\]
over all $n\geq2;$ the $A_{\infty}$-algebra operations $\Phi^{n}$ on $A\otimes
B$ are given by%
\[
\Phi^{n}=\left(  \sigma_{2,n}\right)  ^{\ast}\text{$\iota$}\left(  \zeta
_{A}\otimes\zeta_{B}\right)  \Delta_{K}\left(  e^{n-2}\right)  .
\]
Dually, the \underline{tensor product of $A_{\infty}$-coalgebras} $\left(
A,\xi_{A}\right)  $ and $\left(  B,\xi_{B}\right)  $ is given by%
\[
\left(  A,\xi_{A}\right)  \otimes\left(  B,\xi_{B}\right)  =\left(  A\otimes
B,\xi_{A\otimes B}\right)  ,
\]
where $\xi_{A\otimes B}$ is the sum of the compositions
\[%
\begin{array}
[c]{ccc}%
C_{\ast}(K_{n}) & \overset{{\LARGE \xi}_{A\otimes B}}{\longrightarrow} &
Hom(A\otimes B,(A\otimes B)^{\otimes n})\\
&  & \\
_{\strut\strut{\LARGE \Delta}_{K}}\text{$\ $}\downarrow\text{$\ $\ } &  &
\text{$\ \ \ \ \ \ $}\uparrow\text{$\ (\sigma_{n,2})_{\ast}\iota$}\\
&  & \\
C_{\ast}(K_{n})\otimes C_{\ast}(K_{n}) & \underset{{\LARGE \xi}_{A}%
{\LARGE \otimes\xi}_{B}}{\longrightarrow} & Hom(A,A^{\otimes n})\otimes
Hom\left(  B,B^{\otimes n}\right)
\end{array}
\]
over all $n\geq2;$ the $A_{\infty}$-coalgebra operations $\Psi^{n}$ on
$A\otimes B$ are given by
\[
\Psi^{n}=\left(  \sigma_{n,2}\right)  _{\ast}\text{$\iota$}\left(  \xi
_{A}\otimes\xi_{B}\right)  \Delta_{K}\left(  e^{n-2}\right)  .
\]

\end{definition}

\begin{example}
\label{iterate}If $\left(  A,\psi^{n}\right)  _{n\geq1}$ is an $A_{\infty}%
$-coalgebra,$\ $the following $A_{\infty}$ operations arise on $A\otimes A$:
\[%
\begin{array}
[c]{ll}%
\Psi^{1}= & \psi^{1}\otimes1+1\otimes\psi^{1}\\
& \\
\Psi^{2}= & \sigma_{2,2}\left(  \psi^{2}\otimes\psi^{2}\right) \\
& \\
\Psi^{3}= & \sigma_{3,2}\left(  \psi_{0}^{2}\psi_{0}^{2}\otimes\psi^{3}%
+\psi^{3}\otimes\psi_{1}^{2}\psi_{0}^{2}\right) \\
& \\
\Psi^{4}= & \sigma_{4,2}\left(  \psi_{0}^{2}\psi_{0}^{2}\psi_{0}^{2}%
\otimes\psi^{4}+\psi^{4}\otimes\psi_{2}^{2}\psi_{1}^{2}\psi_{0}^{2}+\psi
_{0}^{3}\psi_{0}^{2}\otimes\psi_{1}^{2}\psi_{0}^{3}\right. \\
&
\begin{array}
[c]{l}%
\hspace{0.3in}\left.  +\psi_{0}^{3}\psi_{0}^{2}\otimes\psi_{1}^{3}\psi_{0}%
^{2}+\psi_{1}^{2}\psi_{0}^{3}\otimes\psi_{1}^{3}\psi_{0}^{2}-\psi_{0}^{2}%
\psi_{0}^{3}\otimes\psi_{2}^{2}\psi_{0}^{3}\right)
\end{array}
\\
\text{$\ \ $}\vdots & \hspace{0.5in}\hspace{0.5in}\hspace{0.5in}\vdots
\end{array}
\vspace*{0.1in}%
\]

\end{example}

Note that the compositions in Definition \ref{tensor product} only use the
operations $\psi^{n}$ and not the quadratic relations (\ref{A-infty-coalg}).
Indeed, one can iterate an arbitrary family of operations $\left\{  \psi
^{n}\right\}  $ as in Example (\ref{iterate}) to produce iterated operations
$\Psi^{n}:A^{\otimes k}\rightarrow\left(  A^{\otimes k}\right)  ^{\otimes n}$
whether or not $\left(  A,\psi^{n}\right)  $ is an $A_{\infty}$-coalgebra. Of
course, the $\Psi^{n}$'s define an $A_{\infty}$-coalgebra structure on
$A^{\otimes k}$ whenever $d_{\widetilde{\Omega}\left(  A^{\otimes k}\right)
}^{2}=0,$ and we make extensive use of this fact in the sequel \cite{SU3}.
Finally, since $\Delta_{K}$ is homotopy coassociative (not strict), the tensor
product only iterates up to homotopy. In the sequel we always coassociate on
the extreme left.

\section{Appendix}

For completeness, we review the definitions of the functors given\ by
Kadeishvili and Saneblidze in \cite{KS1}, \cite{KS2} from the category of
1-reduced simplicial sets to the category of cubical sets and from the
category of 1-reduced cubical sets to the category of permutahedral sets.

\subsection{The cubical set functor $\mathbf{\Omega} X$}

Given a 1-reduced simplicial set $X=\{X_{n},\partial_{i},$\linebreak%
$s_{i}\}_{n\geq0},$ define the graded set $\mathbf{\Omega}X$ as follows: Let
$X^{c}$ be the graded set of formal expressions%
\[
X_{n+k}^{c}=\{\eta_{i_{k}}\cdots\eta_{i_{1}}\eta_{i_{0}}(x)\text{ }|\,x\in
X_{n}\}_{n\geq0;k\geq0},
\]
where $\eta_{i_{0}}=1,$ $i_{1}\leq\cdots\leq i_{k},\,1\leq i_{j}\leq
n+j-1,\,1\leq j\leq k,$ and let $\bar{X}^{c}=s^{-1}(X_{>0}^{c})$ be the
desuspension of $X^{c}.$ Let $\Omega^{\prime}X$ be the free graded monoid
generated by $\bar{X}^{c};$ denote elements of $\Omega^{\prime}X$ by $\bar
{x}_{1}\dotsm\bar{x}_{k},$ where $x_{j}\in X_{m_{j}+1},\ m_{j}\geq0.$ The
total degree $m=\left\vert \bar{x}_{1}\dotsm\bar{x}_{k}\right\vert =\sum
|\bar{x}_{j}|$ and we write $\bar{x}_{1}\dotsm\bar{x}_{k}\in(\Omega^{\prime
}X)_{m}.$ The product of two elements $\bar{x}_{1}\dotsm\bar{x}_{k} $ and
$\bar{y}_{1}\dotsm\bar{y}_{\ell}$ is given by concatenation $\bar{x}_{1}%
\dotsm\bar{x}_{k}\bar{y}_{1}\dotsm\bar{y}_{\ell};$ the only relation on
$\Omega^{\prime}X$ is strict associativity. Let $\mathbf{\Omega}X$ be the
graded monoid obtained from $\Omega^{\prime}X$ via%
\[
\mathbf{\Omega}X=\Omega^{\prime}X\diagup{\Huge \sim}\text{ ,}%
\]
where $\overline{\eta_{n}(x)}\sim\overline{s_{n}(x)}$ for $x\in X_{>0},$ and
$\bar{x}_{1}\cdots\overline{\eta_{m_{i}+1}(x_{i})}\cdot\bar{x}_{i+1}\cdots
\bar{x}_{k}\sim\bar{x}_{1}\cdots\bar{x}_{i}\cdot\overline{\eta_{1}(x_{i+1}%
)}\cdots\bar{x}_{k}$ for $x_{i}\in X_{m_{i}+1}^{c},$ $i<k.\ $Let $MX$ denote
the free monoid generated by $\bar{X}=s^{-1}(X_{>0});$ there is an inclusion
of graded modules $MX\subset\Omega^{\prime}X.$

Apparently $\Omega^{\prime}X$ canonically admits the structure of a cubical
set. Denote the components of Alexander-Whitney diagonal by
\[
\nu_{i}:X_{n}\rightarrow X_{i}\times X_{n-i},
\]
where$\ \nu_{i}(x)=\partial_{i+1}\dotsm\partial_{n}(x)\times\partial_{0}%
\dotsm\partial_{i-1}(x),\ 0\leq i\leq n,$ and let $x^{n}\in X_{n}$ denote an
$n$-simplex simplex. Then
\[
\nu_{i}(x^{n})=(x^{\prime})^{i}\times(x^{\prime\prime})^{n-i}\in X_{i}\times
X_{n-i}%
\]
for all $n>0.$ For $1\leq i\leq n-1,$ define face operators $d_{i}^{0}%
,d_{i}^{1}:(\mathbf{\Omega}X)_{n-1}\rightarrow(\mathbf{\Omega}X)_{n-2}$ on a
(monoidal) generator $\overline{x^{n}}\in\bar{X}_{n}\subset\bar{X}_{n}^{c}$ by%
\[
d_{i}^{0}\left(  \overline{x^{n}}\right)  =\overline{(x^{\prime})^{i}}%
\cdot\overline{(x^{\prime\prime})^{n-i}}\text{ \ and \ }d_{i}^{1}\left(
\overline{x^{n}}\right)  =\overline{\partial_{i}\left(  x^{n}\right)  },
\]
and extend to elements $\bar{x}_{1}\dotsm\bar{x}_{k}\in MX$ via
\[%
\begin{array}
[c]{l}%
d_{i}^{0}\left(  \overline{x}_{1}\cdots\overline{x}_{k}\right)  =\overline
{x}_{1}\cdots\overline{(x_{q}^{\prime})^{j_{q}}}\cdot\overline{(x_{q}%
^{\prime\prime})^{m_{q}-j_{q}+1}}\cdots\overline{x}_{k},\\
\\
d_{i}^{1}\left(  \overline{x}_{1}\cdots\overline{x}_{k}\right)  =\overline
{x}_{1}\cdots\overline{\partial_{j_{q}}\left(  x_{q}\right)  }\cdots
\overline{x}_{k},
\end{array}
\]
where $i=m_{(q-1)}+j_{q}\leq m_{(q)},\ 1\leq i\leq n-1,\ 1\leq q\leq k.$ Then
the defining identities for a cubical set involving $d_{i}^{0}$ and $d_{i}%
^{1}$ can easily be checked on $MX.$ In particular, the simplicial relations
between the $\partial_{i}$'s imply the cubical relations between $d_{i}^{1}%
$'s; the associativity relations between $\nu_{i}$'s imply the cubical
relations between $d_{i}^{0}$'s, and the commutativity relations between
$\partial_{i}$'s and $\nu_{j}$'s imply the cubical relations between
$d_{i}^{1}$'s and $d_{j}^{0}$'s. Next, define degeneracy operators $\eta
_{i}:(\mathbf{\Omega}X)_{n-1}\rightarrow(\mathbf{\Omega}X)_{n}$ on a
(monoidal) generator $\overline{x}\in(\bar{X}^{c})_{n-1}$ by
\[
\eta_{i}(\overline{x})=\overline{\eta_{i}(x)};
\]
and extend to elements $\bar{x}_{1}\dotsm\bar{x}_{k}\in\mathbf{\Omega}X$ via
\[%
\begin{array}
[c]{l}%
\eta_{i}\left(  \overline{x}_{1}\cdots\overline{x}_{k}\right)  =\overline
{x}_{1}\cdots\eta_{j_{q}}(\overline{x_{q}})\cdots\overline{x}_{k},\\
\\
\eta_{n}\left(  \overline{x}_{1}\cdots\overline{x}_{k}\right)  =\overline
{x}_{1}\cdots\overline{x}_{m_{k-1}}\cdot\eta_{m_{k}+1}\left(  \overline{x}%
_{k}\right)  ,
\end{array}
\]
where $i=m_{(q-1)}+\ j_{q}\leq m,\ 1\leq i\leq n-1,$ $1\leq q\leq k,$ and
extend face operators on degenerate elements inductively so that the defining
identities of a cubical set are satisfied. Then in particular, the following
identities hold for all $x^{n}\in X_{n}$:
\[%
\begin{array}
[c]{l}%
d_{1}^{0}\left(  \overline{x^{n}}\right)  =\overline{(x^{\prime})^{1}}%
\cdot\overline{(x^{\prime\prime})^{n-1}}=e\cdot\overline{(x^{\prime\prime
})^{n-1}}=\overline{(x^{\prime\prime})^{n-1}}=\overline{\partial_{0}\left(
x^{n}\right)  },\\
\\
d_{n-1}^{0}\left(  \overline{x^{n}}\right)  =\overline{(x^{\prime})^{n-1}%
}\cdot\overline{(x^{\prime\prime})^{1}}=\overline{(x^{\prime\prime})^{n-1}%
}\cdot e=\overline{(x^{\prime})^{n-1}}=\overline{\partial_{n}\left(
x^{n}\right)  },
\end{array}
\]
where $e\in(\mathbf{\Omega}X)_{0}$ denotes the unit. It is easy to see that
the cubical set $\{\mathbf{\Omega}X,d_{i}^{0},d_{i}^{1},$ $\eta_{i}\}$ depends
functorially on $X.$

\subsection{The permutahedral set functor ${\mathbf{\Omega}}Q$}

Let $Q=(Q_{n},d_{i}^{0},d_{i}^{1},\eta_{i})_{n\geq0}$ be a 1-reduced cubical
set. Recall that the diagonal
\[
\Delta:C_{\ast}(Q)\rightarrow C_{\ast}(Q)\otimes C_{\ast}(Q)
\]
on $C_{\ast}(Q)$ is defined on $a\in Q_{n}$ by
\[
\Delta(a)=\Sigma\ (-1)^{\epsilon}d_{B}^{0}(a)\otimes d_{A}^{1}(a),
\]
where $d_{B}^{0}=d_{j_{1}}^{0}\cdots d_{j_{q}}^{0},$\ $d_{A}^{1}=d_{i_{1}}%
^{1}\cdots d_{i_{p}}^{1},$ summation is over all shuffles $\left(  A;B\right)
=\left(  i_{1}<\cdots<i_{q};j_{1}<\cdots<j_{p}\right)  $ of $\underline{n}$
and $\epsilon$ is the sign of the shuffle. The primitive components of the
diagonal are given by the extreme cases $A=\varnothing$ and $B=\varnothing$.

Let $\bar{Q}=s^{-1}(Q_{>0})$ denote the desuspension of $Q,$ let
$\mathbf{\Omega}^{\prime\prime}Q$ be the free graded monoid generated by
$\bar{Q}$ with the unit $e\in\bar{Q}_{1}\subset\mathbf{\Omega}^{\prime\prime
}Q$ and let $\Upsilon$ be the set of formal expressions
\[
\Upsilon=\{\varrho_{M_{k}|N_{k}}((\cdots\varrho_{M_{2}|N_{2}}(\varrho
_{M_{1}|N_{1}}(\bar{a}_{1}\cdot\bar{a}_{2})\cdot\bar{a}_{3})\cdots)\cdot
\bar{a}_{k+1})|\,a_{i}\in Q_{r_{i}}\}_{r_{i}\geq1;k\geq2},
\]
$M_{i}|N_{i}\in{\mathcal{P}}_{r_{(i)},r_{i+1}}(r_{(i+1)})$ or $M_{i}|N_{i}%
\in{\mathcal{P}}_{r_{i+1},r_{(i)}}(r_{(i+1)}),$ $\,1\leq i\leq k. $ Note that
one or more of the $a_{i}$'s can be the unit $e$. Adjoin the elements of
$\Upsilon$ to $\mathbf{\Omega}^{\prime\prime}Q$ and obtain the graded monoid
$\mathbf{\Omega}^{\prime}Q\ $and let $\mathbf{\Omega}Q$ be the monoid
\[
\mathbf{\Omega}Q=\mathbf{\Omega}^{\prime}Q\diagup{\Huge \sim}\text{,}%
\]
where $\varrho_{M|N}(\bar{a}\cdot\bar{b})\sim\varrho_{N|M}(\bar{a}\cdot\bar
{b}),\,$ $\varrho_{j|\underline{n}\setminus j}(e\cdot\bar{a})\sim
\varrho_{\underline{n}\setminus j|j}(\bar{a}\cdot e)\sim\overline{\eta_{j}%
(a)},\,a,b\in Q_{>0},$ and $\bar{a}_{1}\cdots{\varrho_{\underline{r_{i}}%
|r_{i}+1}(\bar{a}_{i}\cdot e)}\cdot\bar{a}_{i+2}\cdots\bar{a}_{k+1}\sim\bar
{a}_{1}\cdots\bar{a}_{i}\cdot\varrho_{1|\underline{r_{i+2}+1}\setminus
1}(e\cdot\bar{a}_{i+2})\cdots\bar{a}_{k+1}$ for $a_{i}\in Q_{r_{i}}%
,a_{i+1}=e,\,1\leq i\leq k.$ Then $\mathbf{\Omega}Q$ is canonically a
multipermutahedral set in the following way: First, define the face operator
$d_{A|B}$ on a monoidal generator $\bar{a}\in\bar{Q}_{n}$ by
\[
d_{A|B}(\bar{a})=\overline{d_{B}^{0}(a)}\cdot\overline{d_{A}^{1}%
(a)},\ \ \ A|B\in\mathcal{P}_{\ast,\ast}{(n)}.
\]
Next, use the formulas in the definition of a singular multipermutahedral set
(\ref{multip}) to define $d_{A|B}$ and $\varrho_{M|N}$ on decomposables. In
particular, the following identities hold for $1\leq i\leq n$:%
\[
d_{i|\underline{n+1}\setminus i}\left(  \overline{a}\right)  =\overline
{d_{i}^{1}(a)}\text{ \ and \ }d_{\underline{n+1}\setminus i|i}\left(
\overline{a}\right)  =\overline{d_{i}^{0}(a)}.
\]
It is easy to see that $(\mathbf{\Omega}Q,\,d_{A|B},\,\varrho_{M|N})$ is a
multipermutahedral set that depends functorially on $Q.$

\begin{remark}
\label{lift}The fact that the definition of $\mathbf{\Omega}Q$ uses all
cubical degeneracies is justified geometrically by the fact that a degenerate
singular $n$-cube in the base of a path space fibration lifts to a singular
$(n-1)$-permutahedron in the fibre, which is degenerate with respect to
Milgram's projections \cite{Milgram} (c.f., the definition of the cubical set
$\mathbf{\Omega}X$ on a simplicial set $X$).
\end{remark}


\begin{thebibliography}{99}                                                                                               %


\bibitem {Adams}J. F. Adams, On the cobar construction, Proc. Nat. Acad. Sci.
(USA), 42 (1956), 409-412.

\bibitem {Baues1}H. J. Baues, The cobar construction as a Hopf algebra,
Invent. Math., 132 (1998), 467-489.

\bibitem {CM}G. Carlsson and R. J. Milgram, Stable homotopy and iterated loop
spaces, Handbook of Algebraic Topology (Edited by I. M. James), North-Holland
(1995), 505-583.

\bibitem {Coxeter}H.S.M. Coxeter, W.O.J. Moser, Generators and relations for
discrete groups, Springer-Verlag, 1972.

\bibitem {Gaberdiel-Zwiebach}Matthias R. Gaberdiel and Barton Zwiebach, Tensor
constructions of open string theories I: Foundations, Nucl.Phys. B505 (1997), 569-624.

\bibitem {Iwase}N. Iwase and M. Mimura, Higher homotopy associativity, Lecture
Notes in Math., 1370 (1986), 193-220.

\bibitem {Jones}D. W. Jones, A general theory of polyhedral sets and
corresponding T-complexes, Dissertationes Mathematicae, CCLXYI, Warszava (1988).

\bibitem {KS1}T. Kadeishvili and S. Saneblidze, A cubical model of a
fibration, preprint, AT/0210006.

\bibitem {KS2}--------------------, The twisted Cartesian model for the double
path space fibration, preprint, AT/0210224.

\bibitem {Kan}D. M. Kan, Abstract homotopy I, Proc. Nat. Acad. Sci. U.S.A., 41
(1955), 1092-1096.

\bibitem {Lada}T. Lada and M. Markl, Strongly homotopy Lie algebras,
\textit{Communications in Algebra} \textbf{23} (1995), 2147-2161.

\bibitem {Loday}J.-L. Loday and M. Ronco, Hopf algebra of the planar binary
trees, \textit{Adv. in Math. }\textbf{139, }No. 2 (1998), 293-309.

\bibitem {MacLane}S. Mac Lane, ``Homology,'' Springer-Verlag, Berlin/New York, 1967.

\bibitem {Milgram}R. J. Milgram, Iterated loop spaces, Ann. of Math., 84
(1966), 386-403.

\bibitem {SU1}S. Saneblidze and R. Umble, A diagonal on the associahedra,
preprint, math. AT/0011065.

\bibitem {SU3}-------------------, The biderivative and $A_{\infty}%
$-bialgebras, J. Homology, Homotopy and Appl., to appear; preprint, math.AT/0406270.

\bibitem {Smith}J. R. Smith, \textquotedblleft Iterating the cobar
construction,\textquotedblright\ Memiors of the Amer. Math. Soc. \textbf{109},
Number 524, Providence, RI, 1994.

\bibitem {stasheff}J. D. Stasheff, Homotopy associativity of $H$-spaces I, II,
\textit{Trans. Amer. Math. Soc.,} \textbf{108} (1963), 275-312.

\bibitem {tonks}A. Tonks, Relating the associahedron and the permutohedron, In
\textquotedblleft Operads: Proceedings of the Renaissance Conferences
(Hartford CT / Luminy Fr 1995)\textquotedblright\ Contemporary Mathematics,
202 (1997), pp.33-36$.$

\bibitem {Zigler}G. Ziegler, \textquotedblleft Lectures on
Polytopes,\textquotedblright\ GTM 152, Springer-Verlag, New York, 1995.
\end{thebibliography}
\end{document}